\documentclass{kms-j}
\usepackage{amsmath,amssymb} 
\usepackage{oldgerm} 
\usepackage{mathrsfs} 

\title{Generalized group algebras and generalized
measure algebras on non-discrete locally compact
abelian groups}

\begin{document}

\author{Jyunji Inoue and Sin-Ei Takahasi}
\date{} 

\maketitle

\begin{abstract} Let $G$ be a non-discrete LCA
group with the dual group $\Gamma$. We define
generalized group algebra, ${\mathcal L}^1(G)$, and
generalized measure algebra, ${\mathcal M}(G),$ on $G$ as generalizations of the group algebra 
$L^1(G)$ and the measure algebra $M(G)$, respectively. 
Generalized Fourier transforms of elements of
${\mathcal L}^1(G)$ and generalized Fourier-Stieltjes transforms of elements of ${\mathcal M}(G)$ are also defined as generalizations of the Fourier
transforms and the Fourier-Stieltjes transforms,
respectively. The image ${\mathcal A}(\Gamma)$ of
${\mathcal L}^1(G)$ by the generalized Fourier
transform becomes a function algebra on $\Gamma$
with norm inherited from ${\mathcal L}^1(G)$ through
this transform. It is shown that ${\mathcal A}(\Gamma)$ is a natural Banach function algebra
on\, $\Gamma$\,which is BSE and BED. It turns out
that ${\mathcal L}^1(G)$ contains all Rajchman
measures. 
Segal algebras in ${\mathcal L}^1(G)$ are defined and investigated. It is shown that there exists the
smallest isometrically character invariant Segal
algebra in ${\mathcal L}^1(G)$, which (eventually)
coincides with the smallest isometrically
character invariant Segal algebra in $L^1(G)$, the
Feichtinger algebra of $G$. 
A notion of locally bounded elements of 
${\mathcal M}(G)$ is introduced and investigated. It is shown that for each locally bounded element 
$\mu$ of ${\mathcal M}(G)$ there corresponds a unique Radon measure $\iota \mu$ on $G$ which characterizes $\mu$. We investigate  the multiplier algebra
$\mathbb{M}({\mathcal L}^1(G))$ of ${\mathcal L}^1(G)$, and obtain a result that there is a natural
continuous isomorphism from $\mathbb{M}
({\mathcal L}^1(G))$ into $A(G)^*$, the algebra of
pseudomeasures on $G$. When $G$ is compact, this
map becomes surjective and isometric.
\end{abstract}


\renewcommand{\thefootnote}{\fnsymbol{footnote}}
\footnote[0]{\!\!\!\!\!\!\!\!2020 {\it Mathematics
subject classification}: Primary 43A20; Secondary
42A16, 43A25 }



{\bf \S1.\ Introduction and preliminaries} \ In
this paper, $G$ denotes a non-discrete LCA group
with the dual group $\Gamma$. $(L^1(G), \|\cdot
\|_1)$ and $(M(G), \|\cdot \|)$ denote the usual
group algebra and measure algebra on $G$ with
convolution as multiplication, respectively. Haar
measures $m_G$ and $m_\Gamma$ of $G$ and $\Gamma$,
respectively, are normalized so that the inversion
theorem is valid. Infinitesimal increments $dx$
and $d\gamma$ of $m_G$ and $m_\Gamma$,
respectively, are also used. Members of $L^1(G)$
are identified with the absolutely continuous
measures with respect to $m_G$. $C_b(G)$ and
$C_0(G)$ denote the space of all bounded
continuous functions on $G$, and the space of all
continuous functions on $G$ vanishing at infinity,
respectively. $C_c(G)$ denotes the subspace of
$C_0(G)$ consisting of all functions with compact
supports.

For a character $\gamma\in \Gamma$, we write
$(x,\gamma)$ instead of $\gamma(x)$. Hence for $
f\in L^1(G)$ and $\mu\in M(G)$, the Fourier
transform\, (resp.\,the inverse Fourier transform)
and the Fourier-Stieltjes transform\,(resp.\,the
inverse Fourier-Stieltjes transform) of $f$ and
$\mu$ are given by
\begin{eqnarray*}
 \hat{f}(\gamma)&=&\int_G (-x,
\gamma)f(x)dx\ \ \biggr({\rm resp.}\ \check{f}
(\gamma)=\int_G (x, \gamma)f(x)dx\biggr), \\
\hat{\mu}(\gamma)&=&\int_G (-x, \gamma)d\mu(x)
\ \biggr({\rm resp.}\ \check{\mu}(\gamma)=\int_G (x, \gamma)d\mu(x)\biggr)\ \ (\gamma\in \Gamma).
\end{eqnarray*}

$M_0(G)$ denotes the subalgebra of all measures
$\mu\in M(G)$ whose Fourier-Stieltjes transforms
$\hat{\mu}$ vanish at infinity. Measures in
$M_0(G)$ are called Rajchman measures ([12]).

$(A(\Gamma), \|\cdot \|_A)$ and $(B(\Gamma),
\|\cdot\|_B)$ denote the Fourier algebra and the
Fourier-Stieltjes algebra of $\Gamma$,
respectively: $A(\Gamma)=\{\hat{f}: f\in L^1(G)
\}$ and $B(\Gamma)=\{\hat{\mu}: \mu\in M(G)\}$,
and norms $\|\cdot \|_A$ and $\|\cdot \|_B$ are
inherited from $L^1(G)$ and from $M(G)$ through
the Fourier and the Fourier-Stieltjes transform,
respectively.
  
We denote by ${\mathcal K}(G)\,({\rm resp.}\, 
{\mathcal K}(\Gamma))$ the set of all compact subsets of
$G\,({\rm resp.}\, \Gamma)$ which have dense
nonempty interiors. ${\mathcal K}(G)\,({\rm resp.}\,
{\mathcal K}(\Gamma))$ forms a directed set with
respect to the inclusion order. For $K\in 
{\mathcal K}(G)$, let $C_c^K(G)=\{h\in C_c(G): {\rm
supp}(h)\subset K\}$. $\textfrak{M}(G)$ denotes
the space of Radon measures on $G$, 
 namely, the space of
continuous linear functionals of $C_c(G)$ which is
the inductive limit of the spaces $\{(C_c^K(G),
\|\cdot\|_\infty): K\in {\mathcal K}(G)\}$.

A function $\sigma$ on $\Gamma$ is called a local
$A(\Gamma)$-function if it satisfies $\sigma
\varphi\in A(\Gamma)\ (\varphi\in A(\Gamma)_c)$,
where $A(\Gamma)_c:=\{\varphi\in A(\Gamma): {\rm
supp}(\varphi)\ {\rm is\ compact}\}$. The space of
all local $A(\Gamma)$-functions is denoted by
$A(\Gamma)_{loc}$ ([8, \S 5]).

\vspace{0.2cm}

Now, let $(\mathbf{B}, \|\cdot \|_{\mathbf B})$ be
a commutative semisimple Banach algebra with the
Gelfand space $X$, and the Gelfand transform $b\,
\mapsto \hat{b}, {\bf B}\rightarrow C_0(X)$. If we
put $\hat{\bf B}=\{\hat{b}: b\in {\mathbf B}\}, \
\|\hat{b}\|_{\hat{\mathbf B}}=\|b\|_{\mathbf B}\ \
(b\in {\mathbf B})$, then $(\hat{\mathbf B},
\|\cdot \|_{\hat{\mathbf B}})$ becomes a natural
Banach function algebra which is isometrically
isomorphic with ${\bf B}$ through the Gelfand
transform. ${\mathbf B}_c$ is the space of
elements of ${\mathbf B}$ whose Gelfand transforms have compact supports.

A multiplier of ${\mathbf B}$ is a linear operator
$T$ of $\mathbf{B}$ to itself which satisfies
$T(ab)=aT(b)\ \ (a,b\in {\mathbf B})$. The set of
all multipliers of ${\mathbf B}$ forms a
commutative Banach algebra with respect to the
operator norm, which will be denoted by
$\mathbb{M}({\mathbf B})$. It is well-known that
for each $T\in \mathbb{M}(\mathbf B)$, there
exists a unique $\hat{T}\in C_b(X)$ such that
$\widehat{Tb}=\hat{T}\hat{b}\ \ (b\in
\mathbf{B})$. From this, it follows that
$\mathbb{M}(\hat{\mathbf{B}})=\{\varphi\in C_b(X):
\varphi \hat{b}\in \hat{\mathbf{B}}\ \ (b\in
\mathbf{B})\}$.

When ${\mathbf B}=L^1(G)$, every multiplier of
$L^1(G)$ is given by the convolution of an element
$\mu\in M(G)$, and its multiplier norm coincides
with $\|\mu\|$, so $\mathbb{M}(L^1(G))$ is
identified with $M(G)$, that is,
$\mathbb{M}(L^1(G))=M(G)$ and
$\mathbb{M}(A(\Gamma))=B(\Gamma)$. For notations
and fundamental facts about multipliers, we refer
to [11]. \vspace{0.3cm}

Let $Y$ be a subset of $X$. We denote by ${\rm
span}(Y)$ the linear span of $Y$ in the dual space
$\mathbf{B}^*$ of $ \mathbf{B}$. Therefore each
$p\in {\rm span}(Y)$ can be written uniquely in
the form $p=\sum_{y\in Y}\hat{p}(y)y$ with
$\hat{p}(y)\neq 0$ at most finitely many $y\in Y$.

Under the above notation, the spaces $C_{BSE({\bf
B})}(X)$ and $C^0_{BSE({\bf B})}(X)$ are defined
as follows: 
\begin{eqnarray*} C_{BSE({\bf
B})}(X)&=&\{\sigma\in C_b(X): \|\sigma\|_{BSE({\bf
B})}<\infty\}; \\ C^0_{BSE({\bf
B})}(X)&=&\{\sigma\in C_{BSE({\bf B})}(X), \
\lim_{K\in {\mathcal K}(X)}\|\sigma\|_{BSE({\bf B}),
K}=0\}, 
\end{eqnarray*}
where
 \begin{eqnarray*}
\|\sigma\|_{BSE({\bf
B})}&=&\sup\biggl\{\biggl|\sum_{x\in
X}\hat{p}(x)\sigma(x)\biggr|: p\in {\rm span}(X),
\|p\|_{{\bf B}^*}\le 1\biggr\},\ {\rm and} \\
\|\sigma\|_{BSE({\bf B}), K}&=&\sup\biggl\{
\biggl|\sum_{x\in X\setminus
K}\hat{p}(x)\sigma(x)\biggr|: p\in {\rm
span}(X\setminus K), \|p\|_{{\bf B}^*}\le
1\biggr\}.
\end{eqnarray*}

We say ${\bf B}$ is BSE if ${\mathbb M}(\hat{\bf
B})=C_{BSE({\bf B})}(X)$, and is BED if $\hat{\bf
B}=C^0_{BSE({\bf B})}(X)$ ([7]).
 
\vspace{0.5cm}

Outline of the paper: 

In $\S 2$, for a given relatively compact open
neighborhood $V$ of $0\in \Gamma$, definitions of
generalized measure algebras ${\mathcal M}_{\bar{V}}(G)$ of $G$, and generalized group
algebras ${\mathcal L}^1_{\bar{V}}(G)$ of $G$ are
given. Definitions of generalized
Fourier-Stieltjes transform $\mu\mapsto \hat{\mu},
{\mathcal M}_{\bar{V}}(G)\rightarrow C_b(\Gamma)$ and generalized Fourier transform $f\mapsto \hat{f},
\,{\mathcal L}^1_{\bar{V}}(G) \rightarrow C_0(\Gamma)$ are also given.

In $\S 3$, generalized Fourier-Stieltjes algebra
${\mathcal B}_{\bar{V}}(\Gamma)$ of $\Gamma$, and
generalized Fourier algebra ${\mathcal A}_{\bar{V}}(\Gamma)$ of $\Gamma$ are defined. The
generalized Fourier-Stieltjes transform is an
isometric algebra homomorphism of 
${\mathcal M}_{\bar{V}}(G)$ onto 
${\mathcal B}_{\bar{V}}(\Gamma)$.

The fundamental properties of 
${\mathcal A}_{\bar{V}}(\Gamma)$ are proved in Theorem 3.1.
The generalized Fourier transform is an isometric
algebra homomorphism of ${\mathcal L}^1_{\bar{V}}(G)$
onto ${\mathcal A}_{\bar{V}}(\Gamma)$. A remarkable
fact that $M_0(G)$ is contained in 
${\mathcal L}^1_{\bar{V}}(G)$
 is proved in Theorem 3.4.

In $\S 4$, Segal algebras in 
${\mathcal L}^1_{\bar{V}}(G)$ are defined and investigated.
The smallest isometrically character invariant
Segal algebra in ${\mathcal L}^1_{\bar{V}}(G)$ is
constructed. In the final part of this section, it
is shown that ${\mathcal L}^1_{\bar{V}}(G)$ and ${\mathcal M}_{\bar{V}}(G)$ are determined without depending
on the choice of $V$. So, we can use the notations
${\mathcal L}^1(G)$ and ${\mathcal M}(G)$ instead of the
notations ${\mathcal L}^1_{\bar{V}}(G)$ and 
 ${\mathcal M}_{\bar{V}}(G)$.

$\S 5$ is devoted to prove that ${\mathcal B}(\Gamma)\subsetneqq 
\mathbb{M}({\mathcal
A}(\Gamma))$(the multiplier algebra of ${\mathcal
A}(\Gamma))$.

In $\S 6$, some Segal algebras in 
${\mathcal A}(\Gamma)$ (other than treated in $\S 4$) are
investigated.

In $\S 7$, the notion of locally boundedness for
elements of ${\mathcal M}(G)$ is introduced. We give a theorem (Theorem 7.5) which characterizes locally bounded elements of ${\mathcal M}(G)$ using
transformable Radon measures in [1].
 
In $\S 8$ we study the multiplier algebra
$\mathbb{M}({\mathcal L}^1(G))$ of ${\mathcal L}^1(G)$. It
is shown that there exists a natural continuous
algebra isomorphism of $\mathbb{M}({\mathcal L}^1(G))$ into $A(G)^*$, the algebra of pseudomeasures on $G$. When $G$ is compact, 
this isomorphism becomes surjective and isometric.

In the last part of this section, ''locally
boundedness'' for elements of 
$\mathbb{M}({\mathcal
L}^1(\mathbf{R}^d))$ are discussed. We give a 
definition of locally bounded elements
of $\mathbb{M}({\mathcal L}^1(\mathbf{R}^d))$ and its characterizing theorem\,(Theorem 8.6) which  
is similar to that for locally bounded elements of 
${\mathcal M}(G)$ in $\S 7$.

\vspace{0.5cm}
 
{\bf $\S2.$\,Definitions:\,Generalized measure
algebra and generalized group algebra of
$G$\ \ \rm

In the rest of the paper, the symbol $V$ is used
to denote an open neighborhood of $0\in \Gamma$
with compact closure $\bar{V}$.

\vspace{0.3cm}

{\bf Definition 2.1.}
\ Let $\phi\in C_b(\Gamma)$
and $E\subset \Gamma$. We denote by ${\rm
span}(E)$ the linear span of the elements of $E$
in the dual space $L^{\infty}(G)$ of $L^1(G)$.
Therefore, each element $p\in {\rm span}(E)$ can
be expressed uniquely in the form
$p=\sum_{\gamma\in E}\hat{p}(\gamma)\gamma$ with
$\hat{p}(\gamma)\neq 0$ for at most a finite
number of $\gamma$ in $E$. As usual,
$\|p\|_{\infty}:=\sup\{|p(x)|: x\in G \}$. We use
the notations

$\displaystyle\|\phi\|^E=\sup_{\gamma\in \Gamma}
\|\phi\|^{E, \gamma}$, where $\|\phi\|^{E,
\gamma}=\sup\bigl\{\bigr|\sum_{\gamma'\in
E+\gamma} \hat{p}(\gamma')\phi(\gamma')\bigr|:
p\in {\rm span}(E+\gamma), \|p\|_\infty\le 1
\bigr\}$

\vspace{0.5cm}

\it

{\bf Lemma 2.1.}\
 {\rm (cf. [18, Lemma 1])}
Suppose $\gamma\in \Gamma$ and $E\subset \Gamma$, and let $\psi_1, \psi_2\in C_b(\Gamma)$. Then we have 
\[ {\rm (i)}\ \|\psi_1\psi_2\|^{E, \gamma}\le
\|\psi_1\|^{E, \gamma}\|\psi_2\|^{E, \gamma}.\ \
{\rm (ii)}\ \|\psi_1\psi_2\|^E\le
\|\psi_1\|^E\|\psi_2\|^E. 
\]

\rm 

Proof.\ (i)\ We can suppose $\|\psi_j\|^{E,\gamma} <\infty,\, j=1,2.$ Let $p\in {\rm span}(E+\gamma)$ be given. For any $\varepsilon>0$, we can choose $h\in L^1(G)$ such that $\|h\|_1\le 1$ and 
\[
\biggl\|\sum_{\gamma'\in E+\gamma}
\hat{p}(\gamma')
\psi_1(\gamma')\gamma'\biggr\|_{\infty} \le
\biggl|\sum_{\gamma'\in
E+\gamma}\hat{p}(\gamma')\psi_1(\gamma')
\widehat{h}(\gamma')\biggr|+\varepsilon. \] Choose
$k\in L^1(G)$ such that $\|k\|_1\le 1$ and 
\[
\biggl\|\sum_{\gamma'\in E+\gamma}\hat{p}(\gamma')\hat{h}(\gamma')\gamma'\biggl\|_{\infty}\le
\biggl|\sum_{\gamma'\in E+\gamma}\hat{p}(\gamma')\hat{h}(\gamma')\hat{k}(\gamma')\biggr|+\varepsilon.
\]
 Then we have 
\begin{eqnarray*}
\biggl|\sum_{\gamma'\in
E+\gamma}\hat{p}(\gamma')\psi_1(\gamma')
\psi_2(\gamma')\biggr| &\le& \|\psi_2\|^{E,
\gamma} \biggl\|\sum_{\gamma'\in E+\gamma}
\hat{p}(\gamma')
\psi_1(\gamma')\gamma'\biggr\|_{\infty} \\ &\le&
\|\psi_2\|^{E,\gamma}\biggl(\biggl|
\sum_{\gamma'\in E+\gamma}\hat{p}(\gamma')
\psi_1(\gamma')\hat{h}(\gamma')\biggl|
+\varepsilon\biggr) \\ &\le& \|\psi_2\|^{E,
\gamma} \biggl(\|\psi_1\|^{E, \gamma}
\biggl\|\sum_{\gamma'\in E+\gamma}
\hat{p}(\gamma')\hat{h}(\gamma')\gamma'
\biggr\|_{\infty}+\varepsilon\biggr) \\ &\le&
\|\psi_2\|^{E, \gamma} \biggl(\|\psi_1\|^{E,
\gamma} \biggl(\biggl|\sum_{\gamma'\in
E+\gamma}\hat{p}(\gamma')
\hat{h}(\gamma')\hat{k}(\gamma')
\biggr|+\varepsilon\biggr)+\varepsilon\biggr) \\
&\le& \|\psi_2\|^{E, \gamma} \biggl(\|\psi_1\|^{E,
\gamma} \biggl(\biggl\|\sum_{\gamma'\in
E+\gamma}\hat{p}(\gamma')
\gamma'\biggr\|_{\infty}\|h*k\|_1+\varepsilon\biggr)
+\varepsilon\biggr). 
\end{eqnarray*}
 Since $\|h*k\|_1\le 1$ and $\varepsilon>0$ is arbitrary, we get the desired result.

(ii)\ An easy consequence of (i).\ \ \ $\Box$
\vspace{0.5cm}

Now, for $\mu\in M(G)$, we put $
\|\mu\|_{\bar{V}}:=\sup_{\gamma\in
\Gamma}\|\hat{\mu}\|^{\bar{V}, \gamma}$. It
follows easily that $\|\mu\|_{\bar{V}}\le \|\mu\|
\ (\mu\in M(G))$, and $(M(G), \|\cdot
\|_{\bar{V}})$ is a normed linear space. Moreover,
by Lemma 2.1, $(M(G), \|\cdot \|_{\bar{V}})$
becomes a commutative normed algebra. We denote by $\mathcal{M}_{\bar{V}}(G)$ the completion of this
normed algebra. We use the same symbol $*$ and
$\|\cdot \|_{\bar{V}}$ as in $M(G)$ to represent
the multiplication and norm of 
${\mathcal M}_{\ ar{V}}(G)$, respectively. Then
\[
 {\mathcal M}_{\bar{V}}(G)=\bigl\{\mu=\lim_n \mu_n: \{\mu_n\}_n \ {\rm is\ a\ Cauchy\ sequence\ in\ }\ (M(G), \|\cdot \|_{\bar{V}})\bigr\}, 
\]
with $\|\mu\|_{\bar{V}}=\lim_n\|\mu_n\|_{\bar{V}}$, and
\begin{eqnarray*}
\mu*\nu&=&\lim_{n\rightarrow\infty}\mu_n*\nu_n,
\|\mu*\nu\|_{\bar{V}}=\lim_{n\rightarrow\infty}
\|\mu_n*\nu_n\|_{\bar{V}}\le
\|\mu\|_{\bar{V}}\|\nu\|_{\bar{V}} \nonumber \\ &&
\ \ \ (\mu=\lim_n \mu_n, \nu=\lim_n \nu_n\in 
{\mathcal
M}_{\bar{V}}(G)). 
\end{eqnarray*}

{\bf Definition 2.2.}\ 
(Generalized measure algebra) 
We call the commutative Banach algebra
$({\mathcal M}_{\bar{V}}(G), \|\cdot \|_{\bar{V}})$
defined in the above, generalized measure algebra
of $G$.
\vspace{0.5cm}

{\bf Definition 2.3.}\ 
(Generalized group
algebra)\ We define ${\mathcal L}^1_{\bar{V}}(G)$ as
the closure of $L^1(G)$ in $({\mathcal
M}_{\bar{V}}(G), \|\cdot \|_{\bar{V}})$, that is,
\[{\mathcal L}^1_{\bar{V}}(G)=\bigl\{\lim_n f_n:
\{f_n\}_n\ {\rm is\ a\ Cauchy\ sequence\ in}\
(L^1(G), \|\cdot \|_{\bar{V}})\bigr\}. \] We call
(${\mathcal L}^1_{\bar{V}}(G), \|\cdot \|_{\bar{V}})$
generalized\ group algebra of $G$.
\vspace{0.5cm}

We will show in Remark 4.2 below that all 
${\mathcal
M}_{\bar{V}}(G)$ and all ${\mathcal L}^1_{\bar{V}}(G)$
are determined without depending on the choice of
$V$, but norms $\|\cdot \|_{\bar{V}}$ are
equivalent each other.
\vspace{0.5cm} 

{\bf Definition 2.4.}\
(Generalized Fourier-Stieltjes transform, generalized Fourier transform)\ For $\mu=\lim_n \mu_n\in
 {\mathcal
M}_{\bar{V}}(G),\ \mu_n\in M(G), n=1,2,3,...$,
define \begin{equation} \hat{\mu}(\gamma):=\lim_n
\widehat{\mu_n}(\gamma)\ \ (\gamma\in \Gamma).
\end{equation}

For $f=\lim_n f_n\in {\mathcal L}^1_{\bar{V}}(G),\
f_n\in L^1(G), n=1,2,3,...$, define
\begin{equation}
\hat{f}(\gamma):=\lim_n
\widehat{f_n}(\gamma)\ \ (\gamma\in \Gamma).
\end{equation}

It is clear that the right hand sides of (1) and
(2) have limits which are determined independently
of the choices of $\{\mu_n\}_{n=1}^{\infty}$ and
$\{f_n\}_{n=1}^\infty$ used to express $\mu$ and
$f$, respectively. We call $\mu\mapsto \hat{\mu},
{\mathcal M}_{\bar{V}}(G)\rightarrow C_b(\Gamma)$
generalized Fourier-Stieltjes transform, and
$f\mapsto \hat{f},\ {\mathcal
L}_{\bar{V}}^1(G)\rightarrow C_0(\Gamma)$
generalized Fourier transform,
respectively.
 \vspace{0.5cm}

Note that ${\mathcal M}_{\bar{V}}(G)\supset M(G)$ and, the algebraic operations and the generalized
Fourier-Stieltjes transform of ${\mathcal
M}_{\bar{V}}(G)$ are extensions of the algebraic
operations and the Fourier-Stieltjes transform of
$M(G)$, respectively.

Similarly, ${\mathcal L}^1_{\bar{V}}(G)\supset L^1(G)$ and, the algebraic operations and the generalized
Fourier transform of ${\mathcal L}^1_{\bar{V}}(G)$ are extensions of the algebraic operations and the
Fourier transform of $L^1(G)$, respectively.
\vspace{0.5cm}

{\bf \S 3.\ Generalized measure algebra and generalized group algebra of $G$; generalized 
Fourier algebra  and generalized Fourier-Stieltjes algebra of $\Gamma$}

This section constitutes a core of the present paper. 
\vspace{0.5cm}

{\bf Definition 3.1.}\ 
We define the space
\begin{eqnarray*}
&&C_{BSE}^{\bar{V}}(\Gamma):=\biggl\{\phi\in
C_b(\Gamma): \|\phi\|^{\bar{V}}:=\sup_{\gamma\in
\Gamma} \|\phi\|^{\bar{V},\gamma}<\infty\biggr\},
\\ & & \hspace{1.5cm}{\rm where\ }
\|\phi\|^{\bar{V},\gamma}:=\sup\biggl\{\biggl
|\sum_{\gamma'\in
\bar{V}+\gamma}\hat{p}(\gamma')\phi(\gamma')\biggr
| : p\in {\rm span}(\bar{V}+\gamma),
\|p\|_{\infty}\le 1\biggr\} . 
\end{eqnarray*}

By\ [8,\,Theorems 7.4 and 10.3],
$C_{BSE}^{\bar{V}}(\Gamma)\subset A(\Gamma)_{loc}$ and $(C_{BSE}^{\bar{V}}(\Gamma), \|\cdot
\|^{\bar{V}})$ forms a Banach function algebra on
$\Gamma$. Note that the generalized
Fourier-Stieltjes transform is an isometric
isomorphism of ${\mathcal M}_{\bar{V}}(G)$ into
$C^{\bar{V}}_{BSE}(\Gamma)$.

We define 
\begin{eqnarray*}
 {\mathcal B}_{\bar{V}}(\Gamma) &:=& 
 \bigl\{\phi\in C_{BSE}^{\bar{V}}(\Gamma): \phi=\hat{\mu} \ {\rm
for\ some\ } \mu\in {\mathcal M}_{\bar{V}}(G)\bigr\},
\ {\rm with\ norm}\ \|\cdot \|^{\bar{V}}, 
\\ {\mathcal A}_{\bar{V}}(\Gamma) &:=&\bigl\{\psi\in C_{BSE}^{\bar{V}}(\Gamma): 
\lim_{\Omega\in {\mathcal
K}(\Gamma)} \sup_{\gamma \in \Gamma\setminus
\Omega} \|\psi\|^{\bar{V}, \gamma}=0 \bigr\}.
\end{eqnarray*}

\vspace{0.5cm} 

\it

{\bf Theorem 3.1.}\ The space
$(\mathcal{A}_{\bar{V}}(\Gamma), \|\cdot
\|^{\bar{V}})$ is a Banach function algebra on
$\Gamma$ which satisfies the following:

{\rm (i)}\ $\mathcal{A}_{\bar{V}}(\Gamma)$\ has a
bounded approximate identity composed of elements
in ${\mathcal A}_{\bar{V}}(\Gamma)_c$, and the
generalized Fourier transform, denoted by 
${\mathcal
F}$, 
\[ {\mathcal F}: f\mapsto \hat{f},\ {\mathcal
L}^1_{\bar{V}}(G)\rightarrow {\mathcal
A}_{\bar{V}}(\Gamma) \] is a surjective isometric
algebra homomorphism;

{\rm (ii)}\ $\mathcal{A}_{\bar{V}}(\Gamma)$ is a
natural Banach function algebra, hence is regular;

{\rm (iii)}\ $\mathcal{A}_{\bar{V}}(\Gamma)$\ is
BSE and BED, that is, 
\begin{equation}
C_{BSE({\mathcal A}_{\bar{V}}(\Gamma))}(\Gamma)=
\mathbb{M}({\mathcal A}_{\bar{V}}(\Gamma)),\ {\rm
and}\ C_{BSE({\mathcal
A}_{\bar{V}}(\Gamma))}^0(\Gamma) ={\mathcal
A}_{\bar{V}}(\Gamma); 
\end{equation}

{\rm (iv)} \
$C_{BSE}^{\bar{V}}(\Gamma)=\mathbb{M}({\mathcal
A}_{\bar{V}}(\Gamma))$, and
$\|\phi\|^{\bar{V}}=\|\phi\|_{\mathbb{M}({\mathcal
A}_{\bar{V}}(\Gamma))}\ \ (\phi\in
C_{BSE}^{\bar{V}}(\Gamma)).$
\rm

Proof.\ (i)\ For each $\Omega\in\mathcal{K}(\Gamma)$, choose $u_\Omega\in A(\Gamma)_c$ such that $u_\Omega(\gamma)=1\ \ (\gamma\in \Omega)$ with
$\|u_\Omega\|^{\bar{V}}\le \|u_\Omega\|_A\le 2$
([17, Theorem 2.6.8.]).

That ${\mathcal F}$ is an isometric algebra
homomorphism into ${\mathcal A}_{\bar{V}}(\Gamma)$ is apparent. To see that ${\mathcal F}$ is surjective, let $\psi\in {\mathcal A}_{\bar{V}}(\Gamma)$ and $\varepsilon>0$ be given arbitrarily. Since $\displaystyle\lim_{\Omega\in {\mathcal K}(\Gamma)} \sup_{\gamma\in \Gamma\setminus
\Omega}\|\psi\|^{\bar{V},\gamma}=0$, there is
$\Omega_0\in {\mathcal K}(\Gamma)$ such that
$\|\psi\|^{\bar{V},\gamma}\le \varepsilon/3$ for
all $\gamma\in \Gamma\setminus \Omega_0$. Put
$\Omega_1=\Omega_0+\bar{V}$. For $\Omega\in {\mathcal K}(\Gamma)$ with $\Omega_1\subset \Omega$, we have from Lemma 2.1 
\begin{eqnarray*} 
\|\psi-\psi u_\Omega\|^{\bar{V}} &=& \sup_{\gamma\in
\Gamma}\|\psi(1-u_\Omega)\|^{\bar{V}, \gamma}\le
\sup_{\gamma\in \Gamma}\|\psi\|^{\bar{V}, \gamma}
\|1-u_\Omega\|^{\bar{V},\gamma} \\
&\le&\sup_{\gamma\in \Gamma\setminus
\Omega_0}\|\psi\|^{\bar{V}, \gamma}
\|1-u_\Omega\|^{\bar{V}, \gamma}(\because {\rm if
}\ \gamma\in \Omega_0,
\|1-u_\Omega\|^{\bar{V},\gamma}=0) \\
&\le&\frac{\varepsilon}{3}\cdot 3=\varepsilon\
(\because \|1\|^{\bar{V},\gamma}=1,
\|u_\Omega\|^{\bar{V},\gamma}\le\|u_\Omega\|_A\le 2), 
\end{eqnarray*}
 and hence $\{u_\Omega:
\Omega\in {\mathcal K}(\Gamma)\}$ is a bounded
approximate identity for ${\mathcal
A}_{\bar{V}}(\Gamma)$ composed of elements in
${\mathcal A}_{\bar{V}}(\Gamma)_c$, and since $\psi\in
A(\Gamma)_{loc}$, we have $u_\Omega\psi\in
A(\Gamma)$. This implies ${\mathcal
A}_{\bar{V}}(\Gamma) \subseteq \{\hat{f}: f\in
{\mathcal L}^1_{\bar{V}}(G)\bigr\}$, that is, the
generalized Fourier transform is surjective.

(ii) Let $\tau$ be a nonzero complex homomorphism
of ${\mathcal A}_{\bar{V}}(\Gamma)$. Since $A(\Gamma)$
is dense in ${\mathcal A}_{\bar{V}}(\Gamma)$, $\tau$
induces a nonzero complex homomorphism of
$A(\Gamma)$, and there exists $\gamma_0\in \Gamma$
such that $\tau(\psi)=\psi(\gamma_0)\ \ (\psi\in
A(\Gamma))$. For any $\psi\in {\mathcal
A}_{\bar{V}}(\Gamma)$, if we choose
$\{\psi_n\}_{n=1}^\infty\subset A(\Gamma)$ such
that
$\lim_{n\rightarrow\infty}\|\psi-\psi_n\|^{\bar{V}}=0$,
then we have \[
\tau(\psi)=\lim_{n\rightarrow\infty}\tau(\psi_n)
=\lim_{n\rightarrow\infty}\psi_n(\gamma_0)=\psi(\gamma_0).
\] Thus $\tau(\psi)=\psi(\gamma_0)\ (\psi\in
{\mathcal
A}_{\bar{V}}(\Gamma))$. This is what we need to
show.

(iii)\ \ From (i) and [18, Corollary 5], we have
$\mathbb{M}({\mathcal A}_{\bar{V}}(\Gamma))\subseteq
C_{BSE({\mathcal A}_{\bar{V}}(\Gamma))}(\Gamma)$. To
prove the reverse inclusion, let $\phi\in
C_{BSE({\mathcal A}_{\bar{V}}(\Gamma))}(\Gamma)$. By
[18, Theorem 4 (i)], there exists a bounded net
$\{\phi_{\lambda}\}_{\lambda\in \Lambda}$ in
${\mathcal A}_{\bar{V}}(\Gamma)$, with a bound $C>0$,
such that $\lim_{\lambda\in
\Lambda}\phi_{\lambda}(\gamma)=\phi(\gamma)\
(\gamma\in \Gamma)$.\ For all $\gamma\in \Gamma$
and all $p\in {\rm span}(\bar{V}+\gamma)$ with
$\|p\|_{\infty}\le 1$, we have \[
\biggl|\sum_{\gamma'\in
\bar{V}+\gamma}\hat{p}(\gamma')\phi(\gamma')\biggl|=
\lim_{\lambda\in \Lambda}\biggl|\sum_{\gamma'\in
\bar{V}+\gamma}\hat{p}(\gamma')
\phi_{\lambda}(\gamma')\biggl| \le
\sup_{\lambda\in \Lambda}\|
\phi_{\lambda}\|^{\bar{V}, \gamma}\le
\sup_{\lambda\in
\Lambda}\|\phi_\lambda\|^{\bar{V}} \le C. \] Thus
we have $\sup_{\gamma\in \Gamma}\|\phi\|^{\bar{V},
\gamma} \le C$, that is, $\phi\in
C_{BSE}^{\bar{V}}(\Gamma)$. Let $\psi\in {\mathcal
A}_{\bar{V}}(\Gamma)$ be arbitrary. Using Lemma
2.1, it follows that 
\[ \lim_{\Omega\in {\mathcal
K}(\Gamma)}\sup_{\gamma\in
\Gamma\setminus\Omega}\|\phi\psi\|^{\bar{V},\gamma}\le
\lim_{\Omega\in {\mathcal K}(\Gamma)}\sup_{\gamma\in
\Gamma\setminus\Omega}\|\phi\|^{\bar{V},\gamma}\|\psi\|^{\bar{V},\gamma}
=0, 
\] which implies $\phi\psi\in {\mathcal
A}_{\bar{V}}(\Gamma)$. Hence
$\phi\in\mathbb{M}({\mathcal A}_{\bar{V}}(\Gamma))$,
and we have $C_{BSE({\mathcal
A}_{\bar{V}}(\Gamma))}(\Gamma) \subseteq
\mathbb{M}({\mathcal A}_{\bar{V}}(\Gamma))$.
Consequently we have $\mathbb{M}({\mathcal
A}_{\bar{V}}(\Gamma))=C_{BSE({\mathcal A}_{\bar{V}}(\Gamma))}(\Gamma)$, that is,
${\mathcal A}_{\bar{V}}(\Gamma)$ is BSE. Since ${\mathcal
A}_{\bar{V}}(\Gamma)$ is regular and has a bounded
approximate identity composed of elements in
${\mathcal A}_{\bar{V}}(\Gamma)_c$ from\, (i) and
(ii), it follows that ${\mathcal A}_{\bar{V}}(\Gamma)$
is also BED by [7, Theorem 4.7], and hence the
second part of (3) follows.
 
(iv)\ The following (4) is already shown in the
proof of (iii) : 
\begin{equation}
C_{BSE}^{\bar{V}}(\Gamma)\subseteq
\mathbb{M}({\mathcal
A}_{\bar{V}}(\Gamma))=C_{BSE({{\mathcal
A}_{\bar{V}}(\Gamma}))}(\Gamma)\subseteq
C_{BSE}^{\bar{V}}(\Gamma), \end{equation} and the
first part of (iv) follows from (4).

Let $\phi\in C^{\bar{V}}_{BSE}(\Gamma)$ be
arbitrary. Then \begin{eqnarray*}
\|\phi\|_{\mathbb{M}({\mathcal
A}_{\bar{V}}(\Gamma))}&=&\sup\bigl\{\|\phi\psi\|^{\bar{V}}:
\psi\in {\mathcal A}_{\bar{V}}(\Gamma),
\|\psi\|^{\bar{V}}\le 1\bigr\} \\ &\le&
\sup\bigl\{\|\phi\|^{\bar{V}}\|\psi\|^{\bar{V}}:
\psi\in {\mathcal A}_{\bar{V}}(\Gamma),
\|\psi\|^{\bar{V}}\le 1\bigr\} \le
\|\phi\|^{\bar{V}}. \end{eqnarray*} On the other
hand, let $\varepsilon>0$ and
$\gamma_\varepsilon\in \Gamma$ be such that
$\|\phi\|^{\bar{V}}-\varepsilon\le
\|\phi\|^{\bar{V}, \gamma _\varepsilon}$. Choose
$u\in {\mathcal A}_{\bar{V}}(\Gamma)$ such that
$\|u\|^{\bar{V}}\le 1+\varepsilon, u(\gamma)=1 \
(\gamma\in \bar{V}+\gamma_\varepsilon) \, ({\rm \,
[17, Theorem\, 2.6.8\,])}$. It follows that
\begin{equation} \|\phi\|_{\mathbb{M}({\mathcal
A}_{\bar{V}}(\Gamma))}\ge \frac{\|\phi
u\|^{\bar{V}}} {\|u\|^{\bar{V}}} \ge \frac{\|\phi
u\|^{\bar{V}, \gamma_\varepsilon}}
{\|u\|^{\bar{V}}}=\frac{\|
\phi\|^{\bar{V},\gamma_\varepsilon}}{\|u\|^{\bar{V}}}
\ge
\frac{\|\phi\|^{\bar{V}}-\varepsilon}{1+\varepsilon}.
\end{equation} Since $\varepsilon>0$ in (5) can be
taken arbitrarily small,
$\|\phi\|_{\mathbb{M}({\mathcal A}_{\bar{V}}
(\Gamma))} \ge\|\phi\|^{\bar{V}}$ follows. Thus
the second part of (iv) also holds. \ \ \ \ $\Box$
\vspace{0.5cm}

{\bf Definition 3.2.}\ (Generalized Fourier
algebra).\ By Theorem 3.1 (i), the generalized
Fourier transform is an isometric homomorphism of
${\mathcal L}^1_{\bar{V}}(G)$ onto ${\mathcal
A}_{\bar{V}}(\Gamma)$. We call $({\mathcal
A}_{\bar{V}}(\Gamma), \|\cdot \|^{\bar{V}})$
generalized Fourier algebra of $\Gamma$.
\vspace{0.5cm}

{\bf Definition 3.3.}\ (Generalized
Fourier-Stieltjes algebra).\ By Definitions 2.4
and 3.1, the generalized Fourier-Stieltjes
transform is an isometric homomorphism of ${\mathcal
M}_{\bar{V}}(G)$ onto ${\mathcal
B}_{\bar{V}}(\Gamma)$. We call $({\mathcal
B}_{\bar{V}}(\Gamma), \|\cdot \|^{\bar{V}})$
generalized Fourier-Stieltjes algebra of $\Gamma$.

\vspace{0.5cm}

From Theorem 3.1, we obtain the following
corollary.

 \it

\vspace{0.5cm}

{\bf Corollary 3.2.}\ {\rm (i)}\ $({\mathcal
L}^1_{\bar{V}}(G), \|\cdot \|_{\bar{V}})$ is a
commutative semisimple regular Banach algebra, and
its Gelfand transform coincides with the
generalized Fourier transform.
 
{\rm (ii)}\ $({\mathcal L}^1_{\bar{V}}(G), \|\cdot
\|_{\bar{V}})$ has a bounded approximate identity
$\{e_\lambda\}_{\lambda\in \Lambda}$ such that
$\widehat{e_{\lambda}}$ has compact support for
each $\lambda\in \Lambda$.
 
{\rm (iii)}\ $({\mathcal L}^1_{\bar{V}}(G), \|\cdot
\|_{\bar{V}})$ is BSE and BED.
\vspace{0.5cm}

\it

{\bf Lemma 3.3.}\ $\widehat{k\mu}(\gamma)=
\hat{k}*\hat{\mu}(\gamma)\ \ (\gamma\in \Gamma,
\mu \in M(G), k\in L^1(G)\cap B(G))$.

\rm 

Proof. By the inversion theorem, we have
$\hat{k}\in L^1(\Gamma)$ with $k(x)=\int_\Gamma
(x,\gamma)\hat{k}(\gamma)d\gamma$. Hence
\begin{eqnarray*}
 \widehat{k\mu}(\gamma)&=&\int_G
(-x, \gamma)k(x)d\mu(x) =\int_G (-x,
\gamma)\biggl(\int_\Gamma(x,\gamma')
\hat{k}(\gamma')d\gamma'\biggr)d\mu(x) \\ &=&
\int_\Gamma\biggl(\int_G (-x, \gamma-\gamma')
d\mu(x)\biggr)\hat{k}(\gamma')d\gamma'
=\int_\Gamma
\hat{\mu}(\gamma-\gamma')\hat{k}(\gamma') d\gamma'
\\ &=&\hat{k}*\hat{\mu}(\gamma)\ \ (\gamma\in
\Gamma). \ \ \ \Box \
\end{eqnarray*}

{\bf Theorem 3.4.}\ {\rm (i)}\ Let $0\neq \mu$ be
a measure in $M_0(G)$. For any $\varepsilon>0$,
there exists $\Omega_\varepsilon \in {\mathcal
K}(\Gamma)$ such that $\sup_{\gamma\in
\Gamma\setminus
\Omega_\varepsilon}\|\hat{\mu}\|^{\bar{V}, \gamma}
\le\varepsilon$.

{\rm (ii)}\ $M_0(G)= {\mathcal L}^1_{\bar{V}}(G)\cap
M(G)$. \rm

Proof.\ (i)\ Let $\varepsilon>0$ be arbitrary.
Choose $u_{\bar{V}}\in A(\Gamma)_c$ such that
$u_{\bar{V}}(\gamma)=1\ (\gamma\in \bar{V})$, and
put $e_{\bar{V}}(x)=
\int_\Gamma(x,\gamma)u_{\bar{V}}(\gamma)d\gamma\
(x\in G)$. By the inversion theorem,
$e_{\bar{V}}\in L^1(G) \cap C_0(G)$ with
$\widehat{e_{\bar{V}}}(\gamma)
=u_{\bar{V}}(\gamma)=1\ \ (\gamma\in \bar{V})$.

Choose $K\in {\mathcal K}(G)$ such that
\begin{equation} \int_{G\setminus K} \biggl(\int_G
|e_{\bar{V}}(y-z)| d|\mu|(z)\biggr)dy \le
\varepsilon/2. \end{equation} For each $x\in G$,
$e_{\bar{V}}(x-z)\mu(z)\in M_0(G)$. Indeed, the
Fourier-Stieltjes transform of
$e_{\bar{V}}(x-z)\mu(z)$ is
$\widehat{(e'_{\bar{V}})}_x*\hat{\mu}$ by Lemma
3.3, and since $\widehat{(e'_{\bar{V}})_x}(\gamma)
=(-x,\gamma)\widehat{e_{\bar{V}}}(-\gamma)\in
C_c(\Gamma)$ and $\hat{\mu}\in C_0(\Gamma)$, we
have $\widehat{({e'_{\bar{V}})}_x}*\hat{\mu}\in
C_0(\Gamma)$. Therefore \[ g_\gamma(x):=
(x,\gamma)\int_G (-z,
\gamma)e_{\bar{V}}(x-z)d\mu(z) \ \ \ (\gamma\in
\Gamma) \] vanishes at infinity, and there is
$\Omega_x\in {\mathcal K}(\Gamma)$ such that
$|g_\gamma(x)|\le \varepsilon/ (4m_G(K)) \ \
(\gamma \in \Gamma\setminus \Omega_x)$. Since
$e_{\bar{V}}\in C_0(G)$, we can choose a
neighborhood $U_0$ of $0\in G$ such that
$|e_{\bar{V}}(y)-e_{\bar{V}}(z)|
\le\varepsilon/(4m_G(K)\|\mu\|)\ \ (y-z\in U_0)$.
Then \begin{eqnarray} |g_\gamma(y)| &=& \biggl|(y,
\gamma)\int_G (-z, \gamma)
e_{\bar{V}}(y-z)d\mu(z)\biggl| \nonumber \\ &\le&
\biggl|\int_G (-z, \gamma)(e_{\bar{V}}(y-z)-
e_{\bar{V}}(x-z))d\mu(z)\biggr| \nonumber \\
&&\hspace{2cm}+\biggl|\int_G(-z,\gamma)
e_{\bar{V}}(x-z)d\mu(z)\biggr| \nonumber \\ &\le&
\frac{\varepsilon \|\mu\|}{4m_G(K)\|\mu\|}
+\frac{\varepsilon}{4m_G(K)}=\frac{\varepsilon}
{2m_G(K)}\ \ (y\in x+U_0, \gamma\in
\Gamma\setminus \Omega_x). \end{eqnarray} Choose a
finite sub-covering $\{x_j+U_0\}_{j=1}^n$ of $K$
from the covering $\{x+U_0\}_{x\in K}$\ of $K$.

Suppose $\gamma\in \Gamma\setminus
\cup_{j=1}^n\Omega_{x_j}$. Since each $y\in K$
belongs to some $x_j+U_0$, we have
$|g_\gamma(y)|\le \varepsilon/(2m_G(K))$ by (7).
From this result and (6), it follows that
\begin{eqnarray} \|g_\gamma\|_1&=&
\int_{G\setminus K}|g_\gamma(x)|dx +\int_K
|g_\gamma(x)|dx \nonumber \\ &\le&
\int_{G\setminus K}\biggl(\int_G|e_{\bar{V}}(x-z)|
d|\mu|(z)\biggr)dx +\int_K |g_\gamma(x)|dx
\nonumber \\ &\le& \frac{\varepsilon}{2}+\int_K
\frac{\varepsilon} {2m_G(K)}dx=\varepsilon\ \ \
(\gamma\in \Gamma\setminus
(\cup_{j=1}^n\Omega_{x_j})). \end{eqnarray} Since
\[ \widehat{g_\gamma}(\gamma')
=\widehat{e_{\bar{V}}}
(\gamma'-\gamma)\hat{\mu}(\gamma')
=\hat{\mu}(\gamma')\ \ \ (\gamma'\in
\bar{V}+\gamma), \] it follows from (8) that \[
\|\hat{\mu}\|^{\bar{V},
\gamma}=\|\widehat{g_\gamma}\|^{\bar{V},\gamma}
\le \|g_\gamma\|_1\le \varepsilon\ \ \
(\gamma\not\in \cup_{j=1}^n\Omega_{x_j}). \]

(ii). This easily follows from (i).\ \ $\Box$

\vspace{0.5cm}

{\bf Theorem 3.5.}\ {\rm (i)}
\ $L^1(G)\subsetneqq M_0(G)\subsetneqq 
{\mathcal
L}^1_{\bar{V}}(G)$\ \ {\rm (ii)}\ $M(G)\subsetneqq
{\mathcal M}_{\bar{V}}(G)$.

 \rm

Proof.\ (i)\ $L^1(G)\subsetneqq M_0(G)$ follows
from [5, Theorem 7.4.1] and $M_0(G)\subseteq 
{\mathcal L}^1_{\bar{V}}(G)$ follows from Theorem 3.4.
Further, if $M_0(G)={\mathcal L}^1_{\bar{V}}(G),\,
\|\cdot\|$ and $\|\cdot\|_{\bar{V}}$ must be
equivalent. From this we have $L^1(G)={\mathcal
L}^1_{\bar{V}}(G)$ contradicting what we showed
above. Hence $M_0(G)\subsetneqq {\mathcal
L}^1_{\bar{V}}(G)$

(ii)\ If $M(G)={\mathcal M}_{\bar{V}}(G)$, $\|\cdot\|$
and $\|\cdot\|_{\bar{V}}$ are equivalent norms.
This is impossible as is shown in (i). \ \ \
$\Box$
\vspace{0.5cm}

{\bf \S 4. Segal algebras in ${\mathcal
A}_{\bar{V}}(\Gamma), $ I}\ \ The notion of a
Segal algebra in $L^1(G)$ was introduced by H.
Reiter [13], and the theory and results on Segal
algebras developed until the end of 1980's were
contained in the monographs [13, 14, 15], and
[16]. Especially, the invention of the smallest
isometrically character invariant Segal algebra by
H. Feichtinger [3] is famous and important (cf.
[3], [16], [10]).

In [8], the present authors extended the notion of
Reiter's Segal algebras in $L^1(G)$ to the notion
of Segal algebras in non-unital commutative
semisimple Banach algebra $({\bf B}, \| \cdot
\|_{\bf B})$ which satisfies:

($\alpha_{\bf B}$) ${\bf B}$ is regular.

($\beta_{\bf B}$) There\ is\ a bounded approximate
identity of ${\bf B}$ composed of elements in
${\bf B}_c$. \vspace{0.1cm}

\vspace{0.3cm} 

\noindent An ideal ${\mathcal S}$ in ${\bf B}$ is
called a Segal algebra in ${\bf B}$ if ${\mathcal S}$
satisfies:

{\rm (i)}\ ${\mathcal S}$\ {\rm is\ dense\ in}\ ${\bf
B}$;

{\rm (ii)}\ ${\mathcal S}$ is a Banach space under a
norm \ $\|\cdot \|_{\mathcal S}$\ which satisfies
$\|f\|_{\bf B}\le\|f\|_{\mathcal S}\ \ (f\in {\mathcal
S})$;

{\rm (iii)}\ $\|fg\|_{\mathcal S}\le \|f\|_{\mathcal
S}\|g\|_{\bf B}\ \ (f\in {\mathcal S}, g\in {\bf B})$;

{\rm (iv)}\ ${\mathcal S}$ has approximate units.
\vspace{0.2cm}

When ${\bf B}=L^1(G)$, Segal algebras defined
above coincide with Reiter's Segal algebras
([8]).\, Since $B={\mathcal
L}^1_{\bar{V}}(G)\,({\rm resp}. {\mathcal
A}_{\bar{V}}(\Gamma), A(\Gamma))$ satisfies
conditions ($\alpha_{\bf B}$) and ($\beta_{\bf
B}$) by Corollary 3.2, we can define Segal
algebras in ${\mathcal L}^1_{\bar{V}}(G)\,({\rm resp}.
\,{\mathcal A}_{\bar{V}}(\Gamma), \,A(\Gamma))$.
\vspace{0.5cm}

{\bf Definition 4.1.}\ Let ${\mathfrak S}_{\bar{
V}}(G)$ and $\tilde{\mathfrak S}_{\bar{V}}(G)$ be
subsets of $L^1(G)$ defined by \begin{eqnarray*}
&&{\mathfrak S}_{\bar{V}}(G):=\bigl\{f\in L^1(G):
{\rm supp}(\hat{f})\subseteq \bar{V}\bigr\}, \\ &&
\tilde{\mathfrak S}_{\bar{V}}(G):=\bigl\{f\in
L^1(G) : \, {\rm there\, exists}\, \gamma\in
\Gamma\, {\rm such\, that}\ {\rm
supp}(\hat{f})\subseteq \bar{V}+\gamma\bigr\}.
\end{eqnarray*}
For a sequence
$\{f_n\}_{n=1}^\infty$ of elements in
$\tilde{\mathfrak S}_{\bar{V}}(G)$ with
$\displaystyle
\sum_{n=1}^{\infty}\|f_n\|_1<\infty$, there is
$f\in L^1(G)$ such that
$\|f-\sum_{n=1}^Nf_n\|_1\rightarrow 0\,
(N\rightarrow \infty)$. Hence
$f=\sum_{n=1}^{\infty}f_n$, which is called a
$\bar{V}$-representation. 
Define 
\begin{eqnarray*}
&&{\mathcal S}^1_{\bar{V}}(G):=\bigl\{f\in L^1(G):\ f
\ {\rm has\ at\ least\ one\ }\bar{V}-{\rm
representation}\bigr\}, 
\\ && \|f\|_{{\mathcal
S}^1_{\bar{V}}}:=\inf\biggl\{\sum_{n=1}^{\infty}\|f_n\|_1:
f=\sum_{n=1}^{\infty}f_n\ (\bar{V}-{\rm
representation})\biggr\}\ \ \ (f\in {\mathcal
S}^1_{\bar{V}}(G)). \end{eqnarray*}

$({\mathcal S}^1_{\bar{V}}(G), \|\cdot \|_{{\mathcal
S}^1_{\bar{V}}})$ becomes a Segal algebra which is
isometrically character invariant\, (i.e.,
$\gamma\in \Gamma$ and $ f\in {\mathcal
S}^1_{\bar{V}}(G)$ imply $(x,\gamma)f\in
S^1_{\bar{V}}(G)$ with $\|(x,\gamma)f\|_{{\mathcal
S}^1_{\bar{V}}}=\|f\|_{{\mathcal S}^1_{\bar{V}}}$),
and also is the smallest one in the family of all
isometrically character invariant Segal algebras
([3]). ${\mathcal S}^1_{\bar{V}}(G)$ is called the
Feichtinger algebra of $G$ ([16, \S 6.2]).
 
\vspace{0.5cm}

{\bf Definition 4.2.}\ Let
$\Lambda_{\bar{V}}(\Gamma)$ and
$\tilde{\Lambda}_{\bar{V}}(\Gamma)$ be subsets of
$A(\Gamma)$ defined by \begin{eqnarray*}
&&\Lambda_{{\bar{V}}}(\Gamma) :=\bigl\{\psi\in
A(\Gamma): {\rm supp}(\psi)\subseteq
\bar{V}\bigr\}, \\ &&
\tilde{\Lambda}_{\bar{V}}(\Gamma): =\bigl\{\psi\in
A(\Gamma) : \, {\rm there\, exists}\,\gamma\in
\Gamma\, {\rm such\, that}\ {\rm
supp}(\psi)\subseteq \bar{V}+\gamma\bigr\}.
\end{eqnarray*} For a sequence
$\{\psi_n\}_{n=1}^\infty $ of elements in
$\tilde{\Lambda}_{\bar{V}}(\Gamma)$ with
$\displaystyle
\sum_{n=1}^{\infty}\|\psi_n\|_A<\infty$, there is
$\psi\in A(\Gamma)$ such that
$\|\psi-\sum_{n=1}^N\psi_n\|_A\rightarrow 0
(N\rightarrow \infty)$. Hence
$\psi=\sum_{n=1}^{\infty}\psi_n$, which is called
a $\bar{V}$-representation. Define
\begin{eqnarray*}
 &&{\mathcal S}_{\bar{V}}(\Gamma):
=\bigl\{\psi\in A(\Gamma): \psi\ {\rm has\ at\
least\ one\ }\bar{V}-{\rm representation}\bigr\},
\\ &&\|\psi\|_{{\mathcal S}_{\bar{V}}}:
=\inf\biggl\{\sum_{n=1}^{\infty}\|\psi_n\|_A:
\psi=\sum_{n=1}^\infty \psi_n\ (\bar{V}-{\rm
representation})\biggr\} \ \ (\psi\in {\mathcal
S}_{\bar{V}}(\Gamma)). 
\end{eqnarray*}

Since the Fourier transform is an isometric
homomorphism of $L^1(G)$ onto $A(\Gamma)$ with
''$\widehat{(x,\gamma_0)f}(\gamma)=\hat{f}(\gamma-\gamma_0)
\ (f\in L^1(G), \gamma_0, \gamma\in \Gamma)$'', it
readily follows that $({\mathcal S}_{\bar{V}}(\Gamma),
\|\cdot \|_{{\mathcal S}_{\bar{V}}})$ is an
isometrically translation invariant Segal algebra
in $A(\Gamma)$ which is the smallest one in the
family of all isometrically translation invariant
Segal algebras in $A(\Gamma)$.
 
\vspace{0.5cm}

{\bf Definition 4.3.}\ Let $\Lambda_{\bar{V},
{\mathcal A}}(\Gamma)$ and $\tilde{\Lambda}_{\bar{V},
{\mathcal A}}(\Gamma)$ be subsets of ${\mathcal
A}_{\bar{V}}(\Gamma)$ defined by \begin{eqnarray*}
&&\Lambda_{{\bar{V}, {\mathcal A}}}(\Gamma)
:=\bigl\{\psi\in {\mathcal A}_{\bar{V}}(\Gamma): {\rm
supp}(\psi)\subseteq \bar{V}\bigr\}, \\ &&
\tilde{\Lambda}_{\bar{V}, {\mathcal A}}(\Gamma):
=\bigl\{\psi\in {\mathcal A}_{\bar{V}}(\Gamma) :\ {\rm
there\ exists}\ \gamma\in \Gamma\ {\rm such\
that}\ {\rm supp}(\psi)\subseteq
\bar{V}+\gamma\bigr\}. \end{eqnarray*} For a
sequence $\{\psi_n\}_{n=1}^\infty $ of elements in
$\tilde{\Lambda}_{\bar{V}, {\mathcal A}}(\Gamma)$ with
$\displaystyle
\sum_{n=1}^{\infty}\|\psi_n\|^{\bar{V}}<\infty$,
there is $\psi\in {\mathcal A}_{\bar{V}}(\Gamma)$ such
that
$\|\psi-\sum_{n=1}^N\psi_n\|^{\bar{V}}\rightarrow
0 (N\rightarrow \infty)$. Hence
$\psi=\sum_{n=1}^{\infty}\psi_n$, which we call a
$\bar{V}$-representation. Define \begin{eqnarray*}
&&{\mathcal S}_{\bar{V}, {\mathcal A}}(\Gamma):
=\bigl\{\psi\in {\mathcal A}_{\bar{V}}(\Gamma): \psi\
{\rm has\ at\ least\ one\ }\bar{V}-{\rm
representation}\bigr\}, \\ &&\|\psi\|_{{\mathcal
S}_{\bar{V}, {\mathcal A}}}:
=\inf\biggl\{\sum_{n=1}^{\infty}\|\psi_n\|^{\bar{V}}:
\psi=\sum_{n=1}^\infty \psi_n\ (\bar{V}-{\rm
representation})\biggr\} \ \ (\psi\in {\mathcal
S}_{\bar{V}, {\mathcal A}}(\Gamma)). \end{eqnarray*}
It is easy to see that ${\mathcal S}_{\bar{V}, {\mathcal 
A}}(\Gamma)$ is a normed linear space with norm
$\|\ \|_{{\mathcal S}_{\bar{V}, {\mathcal A}}}$, which
dominates $\|\cdot \|^{\bar{V}}$. \vspace{0.5cm}

\it

{\bf Lemma 4.1.}\ $({\mathcal S}_{\bar{V}, {\mathcal
A}}(\Gamma), \|\cdot \|_{{\mathcal S}_{\bar{V}, {\mathcal
A}}})$ is a Banach space.

\rm Proof. Let $\{\phi_n\}_{n=1}^\infty$ be a
Cauchy sequence in ${\mathcal S}_{\bar{V}, {\mathcal A}}(\Gamma)$. We will show that there is $\phi\in
{\mathcal S}_{\bar{V}, {\mathcal A}}(\Gamma)$ such that
$\|\phi-\phi_n\|_{{\mathcal S}_{\bar{V}}, 
{\mathcal A}}\rightarrow 0 (n\rightarrow \infty)$. To do
this we may suppose that 
$\|\phi_{n+1}-\phi_n\|_{{\mathcal S}_{\bar{V}, 
{\mathcal A}}}\le 1/(n+1)^2\,(n=1,2,3,...)$.} Put
$\psi_n:=\phi_n-\phi_{n-1}\ (n=1,2,3,...)$, where
$\phi_0=0$. Since $\|\psi_{n+1}\|^{\bar{V}}\le
\|\psi_{n+1}\|_{{\mathcal S}_{\bar{V}, {\mathcal A}}}\le
1/(n+1)^2\ (n\ge 1),\ \sum_{n=1}^\infty\psi_n$ is
a convergent series in ${\mathcal A}_{\bar{V}}(\Gamma)$ converging to a function,
say $\phi\in {\mathcal A}_{\bar{V}}(\Gamma)$. For each
$n\ge 1$, let $\psi_n=\sum_{m=1}^\infty
\psi_{n,m}$ be a $\bar{V}$-representation such
that $\sum_{m=1}^\infty \|\psi_{n,
m}\|^{\bar{V}}\le \|\psi_n\|_{{\mathcal S}_{\bar{V},
{\mathcal A}}} +1/n^2$. Then $\phi=\sum_{1\le m,
n<\infty} \psi_{n, m}$ is a
$\bar{V}$-representation since 
\[ \sum_{1\le m, n<\infty} \|\psi_{n,
m}\|^{\bar{V}} \le \sum_{n=1}^\infty
\bigl(\|\psi_n\|_{{\mathcal S}_{\bar{V}, {\mathcal 
A}}}+1/n^2\bigr)\le \|\psi_1\|_{{\mathcal S}_{\bar{V},
{\mathcal A}}}+1+2\sum_{n=2}^\infty
\frac{1}{n^2}<\infty. 
\]
 Hence $\phi\in
S_{\bar{V}, {\mathcal A}}(\Gamma)$ with
\begin{eqnarray*}
 \|\phi-\phi_N\|_{{\mathcal
S}_{\bar{V}, {\mathcal A}}} &=&\|\phi-\sum_{n=1}^N
\psi_n\|_{{\mathcal S}_{\bar{V}, {\mathcal A}}}
\le\sum_{n=N+1}^\infty\sum_{m=1}^\infty \|\psi_{n,
m}\|^{\bar{V}} \\ &\le& \sum_{n=N+1}^\infty
(\|\psi_n\|_{{\mathcal S}_{\bar{V}, {\mathcal A}}}
+\frac{1}{n^2}) \rightarrow 0 (N\rightarrow
\infty). \ \ \ \ \ \Box \end{eqnarray*}

\vspace{0.5cm}

\it

{\bf Lemma 4.2} {\rm (i)}
$\Lambda_{\bar{V}}(\Gamma)= \Lambda_{\bar{V},
{\mathcal A}}(\Gamma)$, and there exists a constant
$C_1>0$ such that \[\|\phi\|_A\le
C_1\|\phi\|^{\bar{V}}\ \ (\phi\in
\Lambda_{\bar{V}}(\Gamma)). \]

{\rm (ii)}\ If $({\mathcal S}, \|\cdot \|_{\mathcal S})$
is an isometrically translation invariant Segal
algebra in ${\mathcal A}_{\bar{V}}(\Gamma)$, we have
$\Lambda_{\bar{V}, {\mathcal A}}(\Gamma)\subset {\mathcal 
S}$, and there exists a constant $C_2>0$ such that
\[\|\phi\|_{\mathcal S}\le C_2\|\phi\|^{\bar{V}} \ \
(\phi\in \Lambda_{\bar{V}, {\mathcal A}}(\Gamma)). \]

\rm 

Proof.\ (i)\ Note that $\Lambda_{\bar{V}}(\Gamma)$
and $\Lambda_{\bar{V}, {\mathcal A}}(\Gamma)$ are
closed subspaces of $A(\Gamma)$ and ${\mathcal 
A}_{\bar{V}}(\Gamma)$, respectively. The inclusion
$\subseteq$ is trivial. To show $\supseteq$,
suppose $\phi\in \Lambda_{\bar{V}, {\mathcal A}}(\Gamma)$. Take $\xi\in A(\Gamma)_c$ such that
$\xi(\gamma)=1\ (\gamma\in \bar{V})$. By [8,
Theorems 7.4], we have $\phi\in
C_{BSE}^{\bar{V}}(\Gamma)\subset A(\Gamma)_{loc}
$, and hence there is $\psi\in A(\Gamma)$ such
that $\phi(\gamma)=\psi(\gamma)\ (\gamma\in {\rm
supp}(\xi))$ from [8, Proposition 7.2]. Then
$\phi=\psi\xi\in A(\Gamma)$ and so $\phi\in
\Lambda_{\bar{V}}(\Gamma)$.
 
Hence $\Lambda_{\bar{V}}(\Gamma)=
\Lambda_{\bar{V}, {\mathcal A}}(\Gamma)$, and the
identity map of $\Lambda_{\bar{V}}(\Gamma)$ to
$\Lambda_{\bar{V}, {\mathcal A}}(\Gamma)$ is norm
decreasing. There exists, by the open mapping
theorem, a constant $C_1>0$ such that
$\|\phi\|_A\le C_1\|\phi\|^{\bar{V}}\ \ (\phi\in
\Lambda_{\bar{V}}(\Gamma)).$

(ii)\ If we put ${\mathcal S}_{\bar{V}}=\{\phi\in
{\mathcal S}: {\rm supp}(\phi)\subset \bar{V}\}$, we
have ${\mathcal S}_{\bar{V}}=\Lambda_{\bar{V},
 {\mathcal A}}(\Gamma)$ by [8, Lemma 3.4], and the identity
map of ${\mathcal S}_{\bar{V}}$ onto
$\Lambda_{\bar{V}, {\mathcal A}}(\Gamma)$ is norm
decreasing. There is, by the open mapping theorem,
$C_2>0$ such that $\|\phi\|_{\mathcal S}\le
C_2\|\phi\|^{\bar{V}}\ (\phi \in \Lambda_{\bar{V},
{\mathcal A}}(\Gamma))$.\ $\Box$
 
\vspace{0.5cm}

\it

{\bf Theorem 4.3.} {\rm (i)} ${\mathcal S}_{\bar{V}}(\Gamma)={\mathcal S}_{\bar{V}, {\mathcal 
A}}(\Gamma)$, and the identity map of ${\mathcal 
S}_{\bar{V}}(\Gamma)$ to ${\mathcal S}_{\bar{V}, {\mathcal 
A}}(\Gamma)$ is a bicontinuous algebra
homomorphism.

{\rm (ii)}\ $({\mathcal S}_{\bar{V}, {\mathcal A}}(\Gamma), \|\cdot \|_{{\mathcal S}_{\bar{V}, 
{\mathcal A}}})$ is a Segal algebra in 
${\mathcal A}_{\bar{V}}(\Gamma)$.

{\rm (iii)}\ ${\mathcal S}_{\bar{V}, {\mathcal 
A}}(\Gamma)$ is isometrically translation
invariant, which is the smallest one in the family
of all isometrically translation invariant Segal
algebras in ${\mathcal A}_{\bar{V}}(\Gamma)$.

\rm 

Proof.\ (i)\ Let $\phi\in S_{\bar{V}}(\Gamma)$
with a ${\bar{V}}$-representation
$\phi=\sum_{n=1}^\infty \phi_n. $ Since $\phi_n\in
\tilde{\Lambda}_{\bar{V},{\mathcal A}} (\Gamma)$ with
$\|\phi_n\|^{\bar{V}}\le \|\phi_n\|_A,
\,n=1,2,3,...$, we have
$\sum_{n=1}^\infty\|\phi_n\|^{\bar{V}}<\infty$.
Therefore $\phi\in {\mathcal S}_{\bar{V}, {\mathcal 
A}}(\Gamma)$, and ${\mathcal S}_{\bar{V}}(\Gamma)
\subseteq {\mathcal S}_{\bar{V}, {\mathcal A}}(\Gamma)$
holds.

To prove ${\mathcal S}_{\bar{V}}(\Gamma) \supseteq
{\mathcal S}_{\bar{V}, {\mathcal A}}(\Gamma)$, let
$\psi\in {\mathcal S}_{\bar{V}, {\mathcal A}}(\Gamma)$
with a $\bar{V}$-representation
$\psi=\sum_{n=1}^\infty \psi_n.$ By Lemma 4.2 (i)
and the isometrically translation invariance of
$A(\Gamma)$ and ${\mathcal A}_{\bar{V}}(\Gamma)$,
respectively, we have $\|\psi_n\|_A\le
C_1\|\psi_n\|^{\bar{V}}, n=1,2,3,...$. Hence
$\sum_{n=1}^\infty\|\psi_n\|_A<\infty$, and
$\psi\in {\mathcal S}_{\bar{V}}(\Gamma)$ with
$\bar{V}$-representation
$\psi=\sum_{n=1}^\infty\psi_n$. Therefore ${\mathcal 
S}_{\bar{V}}(\Gamma)\supseteq {\mathcal S}_{\bar{V},
{\mathcal A}}(\Gamma)$, and ${\mathcal 
S}_{\bar{V}}(\Gamma)= {\mathcal S}_{\bar{V}, {\mathcal 
A}}(\Gamma)$ holds. The identity map of 
${\mathcal 
S}_{\bar{V}}(\Gamma)$ onto ${\mathcal S}_{\bar{V},
{\mathcal A}}(\Gamma)$ is a bicontinuous homomorphism
by the open mapping theorem.

(ii)\ Let $\psi\in {\mathcal S}_{\bar{V}, {\mathcal 
A}}(\Gamma)$ with a $\bar{V}$-representation
$\psi=\sum_{n=1}^\infty \psi_n$. Then
$\|\psi\|^{\bar{V}} \le\sum_{n=1}^\infty
\|\psi_n\|^{\bar{V}}$. Taking the infimum of all
over the $\bar{V}$-representations of $\psi$, we
obtain $\|\psi\|^{\bar{V}} \le \|\psi\|_{{\mathcal 
S}_{\bar{V}, {\mathcal A}}}$.

Next, we show that ${\mathcal S}_{\bar{V}, 
{\mathcal 
A}}(\Gamma)$ is a Banach ideal of ${\mathcal 
A}_{\bar{V}}(\Gamma)$. ${\mathcal S}_{\bar{V}, 
{\mathcal 
A}}(\Gamma)$ is a Banach space from Lemma 4.1. For
$\phi\in {\mathcal A}_{\bar{V}}(\Gamma)$ and $\psi\in
{\mathcal S}_{\bar{V}, {\mathcal A}}(\Gamma)$ with a
$\bar{V}$-representation $\psi=\sum_{n=1}^\infty
\psi_n$, $\phi\psi=\sum_{n=1}^\infty\phi\psi_n$ is
a $\bar{V}$-representation, hence $\phi\psi\in
{\mathcal S}_{\bar{V}, {\mathcal A}}(\Gamma)$.
Furthermore, we have \[ \|\phi \psi\|_{{\mathcal 
S}_{\bar{V}, {\mathcal A}}}\le\sum_{n=1}^\infty
\|\phi\psi_n\|^{\bar{V}}\le \sum_{n=1}^\infty
\|\phi\|^{\bar{V}} \|\psi_n\|^{\bar{V}}. \] Taking
the infimum of all over the
$\bar{V}$-representations of $\psi$, we obtain
$\|\phi\psi\|_{{\mathcal S}_{\bar{V},{\mathcal 
A}}}\le\|\phi\|^{\bar{V}} \|\psi\|_{{\mathcal 
S}_{\bar{V}, {\mathcal A}}} $.

Since ${\mathcal S}_{\bar{V}}(\Gamma) $ is dense in
$A(\Gamma)$, and $A(\Gamma)$ is dense in 
${\mathcal 
A}_{\bar{V}}(\Gamma)$, it follows that 
${\mathcal S}_{\bar{V}}(\Gamma)$ is dense in 
${\mathcal 
A}_{\bar{V}}(\Gamma)$, and so (i) implies 
${\mathcal 
S}_{\bar{V}, {\mathcal A}}(\Gamma)$ is dense in ${\mathcal 
A}_{\bar{V}}(\Gamma)$.

Let $\psi\in {\mathcal S}_{\bar{V}, {\mathcal A}}(\Gamma)$
and $\varepsilon>0$ be given, and let
$\psi=\sum_{n=1}^{\infty} \psi_n$ be a
$\bar{V}$-representation. Choose $n_0\in
\mathbf{N}$ so that
$\displaystyle\sum_{n=n_0+1}^{\infty}
\|\psi_n\|^{\bar{V}} <\varepsilon/3$, and put
$\Omega_{n_0} =\cup_{n=1}^{n_0}{\rm
supp}(\psi_n)\in {\mathcal K}(\Gamma)$. Let
$u_{\Omega_{n_0}}\in A(\Gamma)_c$ with
$\|u_{\Omega_{n_0}}\|_A\le 2$ and
$u_{\Omega_{n_0}}(\gamma)=1 \ \ (\gamma\in
\Omega_{n_0})$. From [8, Lemma 3.4] and (i), we
have $u_{\Omega_{n_0}}\in {\mathcal 
S}_{\bar{V}}(\Gamma)= {\mathcal S}_{\bar{V},{\mathcal 
A}}(\Gamma)$. Also, since
\[\displaystyle\psi-u_{\Omega_{n_0}}\psi
=\sum_{n=1}^\infty\psi_n-\sum_{n=1}^\infty
u_{\Omega_{n_0}}\psi_n =\sum_{n=n_0+1}^\infty
(1-u_{\Omega_{n_0}})\psi_n \] is a
$\bar{V}$-representation, it follows that
\begin{eqnarray*}
\|\psi-u_{\Omega_{n_0}}\psi\|_{{\mathcal 
S}_{\bar{V},{\mathcal A}}} &\le&
\sum_{n=n_0+1}^{\infty}\bigl\|\psi_n
-u_{\Omega_{n_0}}\psi_n \bigr\|^{\bar{V}} 
\\ &\le&
\sum_{n=n_0+1}^{\infty}\|\psi_n\|^{\bar{V}}
(1+\|u_{\Omega_{n_0}}\|^{\bar{V}})\le
\frac{\varepsilon}{3}(1+\|u_{\Omega_{n_0}}\|_A)
\le\varepsilon. 
\end{eqnarray*}
 This shows that
${\mathcal S}_{\bar{V}, {\mathcal A}}(\Gamma)$ has
approximate units. Hence $({\mathcal S}_{\bar{V},{\mathcal 
A}}(\Gamma), \|\cdot \|_{{\mathcal S}_{{\bar{V},
 {\mathcal A}}}})$ is a Segal algebra in ${\mathcal 
A}_{\bar{V}}(\Gamma)$.

(iii)\ Let $\psi\in {\mathcal S }_{\bar{V}, 
{\mathcal 
A}}(\Gamma)$. For $\gamma\in \Gamma$, if
$\psi_\gamma$ denotes the translation of $\psi$ by
$\gamma, \ \psi_\gamma\in {\mathcal S}_{\bar{V}, {\mathcal A}}(\Gamma)$ is obvious, and \begin{eqnarray*}
\|\psi\|_{{\mathcal S}_{\bar{V}, {\mathcal A}}}&=& \inf
\biggl\{\sum_{n=1}^{\infty}\|\psi_n\|^{\bar{V}}:
\psi=\sum_{n=1}^{\infty}\psi_n (\bar{V}-{\rm
representation})\biggr\} \\ &=& \inf
\biggl\{\sum_{n=1}^{\infty}\|(\psi_n)_\gamma\|^{\bar{V}}:
\psi_\gamma=\sum_{n=1}^{\infty}(\psi_n)_\gamma
\bigl(\psi= \sum_{n=1}^\infty \psi_n (\bar{V}-{\rm
representation})\bigr)\biggr\} \\
&\ge&\|\psi_\gamma\|_{{\mathcal S}_{\bar{V}, 
{\mathcal 
A}}}. \end{eqnarray*} From above,
$\|\psi_\gamma\|_{{\mathcal S}_{\bar{V}, {\mathcal A}}}\ge
\|(\psi_\gamma)_{-\gamma}\|_{{\mathcal S}_{\bar{V},
{\mathcal A}}} =\|\psi\|_{{\mathcal S}_{\bar{V}, 
{\mathcal 
A}}}$, and hence $\|\psi_\gamma\|_{{\mathcal 
S}_{\bar{V}, {\mathcal A}}}=\|\psi\|_{{\mathcal 
S}_{\bar{V}, {\mathcal A}}}$. Thus ${\mathcal S}_{\bar{V},
{\mathcal A}}(\Gamma)$ is isometrically translation
invariant.

Let ${\mathcal S}$ be an isometrically translation
invariant Segal algebra in ${\mathcal 
A}_{\bar{V}}(\Gamma)$. We have to show 
${\mathcal 
S}_{\bar{V}, {\mathcal A}}(\Gamma)\subseteq {\mathcal S}$.
Let $\psi\in {\mathcal S}_{\bar{V}, {\mathcal A}}(\Gamma)$
with $\psi=\sum_{n=1}^{\infty}\psi_n
({\bar{V}}$-representation). For each $n$, there
is $\gamma_n\in \Gamma$ such that ${\rm
supp}((\psi_n)_{\gamma_n}) \subseteq \bar{V}$. By
Lemma 4.2\,(ii), $(\psi_n)_{\gamma_n}\in {\mathcal 
S}$, and we also have $\psi_n\in {\mathcal S}$ by the
translation invariance of ${\mathcal S}$, and hence
\begin{equation} \|\psi_n\|^{\bar{V}}
=\|(\psi_n)_{\gamma_n}\|^{\bar{V}} \ge
{C_2}^{-1}\|(\psi_n)_{\gamma_n}\|_{\mathcal S}
={C_2}^{-1}\|\psi_n\|_{\mathcal S}, \ n=1,2,3,...
\end{equation} for some constant $C_2>0$. (9)
implies that $\sum_{n=1}^{\infty}\psi_n$ is a
convergent series in ${\mathcal S}$, and has a limit
$\phi\in {\mathcal S}$. Then 
\begin{eqnarray}
\|\phi-\psi\|^{\bar{V}} &\le&
\|\phi-\sum_{n=1}^N\psi_n\|^{\bar{V}}+
\|\sum_{n=1}^N\psi_n-\psi\|^{\bar{V}} \nonumber \\
&\le& \|\phi-\sum_{n=1}^N\psi_n\|_{\mathcal S}+
\|\sum_{n=1}^N\psi_n-\psi\|^{\bar{V}}.
\end{eqnarray} Since the last part of (10)
converges to 0 as $N\rightarrow \infty$, we have
$\phi=\psi$. Therefore we have $\psi\in {\mathcal S}$.
\ \ \ \ $\Box$
\vspace{0.5cm}

{Remark 4.1.}\ Since the generalized Fourier
transform is an isometric isomorphism of 
${\mathcal 
L}^1_{\bar{V}}(G)$ onto ${\mathcal 
A}_{\bar{V}}(\Gamma)$, there is a one to one
correspondence between the family of all Segal
algebras in ${\mathcal L}^1_{\bar{V}}(G)$ and the
family of all Segal algebras in ${\mathcal 
A}_{\bar{V}}(\Gamma)$ through this transform.
Moreover an isometrically character invariant
Segal algebra in ${\mathcal L}^1_{\bar{V}}(G)$
corresponds to an isometrically translation
invariant Segal algebra in ${\mathcal 
A}_{\bar{V}}(\Gamma)$, and vice versa.

\vspace{0.5cm}

{\bf Definition 4.4.}\ By Theorem 4.3 and Remark
4.1, there exists the smallest isometrically
character invariant Segal algebra in ${\mathcal 
L}^1_{\bar{V}}(G)$, which will be denoted by
${\mathcal S}^1_{\bar{V} ,{\mathcal L}^1}(G)$. Obviously,
by the generalized Fourier transform, ${\mathcal 
S}^1_{ \bar{V} ,{\mathcal L}^1}(G)$ corresponds to
${\mathcal S}_{\bar{V} ,{\mathcal A}}(\Gamma)$.
\vspace{0.5cm}

\it

{\bf Corollary 4.4.}\ ${\mathcal 
S}^1_{\bar{V}}(G)={\mathcal S}^1_{ \bar{V} ,
{\mathcal 
L}^1}(G)$, and the identity map of ${\mathcal 
S}^1_{\bar{V}}(G)$ to ${\mathcal S}^1_{\bar{V} ,
{\mathcal 
L}^1}(G)$ is a bicontinuous algebra homomorphism.

\rm

Proof.\ The generalized Fourier transform maps
${\mathcal S}_{\bar{V}}^1(G)$ and ${\mathcal S}_{\bar{V},
{\mathcal L}^1}^1(G)$ to ${\mathcal S}_{\bar{V}}(\Gamma)$
and ${\mathcal S}_{\bar{V}, {\mathcal A}}(G)$,
respectively, and since ${\mathcal 
S}_{\bar{V}}(\Gamma)={\mathcal S}_{\bar{V}, {\mathcal 
A}}(\Gamma)$, we have ${\mathcal 
S}_{\bar{V}}^1(G)={\mathcal S}_{\bar{V}, {\mathcal 
L}^1}(G)$.
ince the identity map of ${\mathcal 
S}^1_{\bar{V}}(G)$ to ${\mathcal S}^1_{\bar{V}, 
{\mathcal L}^1}(G)$ is a norm decreasing isomorphism, the
corollary follows from the open mapping theorem.
\ \ $\Box$

\vspace{0.5cm}

For a commutative semisimple Banach algebra ${\bf
B}$, ${\mathcal I}({\bf B})$ denotes the family of all
closed ideals of ${\bf B}$. Then, by Theorem 4.3
(i), we have ${\mathcal I}({\mathcal S}_{\bar{V}}(\Gamma))
={\mathcal I}({\mathcal S}_{\bar{V}, {\mathcal A}}(\Gamma))$.
We will show that this fact and [8, Theorem B']
induce a natural one to one correspondence between
${\mathcal I}(A(\Gamma))$ and ${\mathcal I}({\mathcal 
A}_{\bar{V}}(\Gamma))$. To make the discussion
clear, for $I\in {\mathcal I}({\mathcal S}_{\bar{V}}(\Gamma))$, $\iota I$ denotes the same
$I$ considered as an ideal in 
${\mathcal I}({\mathcal S}_{\bar{V}, {\mathcal A}}(\Gamma))$.
\vspace{0.5cm} 

\it 

{\bf Theorem 4.5.}\ The map $I\mapsto
\overline{\iota (I\cap S_{\bar{V}}(\Gamma)})$ {\rm
(} $\bar{J}$ means the closure of $J$ in ${\mathcal 
A}_{\bar{V}}(\Gamma)${\rm )} is a one to one
correspondence between ${\mathcal I}(A(\Gamma))$ and
${\mathcal I}({\mathcal A}_{\bar{V}}( \Gamma))$.

\rm

Proof.\ With the help of [8, Theorem B'], we see
that all of the following mappings:
\begin{eqnarray*} 
&& I\mapsto I\cap {\mathcal 
S}_{\bar{V}}(\Gamma): {\mathcal 
I}(A(\Gamma))\rightarrow {\mathcal I}({\mathcal 
S}_{\bar{V}}(\Gamma)), 
\\ 
&& J\mapsto \iota{J}:
{\mathcal I}({\mathcal S}_{\bar{V}}(\Gamma)) \rightarrow
{\mathcal I}({\mathcal S}_{\bar{V}, {\mathcal A}}(\Gamma)), 
\\
&& J\mapsto \bar{J}: {\mathcal I}({\mathcal S}_{\bar{V}, {\mathcal A}}(\Gamma)) \rightarrow {\mathcal I}({\mathcal A}_{\bar{V}}(\Gamma)), 
\end{eqnarray*} 
are
bijective. Combining these three mappings, we
obtain the desired result.\ \ \ $\Box$

\vspace{0.5cm}

Remark 4.2.\ \, Let $V_1$ and $V_2$ be open
neighborhoods of $0\in \Gamma$ with compact
closures ${\bar{V}}_1$ and $\bar{V}_2$,
respectively. Since both ${\mathcal 
S}_{\bar{V}_1}(\Gamma)$ and ${\mathcal 
S}_{\bar{V}_2}(\Gamma)$ are the smallest
isometrically translation invariant Segal algebras
in $A(\Gamma)$, it follows that ${\mathcal 
S}_{\bar{V}_1}(\Gamma)={\mathcal 
S}_{\bar{V}_2}(\Gamma)$. From [8,Theorem 10.3],
\[C^{\bar{V}_1}_{BSE}(\Gamma)=\mathbb{M}({\mathcal 
S}_{\bar{V}_1}(\Gamma))=\mathbb{M}({\mathcal 
S}_{\bar{V}_2}(\Gamma))=C_{BSE}^{\bar{V}_2}(\Gamma),
\] 
and $\|\ \|^{\bar{V_1}}$ and
$\|\ \|^{\bar{V_2}}$ are equivalent by the
uniqueness of the complete norm topology in a
commutative semisimple algebra. From the
relations $\|\mu\|_{\bar{V}_j}
=\|\hat{\mu}\|^{\bar{V}_j}\ (\mu\in M(G)), j=1,2$,
it follows that $\|\cdot \|_{\bar{V}_1}$ and
$\|\cdot \|_{\bar{V}_2}$ defined in $M(G)$ are
equivalent, and hence ${\mathcal 
M}_{\bar{V}_1}(G)={\mathcal M}_{\bar{V}_2}(G)$. Then
${\mathcal L}^1_{\bar{V}_1}(G)={\mathcal 
L}^1_{\bar{V}_2}(G)$ follows from the definition
of ${\mathcal L}^1_{\bar{V}_j}(G), j=1,2$. Thus 
${\mathcal 
B}_{\bar{V}_1}(\Gamma)= {\mathcal 
B}_{\bar{V}_2}(\Gamma)$ and ${\mathcal 
A}_{\bar{V}_1}(\Gamma)={\mathcal 
A}_{\bar{V}_2}(\Gamma)$ also hold. So, we can use
notations \[ ({\mathcal L}^1(G), \|\cdot
\|_{\bar{V}}), ({\mathcal M}(G), \|\cdot
\|_{\bar{V}}), ({\mathcal A}(\Gamma), \|\cdot
\|^{\bar{V}}), \ {\rm and}\ ({\mathcal B}(\Gamma),
\|\cdot \|^{\bar{V}}), \] instead of the notations
\[ ({\mathcal L}^1_{\bar{V}}(G), \|\cdot
\|_{\bar{V}}), ({\mathcal M}_{\bar{V}}(G), \|\cdot
\|_{\bar{V}}), ({\mathcal A}_{\bar{V}}(\Gamma),
\|\cdot \|^{\bar{V}}),\ {\rm and}\ ({\mathcal 
B}_{\bar{V}}(\Gamma), \|\cdot \|^{\bar{V}}). \]
Nevertheless, we will continue to use the notation
$C^{\bar{V}}_{BSE}(\Gamma)$ to distinguish from
the notation $C_{BSE}(\Gamma)$ in [18, p.149].

\vspace{0.5cm}

{\bf \S 5.\ ${\mathcal B}(\Gamma)\subsetneqq
C_{BSE}^{\bar{V}}(\Gamma)$} \ This section is
devoted to prove that ${\mathcal B}(\Gamma)$ is
strictly contained in $C^{\bar{V}}_{BSE}(\Gamma)$.
 \vspace{0.5cm}

\it

{\bf Theorem 5.1.}\ ${\mathcal B} (\Gamma)\subsetneqq
C^{\bar{V}}_{BSE}(\Gamma). $

\vspace{0.3cm}

\rm

For the proof, we prepare necessary definition,
theorems and lemmas.

\vspace{0.5cm}

\rm

{\bf Definition 5.1.}\ Let $E$ be a subset of a
discrete abelian group $\Gamma$. $E$ is called a
Sidon set if there exists a constant $C_E>0$ such
that \[ \sum_{k=1}^n|c_k|\le C_E\|p\|_\infty\ \ \
\biggl(p=\sum_{k=1}^nc_k(x,\gamma_k)\ \in {\rm
span}(E) \biggr) \] holds {\rm ([17, \S 5.7])}.

\vspace{0.5cm}

\it

{\bf Theorem A} {\rm ([4])}.\ The Fourier
transform $L^1(G) \rightarrow C_0(\Gamma)$ is onto
only when $G$ is a finite group.

\vspace{0.5cm}

\it

{\bf Theorem B}\ {\rm ([17, Theorem 5.7.3]). A
subset $E$ of a discrete abelian group $\Gamma$ is
a Sidon set if and only if for any $\phi\in
C_0(E)$, there exists a function $f\in L^1(G)$
such that $\hat{f}(\gamma)=\phi(\gamma)\ \ (\gamma
\in E)$.

\vspace{0.5cm}

{\bf Theorem C} {\rm ([17, Theorem 5.7.4])}.\ A
subset $E$ of a discrete abelian group $\Gamma$ is
a Sidon set if there is a constant $\delta>0$ with
the following property: To every function $\psi$
on $E$ with $\psi(\gamma)=\pm 1$, there
corresponds a function $\xi\in B(\Gamma)$ such
that \[ \sup_{\gamma\in
E}|\psi(\gamma)-\xi(\gamma)|\le 1-\delta. \]

\it

{\bf Lemma 5.2.}\ If $\Gamma$ is an infinite
discrete abelian group, there exists a function
$\psi$ on $\Gamma$ such that $\psi(\gamma)=\pm 1$
and \begin{equation} \sup_{\gamma\in
\Gamma}|\psi(\gamma)-\xi(\gamma)|>1/2 \ \ (\xi\in
B(\Gamma)). \end{equation}

 \rm

Proof. $\Gamma$ is not a Sidon set by Theorems A
and B. So, by Theorem C, for any $\delta>0$, there
exists a function $\psi_\delta$ on $\Gamma$ such
that $\psi_\delta(\gamma)=\pm 1 (\gamma\in
\Gamma)$ and $\sup_{\gamma\in
\Gamma}|\psi_\delta(\gamma)-\xi(\gamma)|>1-\delta$
for every $\xi\in B(\Gamma)$. If we take
$\delta=1/2$ and write $\psi=\psi_{1/2}$, we
obtain (11).\ \ $\Box$ \ \

\vspace{0.5cm}

\it

{\bf Lemma 5.3.}\ Let $\mathbf{R}^d$ be the
$d$-dimensional real group with the dual group
$\hat{\mathbf{R}}^d$, and let 
$V=\{t=(t_1,...,t_d)\in \hat{\mathbf{R}}^d:
|t|=(\sum_{j=1}^d|t_j|^2)^{1/2}<1/4\}$. There
exists $\phi\in
C_{BSE}^{\bar{V}}(\hat{\mathbf{R}}^d)$ such that
$\sup_{t\in
\hat{\mathbf{R}}^d}|\phi(t)-\xi(t)|>1/2$ holds for
every $\xi\in B(\hat{\mathbf{R}}^d).$

\rm

Proof.\ By Lemma 5.2, there exists a function
$\psi$ on $\hat{\mathbf{Z}}^d$ such that
$\psi(n)=\pm 1\ $ and \begin{equation} \sup_{n\in
\hat{\mathbf{Z}}^d}|\psi(n)-\zeta(n)|>1/2\ \
(\zeta\in B(\hat{\mathbf{Z}}^d)). \end{equation}
Choose $h\in C_c(\hat{\mathbf{R}}^d)$ such that
${\rm supp}(h)\subset \{t\in \hat{\mathbf{R}}^d:
|t|< 1/8\} $ with $\|h^2\|_1=1$, and put
$e=h*h^*$. Then $e\in B(\hat{\mathbf{R}}^d)$ with
${\rm supp}(e)\subset V$ and
$e(0)=1$. Now, define a function $\phi$ on
$\hat{\mathbf{R}}^d$ by \[ \phi(t)=\sum_{j\in
\hat{\mathbf{Z}}^d}\psi (j)e(t-j)\ \ \ (t\in
\hat{\mathbf{R}}^d). \] Then $\phi(n)=\psi(n)\
(n\in \hat{\mathbf{Z}}^d)$, and $\phi\in
C_{BSE}^{\bar{V}}(\hat{\mathbf{R}}^d)$ because for
all $s\in \hat{\mathbf{R}}^d,\ (\bar{V}+s)\cap
{\rm supp}(e(t-j))\neq \emptyset$ at most one
$j\in \hat{\mathbf{Z}}^d$, and hence
$\|\phi\|^{\bar{V}} =\sup_{s\in
\hat{\mathbf{R}}^d}\|\phi\|^{\bar{V}, s}\le
\|e\|^{\bar{V}}<\infty$. From (12) and [17,
Theorem 2.7.2], for any $\xi\in
B(\hat{\mathbf{R}}^d)$, we have

\begin{eqnarray*} \sup_{t\in
\hat{\mathbf{R}}^d}|\phi(t)-\xi(t)| &\ge&
\sup_{n\in \hat{\mathbf{Z}}^d}|\phi(n)-\xi(n)| \\
&=&\sup_{n\in
\hat{\mathbf{Z}}^d}|\psi(n)-\zeta(n)|>1/2\
(\zeta=\xi|\hat{\mathbf{Z}}^d\in
B(\hat{\mathbf{Z}}^d)). \end{eqnarray*} This
completes the proof of the lemma.\ \ \ $\Box$

\vspace{0.5cm}

\it

{\bf Lemma 5.4}.\ Let $G_j$ be an LCA group with
the dual group $\Gamma_j$, and let $V_j$ be an
open neighborhood of $0_j\in \Gamma_j$ with
compact closure $\bar{V}_j, \ j=1,2.$ Put
$G=G_1\times G_2$. Then $\Gamma=\Gamma_1\times
\Gamma_2$ is the dual group of $G$, and
$V=V_1\times V_2$ is an open neighborhood of
$0=(0_1, 0_2)\in \Gamma$ with compact closure
$\bar{V}=\bar{V}_1\times \bar{V}_2$. Suppose
$\phi\in C_{BSE}^{\bar{V}_1}(\Gamma_1)$, and
define $\tilde{\phi}\in C_b(\Gamma)$ by \[
\tilde{\phi}(\gamma_1, \gamma_2)=\phi(\gamma_1)\ \
\ ((\gamma_1, \gamma_2)\in \Gamma). \] Then
$\tilde{\phi}\in C_{BSE}^{\bar{V}}(\Gamma)$ with
$\|\tilde{\phi}\|^{\bar{V}}\le
\|\phi\|^{\bar{V}_1}.$

\rm

Proof.\ Let $(\gamma_1, \gamma_2)\in \Gamma$ be
arbitrary, and let $\tilde{p}\in {\rm span}
(\bar{V}+(\gamma_1,\gamma_2))$ such that \[
\tilde{p}(x_1, x_2)= \sum_{k=1}^n c_k((x_1, x_2),
(\gamma_{1,k}+\gamma_1, \gamma_{2,k}+\gamma_2)) \
\ \ ((x_1, x_2)\in G). \] Define $p\in {\rm
span}(\bar{V}_1+\gamma_1)$ by
$p(x_1):=\sum_{k=1}^n c_k(x_1,
\gamma_{1,k}+\gamma_1)$, then \begin{eqnarray*}
\|p\|_\infty&=&\sup_{x_1\in G_1}\biggl|
\sum_{k=1}^nc_k(x_1,\gamma_{1,k}+\gamma_1)\biggl|
\\ &=&\sup_{(x_1, 0)\in G}\biggl|\sum_{k=1}^n
c_k((x_1, 0), (\gamma_{1,k}+\gamma_1,
\gamma_{2,k}+\gamma_2))
\biggr|\le\|\tilde{p}\|_\infty, \end{eqnarray*}
and \begin{eqnarray} \biggl|\sum_{k=1}^nc_k
\tilde{\phi}(\gamma_{1,k}+\gamma_1,\gamma_{2,k}
+\gamma_2)\biggr| &=&
\biggl|\sum_{k=1}^nc_k\phi(\gamma_{1,k}+\gamma_1)
\biggr| \nonumber \\ &\le&
\|p\|_\infty\|\phi\|^{\bar{V}_1,\gamma_1} \le
\|\tilde{p}\|_\infty\|\phi\|^{\bar{V}_1,
\gamma_1}. \end{eqnarray} (13) shows
$\tilde{\phi}\in C_{BSE}^{\bar{V}}(\Gamma)$ with
$\|\tilde{\phi}\|^{\bar{V}}\le
\|\phi\|^{\bar{V_1}}$. \ \ \ $\Box$

\vspace{0.5cm}

\rm 
 
{\bf Proof of Theorem 5.1.}\ By [6, Theorem
24.30], $G$ is topologically isomorphic with
$\mathbf{R}^d\times G_0$, where $d\ge 0$ and $G_0$
is an LCA group which has an open compact subgroup
$H$. We proceed with our proof under the
assumption that $G= \mathbf{R}^d\times G_0$.

First, we consider the case $d=0$. Then we have
$G=G_0$, which has an infinite compact open
subgroup $H$, and the dual group $\Gamma$ has an
open compact subgroup $H^{\perp} =\{\gamma\in
\Gamma: (x,\gamma)=1\ (x\in H)\}$. The dual group
of $H$ is $\Gamma/H^\perp$, which is an infinite
discrete group. We suppose Haar measures $m_H$ and
$m_{H^{\perp}}$ of $H$ and $H^\perp$,
respectively, are normalized so that
$m_H(H)=m_{H^{\perp}}(H^{\perp})=1$, and let
$V=H^\perp$.

By Lemma 5.2, there exists a function $\bar{\phi}$
on $\Gamma/H^{\perp}$ with $\bar{\phi}=\pm 1$ and

\begin{equation} \sup_{\gamma+H^{\perp}\in
\Gamma/H^{\perp}}|\bar{\phi}(\gamma+H^{\perp})
-\zeta(\gamma+H^{\perp})|>1/2\ \ \ (\zeta \in
B(\Gamma/H^\perp)). \end{equation} Put
$\phi=\bar{\phi}\circ \pi$, where $\pi$ denotes
the natural map $\gamma\mapsto \gamma+H^{\perp},
\Gamma\rightarrow \Gamma/H^{\perp}.$ Since $\phi$
is constant $\pm 1$ on each $\bar{V}+\gamma,
\gamma\in \Gamma$, we have $\|\phi\|^{\bar{V}}=1$,
hence $\phi\in C_{BSE}^{\bar{V}}(\Gamma)$. We will
show that the following is true: \begin{equation}
\sup_{\gamma\in
\Gamma}|\phi(\gamma)-\hat{\mu}(\gamma)|>1/2\ \
(\mu\in M(G)). \end{equation} Since
$\|\phi-\hat{\mu}\|^{\bar{V}}\ge \sup_{\gamma\in
\Gamma}|\phi(\gamma)-\hat{\mu}(\gamma)|, \mu\in
M(G)$, (15) implies $\phi\not\in {\mathcal 
B}(\Gamma)$.

To show (15), suppose contrary that there exists
$\mu\in M(G)$ such that \begin{equation}
\sup_{\gamma\in \Gamma}|
\phi(\gamma)-\hat{\mu}(\gamma)|\le 1/2.
\end{equation} Note that $\mu$ is concentrated on
some countable union of cosets, say
$\cup_{n=0}^\infty H+x_n$, where $x_n-x_m\not\in
H\ (m\neq n)$ and $x_0=0$. Put
$\mu_n=(\mu|H+x_n)*\delta_{-x_n}$, which is
concentrated on $H$ for $n\ge
0$, where $\delta_{-x_n}$ denotes the dirac
measure concentrated at $-x_n$. Then
$\mu=\sum_{n=0}^\infty\mu_n*\delta_{x_n}$ and
$\hat{\mu}(\gamma)=\sum_{n=0}^\infty
\widehat{\mu_n}(\gamma) (-x_n, \gamma)\ \ \
(\gamma\in \Gamma)$.

Let $\gamma_1\in \Gamma$ be arbitrary. Since
$1/2\ge |\phi(\gamma)-\hat{\mu}(\gamma)| \ \
(\gamma\in H^\perp+\gamma_1)$ from (16), it
follows that 
\begin{eqnarray}
 1/2&\ge&
\int_{H^\perp+\gamma_1}
|\phi(\gamma)-\hat{\mu}(\gamma)|d\gamma
=\int_{H^\perp}|\phi(\gamma_1+\gamma)-\hat{\mu}
(\gamma_1+\gamma)|d\gamma \nonumber \\ &=&
\int_{H^\perp}\biggl|\phi(\gamma_1+\gamma)
-\sum_{n=0}^\infty \widehat{\mu_n}
(\gamma_1+\gamma)(-x_n,\gamma_1+\gamma)
\biggr|d\gamma \nonumber \\ &\ge&
\biggl|\int_{H^\perp}\biggl(\phi(\gamma_1+\gamma)
-\sum_{n=0}^\infty\widehat{\mu_n}(\gamma_1+\gamma)
(-x_n, \gamma_1+\gamma)\biggr)d\gamma\biggr|
\nonumber \\ &=&
\biggl|\int_{H^\perp}\biggl(\phi(\gamma_1)-\sum_{n=0}^\infty
\widehat{\mu_n}(\gamma_1)(-x_n, \gamma_1)(-x_n,
\gamma)\biggr)d\gamma\biggl| \nonumber \\ &&
(\because \phi(\gamma_1+\gamma)=\phi(\gamma_1), \
\widehat{\mu_n}(\gamma_1+\gamma)=
\widehat{\mu_n}(\gamma_1)\ \ (\gamma\in H^\perp))
\nonumber \\ &=&
|\phi(\gamma_1)-\widehat{\mu_0}(\gamma_1)|
\biggl(\because\int_{H^\perp}(-x_n,\gamma)d\gamma=0\
(n\neq 0)\biggr). 
\end{eqnarray} 
Since $\gamma_1$
is arbitrarily chosen in $\Gamma$, we have from
(17) 
 \begin{equation}
|\phi(\gamma)-\widehat{\mu_0}(\gamma)|\le 1/2\ \ \
(\gamma\in \Gamma). \end{equation} Since
$\mu_0=\mu|H \in M(H)$ with $\widehat{\mu|H}\in
B(\Gamma/H^\perp)$, it follows from (18) that
\begin{equation}
|\bar{\phi}(\gamma+H^\perp)-\widehat{\mu|H}
(\gamma+H^{\perp})| \le 1/2\ \ \
(\gamma+H^\perp\in \Gamma/H^\perp). \end{equation}
(19) contradicts (14). Hence (15) is true, and so
${\mathcal B}(\Gamma)\subsetneqq
C^{\bar{V}}_{BSE}(\Gamma)$ holds.

Next, consider the case $G=\mathbf{R}^d\times G_0$
with $d\ge 1$ , and $G_0$ is an LCA group which
contains an open compact subgroup $H$. Then
$\Gamma=\hat{\mathbf{R}}^d \times \Gamma_0$, where
$\Gamma_0$ is the dual group of $G_0$ which
contains an open compact subgroup $H^\perp$. Let
$V=V_1\times H^{\perp}$, where 
$V_1=\{t\in \hat{\mathbf{R}}^d: |t|<1/4\}$. By
Lemma 5.3, there exists $\phi\in
C_{BSE}^{\bar{V}_1}(\hat{\mathbf{R}}^d)$ such that
\begin{equation}
 \sup_{t\in
\hat{\mathbf{R}}^d}|\phi(t)-\zeta(t)|>1/2 \ \
(\zeta\in B(\hat{\mathbf{R}}^d)). \end{equation}
Let $\tilde{\phi}(t, \gamma)=\phi(t) \ \
((t,\gamma)\in \hat{\mathbf{R}}^d\times
\Gamma_0)$. Then $\tilde{\phi} \in
C_{BSE}^{\bar{V}}(\Gamma)$ by Lemma 5.4. From (20)
and [17, Theorem 2.7.2], for any $\xi\in
B(\Gamma)$ we have \begin{eqnarray*}
\|\tilde{\phi}-\xi\|^{\bar{V}} &\ge&
\sup_{(t,\gamma)\in
\Gamma}|\tilde{\phi}(t,\gamma)-\xi(t,\gamma)| \\
&\ge& \sup_{t\in
\hat{\mathbf{R}}^d}|\tilde{\phi}(t,0)-\xi(t,0)| \\
&=& \sup_{t\in \hat{\mathbf{R}}^d}|\phi(t)
-\zeta(t)|>1/2\ \ \ (\zeta(t) =\xi((t, 0)) \in
B(\hat{\mathbf{R}}^d)). \end{eqnarray*} Hence
$\tilde{\phi}\not\in {\mathcal B}(\Gamma)$, that is,
${\mathcal B}(\Gamma)\subsetneqq
C_{BSE}^{\bar{V}}(\Gamma)$. \ \ $\Box$
\vspace{0.5cm}

\rm 

{\bf \S6.\ Segal algebras in ${\mathcal A}(\Gamma)$,
II} \  In this section, we
investigate Segal algebras in ${\mathcal A}(\Gamma)$
which were not covered in $\S 4$. 
 \vspace{0,5cm}
\it

{\bf Theorem 6.1.} Let $({\mathcal S}, \|\cdot
\|_{\mathcal S})$ be a Segal algebra in 
$({\mathcal 
A}(\Gamma), \|\cdot \|^{\bar{V}})$. Then we have

{\rm (i)}\ $\mathbb{M}({\mathcal S})\supseteq
C_{BSE}^{\bar{V}}(\Gamma)$.

{\rm (ii)}\ If ${\mathcal S}$ is isometrically
translation invariant, we have $\mathbb{M}({\mathcal 
S}) =C^{\bar{V}}_{BSE}(\Gamma)$.

\rm Proof.\ (i)\ Suppose $\phi\in
C_{BSE}^{\bar{V}}(\Gamma)$ and $\psi\in {\mathcal S}$.
There exist $\eta\in {\mathcal A}(\Gamma)$ and
$\tau\in {\mathcal S}$ such that $\psi=\eta\tau$\,([8,
Theorem A']). Since $\phi\in C_b(\Gamma)$ with
$\phi\psi=(\phi\eta)\tau\in {\mathcal S}$, and that
$\phi(\psi_1\psi_2)=(\phi\psi_1)\psi_2\ \ (\psi_1,
\psi_2\in {\mathcal S})$, we have $\phi\in
\mathbb{M}({\mathcal S})$.

(ii)\ We have only to show $\mathbb{M}({\mathcal S})
\subseteq C^{\bar{V}}_{BSE}(\Gamma)$. Let $\phi\in
\mathbb{M}({\mathcal S})$, and choose $u\in {\mathcal 
S}_c$ such that $u(\gamma')=1\ \ (\gamma'\in
\bar{V})$. For each $\gamma\in \Gamma$ it follows
that \[ \|\phi\|^{\bar{V}, \gamma}\le \|\phi
u_\gamma\|^{\bar{V}} \le \|\phi u_\gamma\|_{\mathcal 
S} \le \|\phi\|_{\mathbb{M}({\mathcal 
S})}\|u_\gamma\|_{\mathcal S}
=\|\phi\|_{\mathbb{M}({\mathcal S})}\|u\|_{\mathcal S}. \]
Hence we have $\phi\in C^{\bar{V}}_{BSE}(\Gamma)$,
that is, $\mathbb{M}({\mathcal S}) \subseteq
C^{\bar{V}}_{BSE}(\Gamma)$. \ \ \ $\Box$

\vspace{0.5cm}

{\bf Definition 6.1.}\ Let $\mu$ be a positive
unbounded Radon measure on $\Gamma$ and $1\le
p<\infty$. Define ${\mathcal A}_{p, \mu}(\Gamma)$ and
$\|\cdot \|_{{\mathcal A}_{p, \mu}}$ by \[ {\mathcal 
A}_{p,\mu}(\Gamma):=\bigl\{ \psi\in {\mathcal 
A}(\Gamma): \psi\in L^p(\mu)\bigr\},\ \
\|\psi\|_{{\mathcal A}_{p, \mu}}:=\|\psi\|^{\bar{V}}
+\|\psi\|_{L^p(\mu)}\ \ (\psi\in {\mathcal A}_{p,
\mu}(\Gamma)). 
\]

\vspace{0.5cm}

\it

{\bf Theorem 6.2.}\ $({\mathcal A}_{p, \mu}(\Gamma),
\|\cdot \|_{{\mathcal A}_{p, \mu}})$ is a Segal
algebra in ${\mathcal A}(\Gamma)$ {\rm(cf. [14, $\S
5$, Examples (iv)]). \rm

\rm

Proof.\ That $({\mathcal A}_{p, \mu}(\Gamma), \|\cdot
\|_{{\mathcal A}_{p, \mu}})$ forms a normed linear
space is apparent. To show that $\|\cdot \|_{{\mathcal 
A}_{p,\mu}}$ is a complete norm, let
$\{\psi_n\}_{n=1}^\infty$ be a sequence in 
${\mathcal 
A}_{p,\mu}(\Gamma)$ such that $\sum_{n=1}^\infty
\|\psi_n\|_{{\mathcal A}_{p,\mu}}<\infty$. Since
$\sum_{n=1}^\infty \psi_n$ is a $\|\cdot
\|^{\bar{V}}$-convergent series, there exists
$\psi\in {\mathcal A}(\Gamma)$ such that
$\psi=\sum_{n=1}^\infty \psi_n$. Then
$\|\psi\|_{L^p(\mu)}\le
\sum_{n=1}^\infty\|\psi_n\|_{L^p(\mu)}<\infty$.
Therefore $\psi\in {\mathcal A}(\Gamma)\cap L^p(\mu)$
with \[
\bigl\|\psi-\sum_{n=1}^N\psi_n\bigr\|_{{\mathcal 
A}_{p,\mu}}\le \sum_{n=N+1}^\infty
\|\psi_n\|_{{\mathcal A}_{p,\mu}}\rightarrow
0(N\rightarrow \infty). \] Hence, $\|\cdot
\|_{{\mathcal A}_{p, \mu}}$ is a complete norm which
dominates $\|\cdot \|^{\bar{V}}$.

For $\psi\in {\mathcal A}_{p,\mu}(\Gamma)$ and
$\phi\in {\mathcal A}(\Gamma)$, we have 
\[
\|\phi\psi\|_{{\mathcal A}_{p,
\mu}}=\|\phi\psi\|^{\bar{V}}
+\|\phi\psi\|_{L^p(\mu)} \le
\|\phi\|^{\bar{V}}\|\psi\|^{\bar{V}}
+\|\phi\|^{\bar{V}}\|\psi\|_{L^p(\mu)}
=\|\phi\|^{\bar{V}}\|\psi\|_{{\mathcal A}_{p, \mu}},
\] and hence ${\mathcal 
A}_{p,\mu}(\Gamma)$ is a Banach ideal in
 ${\mathcal A}(\Gamma)$. It is dense in 
 ${\mathcal A}(\Gamma)$
since it contains ${\mathcal A}(\Gamma)_c$.

Finally, we show that it possesses approximate
units. For each $\Omega\in {\mathcal K}(\Gamma)$, let
$u_\Omega \in A(\Gamma)_c$ be such that
$\|u_\Omega\|_A \le 2$ and $u_\Omega(\gamma)=1\ \
(\gamma\in \Omega)$. Suppose that $\psi\in {\mathcal 
A}_{p, \mu}(\Gamma)$ and $\varepsilon>0$ are
given. We can choose $\Omega_1, \Omega_2\in {\mathcal 
K}(\Gamma)$ such that \[ \|\psi-u_{\Omega}\psi
\|^{\bar{V}}\le \varepsilon/2 \ \
(\Omega_1\subseteq \Omega) , \
\|\psi-u_{\Omega}\psi\|_{L^p(\mu)}\le
\varepsilon/2 \ \ (\Omega_2\subseteq \Omega). 
\]
It follows that { \[ \|\psi-u_\Omega\psi\|_{{\mathcal
A}_{p, \mu}} = \|\psi-u_\Omega\psi\|^{\bar{V}}+
\|\psi-u_\Omega\psi\|_{L^p(\mu)} \le \varepsilon \
\ (\Omega_1\cup \Omega_2\subset \Omega). \] Thus,
${\mathcal A}_{p,\mu}(\Gamma)$ possesses approximate
units, and hence it is a Segal algebra in ${\mathcal
A}(\Gamma)$. \ \ \ $\Box$ \vspace{0.5cm}

Let $\textfrak{S}({\mathcal A}(\Gamma))$ denote the
family of all the Segal algebras in ${\mathcal 
A}(\Gamma)$. Then we have the following:

\vspace{0.5cm}

{\bf Proposition 6.3.}\ $\bigcap \bigl\{ {\mathcal S}:
{\mathcal S} \in \textfrak{S} ({\mathcal 
A}(\Gamma))\bigr\}={\mathcal A}(\Gamma)_c.$

\rm

Proof.\ Let $\psi\in {\mathcal A}(\Gamma)\setminus
{\mathcal A}(\Gamma)_c$. We can choose a positive
unbounded Radon measure $\mu$ on $\Gamma$ such
that $\psi\not\in L^1(\mu)$. By Theorem 6.2, the
Segal algebra ${\mathcal A}_{1, \mu}(\Gamma)$ does not
contain $\psi$, and the inclusion $\subseteq$
holds. The reverse inclusion follows from [8,
Lemma 3.4]. } \ \ $\Box$

\vspace{0.5cm}

{\bf Definition 6.2.}\ Let $1\le p<\infty$, and
define ${\mathcal A}_p(\Gamma)$ and $\|\cdot \|_{{\mathcal 
A}_p}$ by \[{\mathcal A}_p(\Gamma):= \bigl\{\psi\in
{\mathcal A}(\Gamma): \psi\in L^p(\Gamma)\bigr\}, \ \
\|\psi\|_{{\mathcal A}_p}: =\|\psi\|^{\bar{V}}
+\|\psi\|_p \ \ (\psi\in {\mathcal A}_p(\Gamma)). 
\]

\it

{\bf Theorem 6.4.} \ For $1\le p<q<\infty$, we
have

{\rm (i)}\ $({\mathcal A}_p(\Gamma), \|\cdot 
\|_{{\mathcal A}_p} )$ is an isometrically translation invariant Segal algebra in ${\mathcal A}(\Gamma)$.

{\rm (ii)}\ ${\mathcal A}_p(\Gamma)\subsetneqq {\mathcal A}_q(\Gamma)$.
\rm

Proof.\ (i)\ That $({\mathcal A}_p(\Gamma),
\|\cdot\|_{{\mathcal A}_p}) {\mathcal A}(\Gamma)$ follows from Theorem 6.2 as a spacial case for $\mu=d\gamma$.
 The isometric translation invariance of 
 $\|\cdot \|_{{\mathcal A}_p}$ is apparent.

(ii)\ Let $u$ be in $A(\Gamma)_c$ such that
$u(\gamma)=1\ \ (\gamma\in \bar{V})$. Choose
$\{\gamma_n\}_{n=1}^\infty \subset \Gamma$ such
that \begin{equation} ({\rm
supp}(u_{\gamma_m})-\bar{V} )\cap ({\rm
supp}(u_{\gamma_n})-\bar{V}) =\emptyset\ \ (m\neq
n), \end{equation} and put
$\psi(\gamma)=\sum_{n=1}^\infty
\frac{1}{n^{1/p}}u(\gamma-\gamma_n) \ \ (\gamma\in
\Gamma)$. From (21), we see that for any
$\gamma\in \Gamma$, $(\bar{V}+\gamma)\cap {\rm
supp}(u_{\gamma_n})\neq \emptyset$ holds at most
one $n\in \mathbf{N}$. One can easily see from
this that $\psi\in C^{\bar{V}}_{BSE}(\Gamma)$.

Let $\varepsilon>0$ be given. Choose $n_0\in
\mathbf{N}$ and $\Omega\in {\mathcal K}(\Gamma)$ such
that \[ \frac{1}{{(n_0+1)}^{1/p}} \|u\|^{\bar{V}}
\le \varepsilon,\ \cup_{j=1}^{n_0}({\rm
supp}(u_{\gamma_j})-\bar{V})\subset \Omega. \] If
$\gamma\in \Gamma\setminus \Omega$, then
$(\bar{V}+\gamma)\cap {\rm supp}(u_{\gamma_n})
=\emptyset\ (1\le n\le n_0)$, hence \[
\|\psi\|^{\bar{V}, \gamma}\le \sup_{n\ge
n_0+1}\frac{1}{n^{1/p}} \|u_{\gamma_n}\|^{\bar{V},
\gamma}
\le\frac{1}{{(n_0+1)}^{1/p}}\|u\|^{\bar{V}} \le
\varepsilon, \] and we have $\psi\in {\mathcal 
A}(\Gamma)$. Also we have from (21) that
 \[
\|\psi\|_p^p=\sum_{n=1}^\infty
\frac{1}{n}\int_\Gamma|u(\gamma-\gamma_n)|^pd\gamma
=\sum_{n=1}^\infty\frac{1}{n} \biggl(\int_\Gamma
|u(\gamma)|^pd\gamma\biggr)=\infty, \] \[
\|\psi\|_q^q=\sum_{n=1}^\infty
\frac{1}{n^{q/p}}\int_\Gamma|u(\gamma-\gamma_n)|^q
d\gamma=\sum_{n=1}^\infty\frac{1}{n^{q/p}}
\biggl(\int_\Gamma
|u(\gamma)|^qd\gamma\biggr)<\infty. \] Hence
$\psi\in {\mathcal A}_q(\Gamma)\setminus {\mathcal 
A}_p(\Gamma)$. \ \ \ $\Box$

\vspace{0.5cm}

{\bf Definition 6.3.}\ Let $1\le p<\infty$, and
put \[ {\mathcal L}^1_p(G):=\{f\in {\mathcal L}^1(G):
\hat{f}\in L^p(\Gamma)\}, \ \|f\|_{{\mathcal L}^1_p} :
=\|f\|_{\bar{V}}+ \|\hat{f}\|_p\ \ (f\in {\mathcal 
L}^1_p(G)). 
\]

\it 

{\bf Corollary\,6.5.}\ {\rm (i)}\ Let $1\le
p<\infty$. Then $({\mathcal L}^1_p(G), \|\cdot
\|_{{\mathcal L}^1_p})$ is an isometrically character
invariant Segal algebra in ${\mathcal L}^1(G)$.

{\rm (ii)}\ If $1\le p<q<\infty$, we have ${\mathcal 
L}^1_p(G) \subsetneqq {\mathcal L}^1_q(G)$.

\rm

Proof. The proof readily follows from Definition
6.3 and Theorem 6.4.} \ \ $\Box$
 
\vspace{0.5cm} 

\rm

{Remark 6.1.\ (i)\ Let $G$ be a non-compact,
non-discrete LCA group, and $1\le p\le 2$. In [11,
\S 6], $A_p(G):=\{f\in L^1(G): \hat{f}\in
L^p(\Gamma)\}$ with norm
$\|f\|_{A_p}=\|f\|_1+\|\hat{f}\|_p\ (f\in
A_p(G))$. $A_p(G)$ is a Segal algebra such that
$\mathbb{M}(A_p(G))=M(G)$ (cf. [14, \S 5, Examples
(vi)], [11 Theorem 6.3.1]).

Define $A_p(\Gamma)=\{\psi\in A(\Gamma): \psi\in
L^p(\Gamma)\}$ with norm $\|\psi\|_{A_p(\Gamma)}
=\|\psi\|_A+\|\psi\|_p\ (\psi\in A_p(\Gamma))$.
Apparently, $A_p(\Gamma)$ is a Segal algebra in
$A(\Gamma)$ since the Fourier transform is an
isometric algebra homomorphism of $A_p(G)$ onto
$A_p(\Gamma)$, and so
$\mathbb{M}(A_p(\Gamma))=B(\Gamma)$.

On the other hand, $\mathbb{M}({\mathcal A}_p(\Gamma))
=C^{\bar{V}}_{BSE}(\Gamma)$ from Theorems 6.1 and
6.4. Since $B(\Gamma) \subsetneqq
C^{\bar{V}}_{BSE}(\Gamma)$ we see that
$A_p(\Gamma)\neq {\mathcal A}_p(\Gamma)$, hence
$A_p(\Gamma)\subsetneqq {\mathcal A}_p(\Gamma)$. This
means at once that $A_p(G)\subsetneqq 
{\mathcal 
L}^1_p(G)$.

(ii)\ In the case that $G$ is a non-discrete,
compact abelian group and $1\le p\le 2$, we have
$A_p(G)={\mathcal L}^1_p(G)$. To see $\supseteq$,
suppose $f\in {\mathcal L}^1_p(G)$. Since $\hat{f}\in
L^p(\Gamma)\subset \mathscr{I}(\Gamma)$ by [1,
Theorem 3.5] we have \ $\check{\hat{f}}\in
L^q(G)\subset L^1(G)$, where $p+q=pq$, so we can
conclude that $f=\check{\hat{f}}\in L^1(G)$, that
is $f\in A_p(G)$. The inclusion $\subseteq$ is
trivial since $L^1(G)\subset {\mathcal L}^1(G)$.

\vspace{0.5cm}

{\bf \S 7.\,Locally bounded elements of 
${\mathcal 
M}(G)$, and transformable Radon measures}
$\textfrak{M}_T(G)$ \ For further development of
our study, it is necessary to use transformable
Radon measures and their Fourier transforms which
were introduced in [1] by L. Argabright and J. Gil
de Lamadrid.

In this section, we suppose that $G$ is
$\sigma$-compact, and let $\{K_n: n=1,2,3,....\}$
be a sequence of elements of ${\mathcal K}(G)$ such
that $K_n\subset K_{n+1}, n=1,2,3,...,$ and $
G=\cup_{n=1}^\infty K_n$.

\vspace{0.5cm}

{\bf Definition 7.1}\,([1, p.8 and p.21]).\, A
Radon measure $\mu\in \textfrak{M}(G)$ is said to
be transformable, expressed by $\mu\in
\textfrak{M}_T(G)$, if there exists a Radon
measure $\hat{\mu}\in \textfrak{M}(\Gamma)$
satisfying \begin{equation} \int_G h*h^*(x)d\mu(x)
=\int_\Gamma
|\check{h}(\gamma)|^2d\hat{\mu}(\gamma) \ \ \
(h\in C_c(G)). 
\end{equation} $\hat{\mu}$ in (22)
is called the Fourier transform of $\mu$. One can
see easily that (22) is equivalent to
\begin{equation} \int_G k(x)d\mu(x)= \int_\Gamma
\check{k}(\gamma) d\hat{\mu}(\gamma)\ \ \ (k\in
C_{c,2}(G)), \end{equation} where $C_{c,2}(G)={\rm
span}\{h*h^*: h\in C_c(G)\} ={\rm span}\{h*k: h,
k\in C_c(G)\}$. We also use the notation
$\mathscr{I}(G):=\{\mu\in \textfrak{M}_T(G):
\hat{\mu}\in \textfrak{M}_T(\Gamma)\}$. For every
$\mu\in \mathscr{I}(G)$, the inversion formula
$\mu=\check{\hat{\mu}} $ holds ([1, Theorem 3.4]).

\vspace{0.5cm}

\it

{\bf Lemma 7.1.}\ Let $\mu\in \textfrak{M}_T(G)$
with the Fourier transform
$\hat{\mu}=\phi(\gamma)d\gamma \in
\textfrak{M}(\Gamma)$ for some
$\phi\in C_b(\Gamma)$. Then we have
\begin{equation} \int_G k(x)d\mu(x)=\int_\Gamma
\check{k}(\gamma)\phi(\gamma)d\gamma\ \ (k\in
C_c(G)\cap B(G)). \end{equation} \rm

Proof.\ Let $k\in C_c(G)\cap B(\Gamma)$ and
$\varepsilon>0$. We can choose $u\in C_c(G)$ such
that \begin{eqnarray} &&\biggl|\int_G
k(x)d\mu(x)-\int_G
k*u(x)d\mu(x)\biggr|\le
\varepsilon,\ \ \\ &&\biggl|\int_\Gamma
\check{k}(\gamma)
\phi(\gamma)d\gamma-\int_\Gamma
\check{k}(\gamma)\check{u}(\gamma)\phi(\gamma)d\gamma\biggr|\le
\varepsilon. \end{eqnarray} Since
$k*u\in C_{c,2}(G)$, we have
from (23), (25), and (26) that \[ \biggl|\int_G
k(x)dx-\int_\Gamma
\check{k}(\gamma)\phi(\gamma)d\gamma\biggr|\le
2\varepsilon. \] Since $\varepsilon>0$ in the
above is arbitrary, we obtain (24). \ \ $\Box$

\vspace{0,5cm}

{\bf Definition 7.2.}\ An element $\mu\in 
{\mathcal 
M}(G)$ is said to be locally bounded if $\mu$ has
at least one "locally bounded" representation: \[
\mu=\lim_{n\rightarrow\infty}\mu_n \ (\mu_n\in
M(G), n=1,2,3,...)\ {\rm in}\ ({\mathcal M}(G),
\|\cdot\|_{\bar{V}}) \] such that
$\{h\mu_n\}_{n=1}^\infty$ is a Cauchy sequence in
$(M(G), \|\cdot\|)$ for every $h\in C_c(G)$. We
write ${\mathcal M}(G)_{(0)}:=\bigl\{\mu\in
 {\mathcal 
M}(G): \mu\ {\rm is \ locally \ bounded}\bigr\}$.

 Note that when $G$ is compact
and $\mu\in {\mathcal M}(G)_{(0)}$ with
$\displaystyle\mu=\lim_{n\rightarrow\infty}\mu_n$,
a locally bounded representation, we have
$\displaystyle 1_G\mu=\lim_{n\rightarrow\infty}
1_G \mu_n=\mu\in M(G)$, hence ${\mathcal 
M}(G)_{(0)}=M(G)$.

\it

\vspace{0.5cm}

\it {\bf Lemma 7.2.}\, If $\mu\in {\mathcal 
M}(G)_{(0)}$ has two locally bounded
representations\ $\displaystyle\mu=
\lim_{n\rightarrow\infty}\mu_{j,n}, j=1,2$, we
have $\|h\mu_{1,n}-h\mu_{2,n}\|\rightarrow 0\,
(n\rightarrow\infty) \ \ \ (h\in C_c(G)). $

\rm

Proof.\ Suppose $h\in C_c(G)$ with $K={\rm
supp}(h)\in {\mathcal K}(G)$. Let $u_K\in L^1(\Gamma)$
be such that $\check{u_K}\in C_c(G)$ with
$\check{u_K}(x)=1\ (x\in K)$. Since
$\check{u_K}\in C_c(G)\cap B(G)$, we have from
Lemma 3.3 that \[ \lim_{n\rightarrow\infty}
\widehat{\check{u_K}\mu_{j,
n}}(\gamma)=\lim_{n\rightarrow\infty}u_K*\widehat{\mu_{j,n}}(\gamma)=u_K*\hat{\mu}(\gamma)
\ \ \ (\gamma\in \Gamma,\ j=1,2), 
\]
 it follows
that $\lim_n
\check{u_K}\mu_{1,n}=\lim_n\check{u_K}\mu_{2,n}$,
and so
$\displaystyle\lim_n\|\check{u_K}\mu_{1,n}-\check{u_K}\mu_{2,n}\|=0.$
Hence we have \begin{eqnarray*}
\|h\mu_{1,n}-h\mu_{2,n}\|&=&\|h\check{u_K}\mu_{1,n}-
h\check{u_K}\mu_{2,n}\| \\ &\le&
\|h\|_\infty\|\check{u_K}\mu_{1,n}-\check{u_K}\mu_{2,n}\|\rightarrow
0 (n\rightarrow\infty). \ \ \ \ \ \ \Box
\end{eqnarray*}

\vspace{0.5cm}

{\bf Definition 7.3.}\ Suppose 
$\mu\in {\mathcal M}(G)_{(0)}$ with a locally bounded representation
\[ (*)\ \ \mu=\lim_{n\rightarrow\infty} \mu_n. 
\]
For each $h\in C_c(G)$, we define $h\mu\in M(G)$
by $\displaystyle
h\mu=\lim_{n\rightarrow\infty}h\mu_n$, which is
independent of the choice of the representation
$(*)$ with the help of Lemma 7.2.

\vspace{0.5cm}

\it

{\bf Lemma 7.3.}\ $\|\rho*\phi\|^{\bar{V}} \le
\|\phi\|^{\bar{V}}\|\rho\|\ \ \ (\phi\in
C_{BSE}^{\bar{V}}(\Gamma), \rho\in M(\Gamma)).$

\rm

Proof.\ Let $\phi\in C_{BSE}^{\bar{V}}(\Gamma)$
and $\rho\in M(\Gamma)$. Then \begin{eqnarray*}
\|\phi*\rho\|^{\bar{V}}&=& \sup_{\gamma\in
\Gamma}\|\phi*\rho\|^{\bar{V},\gamma} =
\sup_{\gamma\in \Gamma}\sup_{p\in {\rm
span}(\bar{V}+\gamma), \|p\|_\infty\le 1}
\biggl|\sum_{\gamma'\in
\bar{V}+\gamma}\hat{p}(\gamma')\phi*\rho(\gamma')\biggr|
\\ &=& \sup_{\gamma\in \Gamma}\sup_{p\in {\rm
span}(\bar{V}+\gamma), \|p\|_\infty\le 1}
\biggl|\sum_{\gamma'\in
\bar{V}+\gamma}\hat{p}(\gamma')\int_\Gamma
\phi(\gamma'-\gamma'')d\rho(\gamma'')\biggr| \\
&=& \sup_{\gamma\in \Gamma}\sup_{p\in {\rm
span}(\bar{V}+\gamma), \|p\|_\infty\le 1}
\biggl|\int_\Gamma \sum_{\gamma'\in
\bar{V}+\gamma}\hat{p}(\gamma')\phi(\gamma'-\gamma'')d\rho(\gamma'')\biggr|
\\ &\le& \int_\Gamma \sup_{\gamma\in
\Gamma}\sup_{p\in {\rm span}(\bar{V}+\gamma),
\|p\|_\infty\le 1} \biggl|\sum_{\gamma'\in
\bar{V}+\gamma}\hat{p}(\gamma')\phi(\gamma'-\gamma'')\biggr|
d|\rho|(\gamma'') \\ &\le&
\|\phi\|^{\bar{V}}\|\rho\|.\ \
\ \ \Box \end{eqnarray*}

\vspace{0,5cm}

{\bf Lemma 7.4.}\ Suppose $\mu\in {\mathcal M}(G)$. If
there exists $\iota\mu\in \textfrak{M}_T(G)$ such
that $d\widehat{\iota\mu} (\gamma)
=\hat{\mu}(\gamma)d\gamma$, we have \[
\hat{k}*\hat{\mu}\in B(\Gamma) \ \ (k\in
C_c(G)\cap B(G)). \] \rm

Proof.\ For $\gamma\in \Gamma$ and $k\in C_c
(G)\cap B(G)$, we have $(-x, \gamma)k(x)\in
C_c(G)\cap B(G)$, and from Lemma 7.1, it follows
that \[ \widehat{k\iota\mu}(\gamma)= \int_G (-x,
\gamma)k(x)d\iota \mu(x)= \int_\Gamma
\check{k}(\gamma'-\gamma)
\hat{\mu}(\gamma')d\gamma'
=\hat{k}*\hat{\mu}(\gamma), \] which implies
$\hat{k}*\hat{\mu}=\widehat{k\iota \mu}\in
B(\Gamma)$. \ \ \ $\Box$

\vspace{0.5cm}

{\bf Lemma 7.5.}\ Suppose $\mu\in {\mathcal M}(G)$. If
$\hat{k}*\hat{\mu}\in B(\Gamma)\ \ (k\in
C_c(G)\cap B(G))$, we have $\mu\in {\mathcal 
M}(G)_{(0)}$.

\rm

Proof.\ \ For each $k\in C_c(G)\cap B(G)$, we have
$\hat{k}*\hat{\mu}\in B(\Gamma)$ by the
assumption, and there exist $\xi\in M(G)$ such
that $\hat{\xi}=\hat{k}*\hat{\mu}$. We denote this
$\xi\in M(G)$ by $k\mu$. For $h, k\in C_c(G)\cap
B(G)$, we have
$\widehat{h(k\mu)}=\widehat{(hk)\mu}
=\hat{h}*\hat{k}*\hat{\mu}$ from Lemma 3.3 and the
definition of $h\mu, (hk)\mu$, it follows that
$h(k\mu)=(hk)\mu$.
 
Let
$\displaystyle\mu=\lim_{n\rightarrow\infty}\mu_n\
(\mu_n\in M(G), n=1,2,3,...)$ with $\|\mu-\mu_n\|_
{\bar{V}}\rightarrow 0(n\rightarrow\infty)$. We
can suppose that each $\mu_n \ (n=1,2,3,...)$ has
compact support. Let $\varepsilon>0$ be given
arbitrarily. Choose $n_1\in \mathbf{N}$ such that
\begin{equation} \|\mu-\mu_n\|_{\bar{V}}\le
\varepsilon/3\ \ (n_1\le n). \end{equation} For
each $n\in \mathbf{N}$, let $u_{K_n}\in
L^1(\Gamma)$ such that $\|u_{K_n}\|_1\le 2,
\check{u_{K_n}}\in C_c(G)$ with
$\check{u_{K_n}}(x)=1\ \ (x\in K_n)$. Choose
$n_2\in \mathbf{N}$ such that $n_1\le n_2$ and
${\rm supp}(\mu_{n_1})\subseteq K_{n_2}$. If
$n_2\le n$, we have from (27) and
Lemma 7.3 that
\begin{eqnarray} \|\check{u_{K_n}}\mu_{n_1}
-\check{u_{K_n}}\mu\|_{\bar{V}}
&=&\|u_{K_n}*\widehat{\mu_{n_1}}-u_{K_n}*\hat{\mu}\|^{\bar{V}}=\|u_{K_n}*(\widehat{\mu_{n_1}}-\hat{\mu})\|^{\bar{V}}
\nonumber \\ &\le&
\|u_{K_n}\|_1\|\widehat{\mu_{n_1}}-\hat{\mu}\|^{\bar{V}}
\le 2\|\mu_{n_1}-\mu\|_{\bar{V}}\le
(2/3)\varepsilon. \end{eqnarray} Since
$\check{u_{K_n}}\mu_{n_1}=\mu_{n_1}\ \ (n_2\le
n)$, we have from (27) and (28) that
\begin{eqnarray*} \|\check{u_{K_n}}
\mu-\mu\|_{\bar{V}} &\le&
\|\check{u_{K_n}}\mu-\check{u_{K_n}}\mu_{n_1}\|_{\bar{V}}+
\|\check{u_{K_n}}\mu_{n_1}-\mu_{n_1}\|_{\bar{V}}+
\|\mu_{n_1}-\mu\|_{\bar{V}} \nonumber \\ &\le&
(2/3)\varepsilon+\varepsilon/3=\varepsilon\ \ \
(n_2\le n). \end{eqnarray*}

Let $h\in C_c(G)$ with $K:={\rm supp}(h)$. Choose
$n_0\in \mathbf{N}$ such that $K\subset K_{n_0}$.
Then we have \begin{equation}
k(\check{u_{K_n}}\mu)=(k\check{u_{K_{n_0}}})\mu=k\mu\
\ (n_0\le n)\ \ (k\in C^K_c(G)\cap B(G)).
\end{equation} Since $C^K_c(G)\cap B(G)$ is dense
in $C^K_c(G)$, we have from (29)
\[h(\check{u_{K_n}}\mu)=h(\check{u_{K_{n_0}}}\mu)
\ (n_0\le n). \] Therefore
$\mu=\lim_{n\rightarrow\infty} \check{u_{K_n}}\mu$
is a locally bounded representation; $\mu\in {\mathcal 
M}(G)_{(0)}$. \ \ $\Box$

{\bf Lemma 7.6.}\ For $\mu\in {\mathcal M}(G)_{(0)}$,
there exists $\iota\mu\in \textfrak{M}_T(G)$ such
that $d\widehat{\iota\mu}(\gamma)=
\hat{\mu}(\gamma)d\gamma, $ and $h\mu=h\iota \mu$
 holds for all $h\in C_c(G)$.

Note that, this $\iota \mu$ is uniquely determined
from [1, Theorem 2.1].

 \rm

Proof.  When $G$ is compact, we
have ${\mathcal M}_{(0)}(G)=M(G)$. So, the lemma is
trivially true with $\mu=\iota \mu$.

We consider, when $G$ is non-compact. Let $
\mu\in {\mathcal M}(G)_{(0)}$ with a locally bounded
representation
$\displaystyle\mu=\lim_{n\rightarrow\infty}\mu_n$.
We see easily that \begin{equation} (hk)\mu=\lim_n
(hk)\mu_n=h\lim_n k\mu_n=h(k\mu) \ \ (h, k\in
C_c(G)). 
\end{equation}
 Further, for any $K\in
{\mathcal K}(G)$, we have 
\[
 \biggl|\int_G d
h\mu(x)\biggr|=\biggl|\int_G d
(h\check{u_K})\mu(x)\biggr| =\biggl|\int_G
h(x)d(\check{u_K}\mu)(x)\biggr| \le
\|h\|_\infty\|\check{u_K}\mu\|\ \ (h\in C^K_c(G)).
\]
 Therefore $h\mapsto \int_G d h\mu(x)$ defines a
continuous linear functional of $C_c(G)$ with
respect to the inductive limit topology, which is
a Radon measure, say $\iota\mu\in
\textfrak{M}(G)$, namely \begin{equation} \int_G
d(h\mu)(x)=\int_G h(x)d\iota \mu(x) \ \ (h\in
C_c(G)). \end{equation} From (30) and (31), it
follows that \begin{eqnarray} \int_G k(x)(h(x)d
\iota \mu)(x)&=&\int_G (kh)(x)d\iota \mu(x)
=\int_G d (kh)\mu(x) \nonumber \\ &=&\int_G
k(x)dh\mu(x)\ \ (h, k\in C_c(G)). \end{eqnarray}
(32) implies $h\iota \mu=h\mu\ \ (h\in C_c(G))$.
Hence, for any $h\in C_c(G)$, we have
\begin{eqnarray*} \int_G h*h^*(x)d\iota\mu(x)
&=&\int_G d(h*h^*\mu)(x)=\lim_{n\rightarrow\infty}
\int_G h*h^*(x)d\mu_n(x) \\ &=& \lim_{n\rightarrow
\infty} \int_\Gamma |\check{h}(\gamma)|^2
\widehat{\mu_n}(\gamma)d\gamma =\int_\Gamma
|\check{h}(\gamma)|^2 \hat{\mu}(\gamma)d\gamma.
\end{eqnarray*} Therefore we have $\iota \mu\in
\textfrak{M}_T(G)$ with $d\widehat{\iota
\mu}(\gamma)=\hat{\mu}(\gamma)d\gamma.$ The
relation $h\mu=h\iota \mu$ for $h\in C_c(G)$ is
already shown above.\ \ \ \ \ $\Box$

\it

\vspace{0,5cm}

\it 

{\bf Theorem 7.7.}\ For $\mu\in {\mathcal M}(G)$, the
following {\rm (a)}, {\rm (b)}, and {\rm (c)} are
equivalent:

\ \ \ {\rm (a)}\ There exists $\iota\mu\in
\textfrak{M}_T(G)$ such that
$d\widehat{\iota\mu}(\gamma)=\hat{\mu}(\gamma)d\gamma.$\

\ \ \ {\rm (b)}\ $\hat{k}*\hat{\mu}\in B(\Gamma) \
\ (k\in C_c(G)\cap B(G))$.

\ \ \ 
{\rm (c)}\ $\mu\in {\mathcal M}(G)_{(0)}$.

\rm 

Proof. \ (a)$\Rightarrow$ (b) follows from Lemma
7.4, (b)$\Rightarrow $ (c) from Lemma 7.5, and
(c)$\Rightarrow$(a) from Lemma 7.6, respectively.
$\Box$

\vspace{0.5cm}

Let $\mu\in {\mathcal 
M}(G)_{(0)}$. The measure $\iota\mu\in
\textfrak{M}_T(G)$ in Theorem 7.7 (a) (which is
uniquely determined by [1, Theorem 2.1]) will be
called the associated Radon measure for $\mu$. 

\vspace{0.5cm}

\it

{\bf Theorem 7.8.}\ Let $\mu\in {\mathcal M}(G)$. The
following {\rm (a)} and {\rm (b)} are equivalent:

{\rm (a)}\
$\hat{\mu}(\gamma)d\gamma\in
\textfrak{M}_T(\Gamma)$.

{\rm (b)}\ $\mu\in {\mathcal M}(G)_{(0)}$ and its
associated Radon measure $\iota
\mu\in\textfrak{M}_T(G)$ is translation bounded.

\rm 

Proof.\ Suppose (a). If $\Gamma$ is non-discrete,
we have $\hat{\mu}(\gamma)d\gamma\in
\mathscr{I}(\Gamma)$ from [9, Theorem 2]. If
$\Gamma$ is discrete, then $G$ is compact, and the
Fourier transform of $\hat{\mu}(\gamma)d\gamma$ is
an element of $M(G)$ which is transformable. So,
$\hat{\mu}(\gamma)d\gamma \in
\mathscr{I}(\Gamma)$. From [1, Theorems 2.5 and
3.4], there exists a
translation bounded Radon measure $\iota \mu\in
\textfrak{M}_T(G)$ such that $d\widehat{\iota
\mu}(\gamma)= \hat{\mu}(\gamma)d\gamma$. Then
$\mu\in {\mathcal M}(G)_{(0)}$ from
Theorem 7.7 with the
associated Radon measure $\iota \mu$ which is
translation bounded; (b) follows

Suppose (b). Then (a) holds from [9, Theorem 2]
and Theorem 7.7. \ \ \ $\Box$

\vspace{0.5cm}

Now, we give two
examples which show that the associated Radon
measures $\iota \mu$ for $\mu\in {\mathcal 
M}(G)_{(0)}$ are not necessarily translation
bounded.

\vspace{0.5cm}

{\bf Example 7.1.}\ ($\mu\in {M}(G)_{(0)}$, with
$\iota \mu\in \textfrak{M}_T(G)$ which is not
translation bounded: the case when $G$ has an open
compact subgroup)

Suppose that $G$ is non-compact and has an
infinite open compact subgroup $H$. Taking a set
of representatives $\{x_n\}_{n=1}^\infty$,  of the
cosets $\{H+x: x\in G\}$, we have
$G=\cup_{n=1}^\infty H+x_n$. Let $V=H^{\perp}$,
the annihilator of $H$, and the Haar measures
$m_H$ and $m_{H^{\perp}}$ of $H$ and $H^{\perp}$,
respectively, are adjusted so that $m_H(H)=1$ and
$m_{H^{\perp}}(H^{\perp})=1$.

Since $\|\hat{f}\|_\infty\le \|f\|_1$ and two
norms $\|\hat{f}\|_\infty$ and $\|f\|_1$ in
$L^1(H)$ are not equivalent, for each $n\in
\mathbf{N}$, we can choose $f_n\in L^1(H)$ such
that $\|f_n\|_1=n, \|\widehat{f_n}\|_\infty\le
1/n^2$. Consider $f_n$ as an element of $M(G)$
which is concentrated on $H$, we put
$\mu_n=\sum_{k=1}^n f_k*\delta_{x_k},
n=1,2,3,....$, where $\delta_{x_k}$ is the unit
point mass concentrated on $x_k$ . Define
$\displaystyle\iota\mu:=\lim_{n\rightarrow \infty}
\mu_n\in \textfrak{M}(G)$ (converges in the vague
topology). Then we have 
\[ \|\mu_n\|_{\bar{V}}\le
\sum_{k=1}^n \|f_k*\delta_{x_k}\|_{\bar{V}} \le
\sum_{k=1}^n\|\widehat{f_k}\|^{\bar{V}}
\|\widehat{\delta_{x_k}}\|^{\bar{V}}
\le\sum_{k=1}^n 1/k^2\ (1\le n). \] Therefore
$\{\mu_n\}_{n=1}^\infty$ forms a Cauchy sequence
in $({\mathcal M}(G), \|\cdot \|_{\bar{V}})$, which
converges to an element, say $\mu\in {\mathcal 
M}(G);\, \mu=\lim_{n\rightarrow\infty} \mu_n$.
 
Let $h\in C_c(G)$ with $K={\rm supp}(h)\in {\mathcal 
K}(G)$, and choose $n_0\in \mathbf{N}$ such that
$K\subset \cup_{j=1}^{n_0} H+x_j$. Then
$h\mu_n=h\mu_{n_0} \ (n_0\le n)$. Hence
$\displaystyle \mu=\lim_{n\rightarrow\infty}\mu_n$
is a locally bounded representation, and $\mu\in
{\mathcal M}(G)_{(0)}$. For each $k\in C_{c,2}(G)$, we
have \[ \int_G k(x)d\iota\mu(x)=\lim_n\int_G
k(x)d\mu_n(x) =\lim_n\int_\Gamma \check{k}(\gamma)
\widehat{\mu_n}(\gamma)d\gamma =\int_\Gamma
\check{k}(\gamma)\hat{\mu}(\gamma)d\gamma,
\] and hence $d\widehat{\iota\mu}(\gamma)
=\hat{\mu}(\gamma)d\gamma.$ Since
$|\iota\mu|(H+x_n)=
\|f_n*\delta_{x_n}\|_1=\|f_n\|_1=n\rightarrow\infty\
(n\rightarrow \infty)$, the associated Radon
measure $\iota\mu$ for $\mu$ is not translation
bounded.

\vspace{0.5cm}

{\bf Example 7.2.}\ ($\mu\in {M}(G)_{(0)}$, with
$\iota \mu\in \textfrak{M}_T(G)$ which is not
translation bounded: the case when $G=\mathbf{R}$)

Let $G=\mathbf{R}$ with $\Gamma=\hat{\mathbf{R}}$,
and $V=(-1/4, 1/4)\subset \hat{\mathbf{R}}.$ For
each $n\in\mathbf{N}$, let $e_n\in
L^1(\mathbf{R})$ be such that $\|e_n\|_1\le 2,
\widehat{e_n}$ has compact support with
$\widehat{e_n}(t)=1\ \ (|t|\le n+1)$.
 
Let $\nu$ be a singular measure (with respect to
$dx$)\ in $M_0(\mathbf{R})$ such that ${\rm
supp}(\nu) \subset [0,\,1]$ and $\|\nu\|=1$. For
each $n\in \mathbf{N}$, there exists $N_n\in
\mathbf{N}$ from Theorem 3.4 (i) such that \[
\|\hat{\nu}\|^{\bar{V}, t}< 1/n^3\ \ (t\in
\hat{\mathbf{R}}, N_n< |t|), \] and put
$\nu_n:=\check{u_{[0,1]}}(\nu-e_{N_n}*\nu)$, where
$u_{[0,1]}$ is an element of
$L^1(\hat{\mathbf{R}})$ such that
$\|u_{[0,1]}\|_1\le 2, \check{u_{[0,1]}}\in
C_c(\mathbf{R})\cap B(\mathbf{R})$ with
$\check{u_{[0,1]}}(x)=1\ (x\in [0,1])$. Since
$\nu, \nu*e_{N_n}\in M(\mathbf{R})$, it follows
from Lemmas 3.3 and 7.3 that 
\begin{eqnarray} \|\widehat{\nu_n}\|^{\bar{V}}
&=&\|u_{[0,1]}*(\hat{\nu}-
\widehat{e_{N_n}}\hat{\nu})\|^{\bar{V}}\le
\|u_{[0,1]}\|_1\|\hat{\nu}-\widehat{e_{N_n}}\hat{\nu}\|^{\bar{V}}
\le 2\sup_{t\in
\hat{\mathbf{R}}}\|\hat{\nu}(1-\widehat{e_{N_n}})\|^{\bar{V},t}
\nonumber \\ &\le&2\sup_{t\in
\hat{\mathbf{R}}}\|\hat{\nu}\|^{\bar{V},t}\|1-\widehat{e_{N_n}}\|^{\bar{V},t}
\le2\sup_{|t|>N_n}\|\hat{\nu}\|^{\bar{V},t}\|1-\widehat{e_{N_n}}\|^{\bar{V},t}
\nonumber \\ &\le&
2\sup_{|t|>N_n}\|\hat{\nu}\|^{\bar{V},t}\|1-\widehat{e_{N_n}}\|^{\bar{V}}\le
2\sup_{|t|> N_n}\|\hat{\nu}\|^{\bar{V},t}
+2\sup_{|t|>
N_n}\|\hat{\nu}\|^{\bar{V},t}\|\widehat{e_{N_n}}\|^{\bar{V}}
\nonumber \\ &\le& 6\sup_{|t|> N_n}
\|\hat{\nu}\|^{\bar{V}, t}\le 6/n^3\ \
(n\in\mathbf{N}). 
\end{eqnarray}

Let $n_0:=\min\bigl\{n\in \mathbf{N}: {\rm
supp}({\check{u_{[0,1]}}})\subset [-n, n]\bigr\}$,
and define \[ \mu_n:=\sum_{j=1}^n
j\nu_j*\delta_{3(j-1)n_0} \in M(\mathbf{R}),\
n=1,2,3,.... \] From (33),
$\{\mu_n\}_{n=1}^\infty$ is a Cauchy sequence in
$(M(\mathbf{R}), \|\cdot \|_{\bar{V}})$. We put
$\mu:=\lim_{n\rightarrow\infty}\mu_n\in {\mathcal 
M}(\mathbf{R})$, and then
$\hat{\mu}(t)=\lim_{n\rightarrow\infty}\widehat{\mu_n}(t)\
\ (t\in\hat{\mathbf{R}})$. Since ${\rm
supp}(j\nu_j* \delta_{3(j-1)n_0})\cap {\rm
supp}(j'\nu_{j'}*\delta _{3(j'-1)n_0})=\emptyset\
(j\neq j')$, there exists $\iota\mu\in
\textfrak{M}(\mathbf{R})$ such that
$\mu_n\rightarrow \iota\mu$\ (converges in the
vague topology). Since
$\{\widehat{\mu_n}\}_{n=1}^\infty$ converges to
$\hat{\mu}$ uniformly on $\hat{\mathbf{R}}$, it
follows that

\begin{eqnarray*} \int_{\mathbf{R}}
h(x)d\iota\mu(x)&=&
\lim_{n\rightarrow\infty}\int_{\mathbf{R}}
h(x)d\mu_n(x)
=\lim_{n\rightarrow\infty}\int_{\hat{\mathbf{R}}}\check{h}(t)
\widehat{\mu_n}(t)dt \\ &=&
\int_{\hat{\mathbf{R}}}\check{h}(t)\hat{\mu}(t)dt\
\ \ (h\in C_{c,2}(\mathbf{R}))\, (\because
\check{h}\in L^1(\hat{\mathbf{R}})).
\end{eqnarray*}

Therefore $\iota\mu\in \textfrak{M}_T(\mathbf{R})$
with $d\widehat{\iota\mu}(t)=\hat{\mu}(t)dt$.
Since $e_{N_n}*\nu\in L^1(\mathbf{R})$ and $\nu$
is a singular measure, it follows that
\begin{eqnarray*}
&&|\iota\mu|\biggl(\biggl[3(n-1)n_0, 1+
3(n-1)n_0\biggr]\biggr) \\ &=&
|n\nu_n*\delta_{3(n-1)n_0}|([3(n-1)n_0,
1+3(n-1)n_0]) \\ &=&
|n(\check{u_{[0,1]}}\nu-\check{u_{[0,1]}}e_{N_n}*\nu)*\delta_{3(n-1)n_0}|([3(n-1)n_0,
1+3(n-1)n_0]) \\
&\ge&|n\nu*\delta_{3(n-1)n_0}|([3(n-1)n_0,
1+3(n-1)n_0]) =n\|\nu\|=n,\ \ n=1,2,3,...,
\end{eqnarray*} which implies that $\iota\mu$ is
not translation bounded.

\vspace{0.5cm}

Next, we give   the definition of '' locally integrability'' for elements of ${\mathcal L}^1(G)$, which is similar to ''locally boundedness'' for elements of ${\mathcal M}(G)$ above. 

\vspace{0.5cm}

{\bf Definition 7.4.}\ 
An element $f\in {\mathcal L}^1(G)$ is 
said to be locally integrable if $f$ has at least 
one ''locally integrable'' representation:
\[
f=\lim_{n\rightarrow\infty}f_n\ (f_n\in L^1(G), n=1,2,3,...)
\, {\rm in}\, ( {\mathcal L}^1(G), \|\cdot\|_{\bar{V}})
\]
such that $\{hf_n\}_{n=1}^\infty$ is a Cauchy sequence 
in $(L^1(G), \|\cdot\|_1)$ for every $h\in C_c(G)$. 
We write 
${\mathcal L}^1(G)_{(0)}:=\bigr\{
f\in {\mathcal L}^1(G): \, f\, {\rm has\ at\ least\ one\ locally
y\ integrable\  representation}\,\bigr\}. $

$L^1_{loc}(G)$ denotes the set of all locally integrable 
functions on $G$, and $L^1_{loc, T}(G)$ is 
the set of all $f\in L^1_{loc}(G)$ such that 
$f(x)dx$ is a transformable Radon measure in the 
sense of [1].   
\vspace{0,5cm}

In the following, we give lemmas and theorems for 
 locally integrability of elements in ${\mathcal L}^1(G)$ without proofs. Their proofs   
are similar to that of the corresponding  lemmas and theorems 
given above for locally boundedness of elements in ${\mathcal M}(G)$.

\vspace{0.5cm}

\it 

{\bf Lemma 7.9.}\, If $f\in {\mathcal L}^1(G)_{(0)}$ has two locally integrable 
representations\ $\displaystyle f=
\lim_{n\rightarrow\infty}f_{j,n}$
,$ j=1,2$, we
have $\|hf_{1,n}-hf_{2,n}\|_1\rightarrow 0\,
(n\rightarrow\infty) \ \ \ (h\in C_c(G)). $

\vspace{0.5cm}

\rm

{\bf Definition 7.5.}\ Suppose that 
$f\in {\mathcal L}^1(G)_{(0)}$ with a locally integrable 
representation 
\[
(*)\ \ f=\lim_{n\rightarrow\infty} f_n.
\] 
For each $h\in C_c(G)$, 
we define $hf\in L^1(G)$ 
by $\displaystyle hf=\lim_{n\rightarrow\infty}hf_n$, which is independent of the choice 
of the representation $(*)$ 
 with the help of Lemma 7.9. 

\vspace{0.5cm}

\it 

{\bf Lemma 7.10.}\  Let $f\in {\mathcal L}^1(G)$. 
If there exists 
$\iota f\in
L^1_{loc, T}(G)$ such that $d\widehat{\iota f}
(\gamma)=\hat{f}(\gamma)d\gamma$,  we have 
\[
\hat{k}*\hat{f}\in A(\Gamma)
 \ \ (k\in C_c(G)\cap B(G)).
\]
\rm

\vspace{0.5cm}

{\bf Lemma 7.11.}\ Suppose $f\in {\mathcal L}^1(G)$.
If $\hat{k}*\hat{f}\in A(\Gamma)\  
\ (k\in C_c(G)\cap B(G))$, we have 
$f\in {\mathcal L}^1(G)_{(0)}$. 

\rm

\vspace{0.5cm}

{\bf Lemma 7.12.}\ For 
$f\in {\mathcal L}^1(G)_{(0)}$
, there exists $\iota f\in
L^1_{loc, T}(G)$ such that 
$d\widehat{\iota f}(\gamma)=
\hat{f}(\gamma)d\gamma, $ 
and $hf=h\iota f$  holds for all $h\in C_c(G)$.

\vspace{0.5cm}

\it 

{\bf Theorem 7.13.}\ For $f\in {\mathcal L}^1(G)$, the
following {\rm (a)}, {\rm (b)}, and {\rm (c)} are equivalent:

\ \ \ {\rm (a)}\ There exists $\iota f\in
L^1_{loc, T}(G)$ such that
$d\widehat{\iota f}(\gamma)=\hat{f}(\gamma)d\gamma.$\

\ \ \ {\rm (b)}\ $\hat{f}*\hat{k}\in 
A(\Gamma)\ \ (k\in C_c(G)\cap B(G))$.  

\ \ \ {\rm (c)}\ $f\in {\mathcal L}^1(G)_{(0)}$.

\vspace{0.5cm}

For $f\in {\mathcal L}^1(G)_{(0)}$, the function 
$\iota f\in L^1_{loc,T}(G)$ in Theorem 7.13 (a) 
(which is uniquely determined by [1, Theorem 2.1])
 will be called the associated
locally integrable function for $f$. 

\vspace{0.5cm}

\it

{\bf Theorem 7.14.}\ For $f\in {\mathcal L}^1(G)$, 
the following {\rm (a)} and {\rm (b)}
 are equivalent:

{\rm (a)}\ $\hat{f}\in L^1_{loc, T}(\Gamma)$.

{\rm (b)}\ $f\in {\mathcal L}^1(G)_{(0)}$ and its 
associated locally integrable function 
$\iota f\in L^1_{loc, T}(G)$ is translation bounded.

\rm 

\vspace{0.5cm}

\it

{\bf Corollary 7.15}.}\  ${\mathcal L}^1_p(G)\subset {\mathcal L}^1(G)_{(0)}\ \ \ (1\le p\le 2)$. 

\rm 

Proof. \ Suppose $1\le p\le 2$ and $f\in {\mathcal L}^1_p(G)$. Since 
$\hat{f} \in L^p(\Gamma)$, we have $\hat{f}\in L^1_{loc, T}(\Gamma)$ 
by [1, Theorem 2.2], 
and $f\in {\mathcal L}^1(G)_{(0)}$ by  Theorem 7.14.
\ \ \ 
$\Box$

\vspace{1cm}
{\bf \S 8.\, Studies of the multiplier algebra
$\mathbb{M}({\mathcal L}^1(G))$; a natural continuous
isomorphism of $\mathbb{M}({\mathcal L}^1(G))$ into
$A(G)^*$, and locally bounded elements of
$\mathbb{M}({\mathcal L}^1(G))$}
 
\vspace{0.5cm}

\rm

{\bf Definition 8.1.}\ Since the generalized
Fourier transform ${\mathcal F}$ is an isometric
Banach algebra homomorphism from ${\mathcal L}^1(G)$
onto ${\mathcal A}(\Gamma)$, we can define the
multiplier algebra $\mathbb{M}({\mathcal L}^1(G))$ as
follows: for each $\phi\in \mathbb{M}({\mathcal 
A}(\Gamma))$, we define a multiplier $F_\phi$ of
${\mathcal L}^1(G)$ by \[ F_\phi*f\in {\mathcal L}^1(G), \
{\rm such\ that\ \ } {\mathcal 
F}(F_\phi*f)=\phi\hat{f}\ \ \ (f\in {\mathcal 
L}^1(G)). \] The set of all such multiplies,
$\bigl\{F_\phi: \phi\in \mathbb{M}({\mathcal 
A}(\Gamma))\bigr\}$, forms the multiplier algebra
$\mathbb{M}({\mathcal L}^1(G))$ of ${\mathcal L}^1(G)$
satisfying \[ \alpha F_{\phi_1}+\beta
F_{\phi_2}=F_{ \alpha \phi_1 +\beta \phi_2} ,\ \
F_{\phi_1}*F_{\phi_2}=F_{\phi_1\phi_2}\ \ \
(\alpha, \beta \in \mathbf{C}, \phi_1, \phi_2\in
\mathbb{M}({\mathcal A}(\Gamma))). \] We denote by
$\bar{\mathcal F}$ the inverse map of $\phi\mapsto
F_\phi$: \[ \bar{\mathcal F}: F_{\phi} \mapsto \phi,\
\mathbb{M}({\mathcal L}^1(G)) \rightarrow
\mathbb{M}({\mathcal A}(\Gamma)). \] For $F_\phi$, the
operator norm is 
\begin{eqnarray*}
\|F_\phi\|_{\mathbb{M}({\mathcal L}^1(G))}&=&
\sup\bigl\{\|F_\phi*f\|_{\bar{V}}: f\in {\mathcal 
L}^1(G), \|f\|_{\bar{V}}\le 1\bigr\} \\ &=&
\sup\bigl\{\|\phi\psi\|^{\bar{V}}: \psi\in {\mathcal 
A}(\Gamma), \|\psi\|^{\bar{V}}\le 1\bigr\}
=\|\phi\|_{\mathbb{M}({\mathcal 
A}(\Gamma))}=\|\phi\|^{\bar{V}}. 
\end{eqnarray*}
Therefore  $\bar{\mathcal F}$ is
an isometric algebra homomorphism of
$\mathbb{M}({\mathcal L}^1(G))$ onto $\mathbb{M}({\mathcal 
A}(\Gamma))$. We also use the simple notation: \[
\bar{\mathcal F}(F)=\hat{F} \ \ \ (F\in
\mathbb{M}({\mathcal L}^1(G))). \]
 
 \vspace{0.5cm}
 
{\bf Definition 8.2.}\ For each $\Omega\in {\mathcal 
K}(\Gamma)$, let $u_\Omega$ be an element of
$L^1(G)$ such that $\|u_\Omega\|_1\le 2,\
\widehat{u_\Omega}$ has compact support, and
$\widehat{u_\Omega}(\gamma)=1\ \ (\gamma\in
\Omega)$.

\vspace{0,5cm}

{\bf Lemma 8.1.}\ Let $\phi\in \mathbb{M}({\mathcal 
A}(\Gamma))$. Then

${\rm (i)}\ \{\phi\widehat{u_\Omega}\}_{\Omega\in
{\mathcal K}(\Gamma)}$ is a bounded net in ${\mathcal 
A}(\Gamma)$, contained in
$A(\Gamma)$, satisfying \[\lim_{\Omega\in {\mathcal 
K}(\Gamma)} (\phi\widehat{u_\Omega})\psi
=\phi\psi,\ \ (\psi\in {\mathcal A}(\Gamma)), 
\] hence
$\phi=\lim_{\Omega\in {\mathcal K}(\Gamma)}
\phi\widehat{u_\Omega}$ is a strongly convergent
bounded net representation in ${\mathcal A}(\Gamma)$.

{\rm (ii)}\ $F_\phi=\lim_{\Omega\in {\mathcal 
K}(\Gamma)} F_\phi*u_\Omega$\, is a strongly
convergent bounded net representation in ${\mathcal 
L}^1(G)$, where $F_\phi*u_\Omega$ is a function in
$L^1(G)$ such that $\widehat{F_\phi*u_\Omega}
=\phi\widehat{u_\Omega}\ \ (\Omega\in {\mathcal 
K}(\Gamma))$.

\rm

Proof.\ (i)\ That
$\{\phi\widehat{u_\Omega}\}_{\Omega\in 
{\mathcal K}(\Gamma)}$ is a bounded net in ${\mathcal 
A}(\Gamma)$ is seen from 
\[
\|\phi\widehat{u_\Omega}\|^{\bar{V}} \le
\|\phi\|^{\bar{V}}
\|\widehat{u_\Omega}\|^{\bar{V}}\le
2\|\phi\|^{\bar{V}}\ \ \ (\Omega\in {\mathcal 
K}(\Gamma)). 
\] Since $\phi\in
\mathbb{M}({\mathcal 
A}(\Gamma))=C^{\bar{V}}_{BSE}(\Gamma)\subset
A(\Gamma)_{loc}$ and $\widehat{u_\Omega}\in 
A(\Gamma)_c$, \,it follows that
$\phi\widehat{u_\Omega}\in A(\Gamma)$ (cf. [8,
Definition 5.1]).

For $\psi\in {\mathcal A}(\Gamma)$ and $\varepsilon>0$
we can choose $\Omega_{\varepsilon}\in {\mathcal 
K}(\Gamma)$ such that $\sup_{\gamma\in
\Gamma\setminus\Omega_\varepsilon}\|\phi\psi
\|^{\bar{V}, \gamma}\le\varepsilon/3$. 
Since
$\gamma+\bar{V}\subset \Omega_\varepsilon +
\bar{V}\ (\gamma\in \Omega_\varepsilon)$, it
follows that 
\begin{eqnarray*}
\|\phi\psi-\phi\widehat{u_\Omega}\psi\|^{\bar{V}}
&=& \sup_{\gamma\in
\Gamma}\|\phi\psi-\phi\widehat{u_\Omega}\psi\|^{\bar{V},\gamma}
=\sup_{\gamma\in \Gamma\setminus
\Omega_\varepsilon}\|\phi\psi-\phi\widehat{u_\Omega}\psi\|^{\bar{V},\gamma}
\\ &\le& \sup_{\gamma\in \Gamma\setminus
\Omega_\varepsilon}\|\phi\psi\|^{\bar{V},\gamma}
+\sup_{\gamma\in \Gamma\setminus
\Omega\varepsilon} \|\phi\psi\|^{\bar{V},\gamma}
\|\widehat{u_\Omega}\|^{\bar{V},\gamma} \\ &\le&
\varepsilon/3+2\varepsilon/3=\varepsilon \ \
(\Omega_\varepsilon +V \subset \Omega).
\end{eqnarray*}

(ii). An easy consequence of (i). \ \ \ $\Box$

\vspace{0.5cm}

Now we consider the dual space $A(G)^*$ of the
Fourier algebra $A(G)$ of $G$. Elements of
$A(G)^*$ are called pseudomeasures (cf. [11,\,\S
4.2]). Since the Fourier transform is an isometric
homomorphism from $L^1(\Gamma)$ onto $A(G)$, for
each pseudomeasure $\tau\in A(G)^*$, there exists
a unique $\xi \in L^{\infty}(\Gamma)$ such that 
 \[
\tau(\hat{f})=\int_\Gamma
f(-\gamma)\xi(\gamma)d\gamma \ \ (f\in
L^1(\Gamma)). \] In this case, $\xi$ is called the
Fourier transform of $\tau$\, (cf.[11,\,p.98]),
and denoted by $\xi=\hat{\tau}$. Of course we have
$\|\tau\|_{A(G)^*}=\|\xi\|_\infty$. Furthermore,
we can define multiplication $*$ in $A(G)^*$ by

"for $\tau_1, \tau_2\in A(G)^*$, define
$\tau_1*\tau_2 \in A(G)^*$ so that
$\widehat{\tau_1*\tau_2}=\widehat{\tau_1}
\widehat{\tau_2}$\ holds".

\noindent Obviously, the Fourier transform is a
surjective isometric algebra homomorphism of
$A(G)^*$ to $L^\infty(\Gamma)$.

For $\mu\in M(G)$,\ define $\tau_\mu(\hat{f}):=
\int_G \hat{f}(-x)d\mu(x)\ \ (f\in L^1(\Gamma))$,
then \begin{eqnarray}
\tau_\mu(\hat{f})&=&\int_G\hat{f}(-x)d\mu(x)
=\int_G \biggl(\int_\Gamma (\gamma,
x)f(\gamma)d\gamma\biggr)d\mu(x) \nonumber
\nonumber \\ &=& \int_\Gamma
f(\gamma)\biggl(\int_G (\gamma,
x)d\mu(x)\biggr)d\gamma \nonumber \\ &=&
\int_\Gamma f(\gamma)\check{\mu}(\gamma)d\gamma=
\int_\Gamma f(-\gamma)d\hat{\mu}(\gamma)d\gamma.
\end{eqnarray} From (34) we see that $\tau_\mu\in
A(G)^*$ which has the Fourier transform
$\hat{\mu}$. Since
$\tau_{\mu*\nu}(\hat{f})=\int_\Gamma f(-\gamma)
\hat{\mu}(\gamma)\hat{\nu}(\gamma)d\gamma\ \ (f\in
L^1(\Gamma))$, we have $\tau_{\mu*\nu}
=\tau_\mu*\tau_\nu$. Since $\|\tau_\mu\|_{A(G)^*}
=\|\hat{\mu}\|_\infty\le \|\mu\|$, $\mu \mapsto
\tau_\mu$ is a continuous algebra isomorphism from
$M(G)$ into $A(G)^*$.
 
Next, consider $\mu\in {\mathcal M}(G)$ with
$\displaystyle \mu=\lim_n\mu_n, (\mu_n\in M(G),
n\ge1)$. Note that $\widehat{\mu_n}$ converges to
$\hat{\mu}$ uniformly on $\Gamma$. We define
$\displaystyle\tau_\mu(\hat{f}):
=\lim_{n\rightarrow\infty}\int_\Gamma
\hat{f}(-x)d\mu_n(x)\ \ (f\in L^1(\Gamma)). $\
Then \begin{eqnarray*}
\tau_\mu(\hat{f})&=&\lim_{n\rightarrow \infty}
\int_G \hat{f}(-x)d\mu_n(x) \\ &=&
\lim_{n\rightarrow \infty} \int_\Gamma f(-\gamma)
\widehat{\mu_n}(\gamma)d\gamma = \int_\Gamma
f(-\gamma)\widehat{\mu}(\gamma)d\gamma \ \ (f\in
L^1(\Gamma)). 
\end{eqnarray*} 
Thus, $\tau_\mu\in
A(G)^*$ which has the Fourier transform
$\hat{\mu}$. Likewise in the case of $M(G)$,
$\mu\mapsto \tau_\mu, {\mathcal M}(G)\rightarrow
A(G)^*$ is a continuous algebra isomorphism.

\vspace{0.5cm}

\it 

{\bf Lemma 8.2.}\ Let \, $F\in
\mathbb{M}({\mathcal L}^1(G))$. Then $F*u_{\Omega}\in
L^1(G)\ (\Omega\in {\mathcal K}(\Gamma))$ 
{\rm (see Lemma 8.1)} and 
\[
\lim_{\Omega\in {\mathcal K}(\Gamma)}\int_G
\hat{f}(-x)F*u_\Omega(x)dx=\int_\Gamma
f(-\gamma)\hat{F}(\gamma)d\gamma\ \ (f\in
L^1(\Gamma)). \] \rm

Proof. Let $f\in L^1(\Gamma)$ and $\varepsilon>0$
be given, and let $\Omega_{f, \varepsilon}\in
{\mathcal K}(\Gamma)$ be chosen such that
$\int_{\Gamma\setminus \Omega_{f, \varepsilon}}
\bigl|f(-\gamma)|d\gamma\le
\varepsilon/(3(1+\|\hat{F}\|_\infty)).$ For
$\Omega\in {\mathcal K}(\Gamma) \ {\rm with}\
\Omega_{f,\varepsilon} \subset \Omega$, it follows
that \begin{eqnarray} &&\biggl|\int_\Gamma
f(-\gamma) \hat{F}(\gamma)d\gamma -\int_\Gamma
f(-\gamma) \hat{F}(\gamma)
\widehat{u_\Omega}(\gamma)d\gamma\biggr| \nonumber
\\ &&\hspace{2cm} \le \int_{\Gamma\setminus
\Omega}\bigl| f(-\gamma)\hat{F}
(\gamma)\bigr|d\gamma
+\|\widehat{u_\Omega}\|_\infty
\int_{\Gamma\setminus \Omega}\bigl|
f(-\gamma)\hat{F}(\gamma) \bigr|d\gamma \nonumber
\\ && \hspace{2cm} \le 3\|\hat{F}\|_\infty
\int_{\Gamma\setminus \Omega_{f,\varepsilon}}\bigl
|f(-\gamma)\bigr|d\gamma<\varepsilon.
\end{eqnarray} On the other hand, we have
\begin{eqnarray} \int_G \hat{f}(-x)F*u_\Omega(x)dx
&=& \int_G
\biggl(\int_\Gamma(x,\gamma)f(\gamma)d\gamma\biggr)
F*u_\Omega(x)dx \nonumber \\ &=& \int_\Gamma
f(-\gamma)\biggl(\int_G (-x,\gamma) F*u_\Omega(x)
dx\biggr)d\gamma \nonumber \\ &=&\int_\Gamma
f(-\gamma)\hat{F}(\gamma)\widehat{u_\Omega}
(\gamma)d\gamma. \end{eqnarray} From (35) and
(36), it follows that \begin{eqnarray*}
\lim_{\Omega\in {\mathcal K}(\Gamma)}\int_G
\hat{f}(-x) F*u_\Omega(x)dx &=&\lim_{\Omega\in
{\mathcal K}(\Gamma)}\int_\Gamma
f(-\gamma)\hat{F}(\gamma)
\widehat{u_\Omega}(\gamma)d\gamma \nonumber
\nonumber \\ &=&\int_\Gamma f(-\gamma)
\hat{F}(\gamma)d\gamma. \ \ \ \Box
\end{eqnarray*}\ \ \

\vspace{0.5cm}

{\bf Definition 8.3.}\ For each $F\in
\mathbb{M}({\mathcal L}^1(G))$, we define \[
\tau_F(\hat{f}) :=\lim_{\Omega\in {\mathcal 
K}(\Gamma)} \int_G \hat{f}(-x)F*u_\Omega(x)dx \ \
(f\in L^1(\Gamma)). \] 
\vspace{0.5cm}

\it 

{\bf Theorem 8.3.}\ {\rm (i)}\ For each $F\in
\mathbb{M}({\mathcal L}^1(G))$, we have $\tau_F\in
A(G)^*$ with the Fourier transform $\hat{F}$, and
the map $F\mapsto \tau_F, \, 
\mathbb{M}({\mathcal 
L}^1(G))\rightarrow A(G)^*$ is a continuous
algebra isomorphism.

{\rm (ii)}\ When $G$ is compact, the above map
$F\mapsto \tau_F$ is a surjective isometric
algebra homomorphism.

\rm

Proof.\ (i)\ By Lemma 8.2 and Definition 8.3,
$\tau_F$ is an element of $A(G)^*$ whose Fourier
transform is $\hat{F}$. That the map $F\mapsto
\tau_F,\, \mathbb{M}({\mathcal L}^1(G)) \rightarrow
A(G)^*$ is a continuous algebra isomorphism is
apparent.

(ii)\ If $G$ is compact, we can take $V=\{0\}$.
Then $({\mathcal A}(\Gamma), \|\cdot\|^{\bar{V}})
=(C_0(\Gamma), \|\cdot\|_\infty)$ and ${\mathbb
M}({\mathcal A}(\Gamma))=C_b(\Gamma)$. Therefore the
map $F\mapsto \tau_F,\, \mathbb{M}({\mathcal 
L}^1(G))\rightarrow A(G)^*$ is surjective and
isometric. \ \ $\Box$

\vspace{0.5cm}

Remark 8.1.\ If $G$ is compact and $1\le p\le2$,
we have ${\mathcal L}^1_p(G))$ from Remark 6.1 (ii) and
$\mathbb{M}({\mathcal A}_p(\Gamma))=\mathbb{M}({\mathcal 
A}(\Gamma))=C_{BSE}^{\bar{V}}(\Gamma)$ from
Theorem 6.1(ii). Therefore, as
multiplier algebras, they are isomorphic:  \[
\mathbb{M}(A_p(G))\cong \mathbb{M}({\mathcal 
L}^1_p(G))\cong \mathbb{M}({\mathcal L}^1(G)). 
\] The
fact that there is a natural continuous algebra
isomorphism from $\mathbb{M}(A_p(G))$ onto
$A(G)^*$ was already known (cf.\,[11, Corollary
6.4.1]).

\vspace{0.5cm}

In the remainder of this section, we investigate a
similar problem for elements of $\mathbb{M}({\mathcal 
L}^1(G))$ as that, ''locally boundedness'', for
elements of ${\mathcal M}(G)$ which is investigated in
$\S 7.$ We treat only the case where
$G=\mathbf{R}^d, \Gamma=\hat{\mathbf{R}}^d$. The
symbols $x\cdot t, |x|,\, |t|,\, $ for
$x=(x_1,...,x_d)\in \mathbf{R}^d,
t=(t_1,...,t_d)\in \hat{\mathbf{R}}^d$ mean
$x\cdot t=\sum_{j=1}^d x_jt_j,\
|x|=\bigl(\sum_{j=1}^d{x_j}^2\bigr)^{1/2},\
|t|=\bigr(\sum_{j=1}^d {t_j}^2\bigr)^{1/2},\ $ and
let $V=\{t\in \hat{\mathbf{R}}^d: |t|<1/4\}$. The
Fourier transforms of $f\in L^1(\mathbf{R}^d)$ and
$\xi\in L^1(\hat{\mathbf{R}}^d)$ are
$\hat{f}(t)=\frac{1}{(2\pi)^d}\int_{\mathbf{R}^d}e^{-ix\cdot
t}f(x)dx$ and $\hat{\xi}(x)
=\int_{\hat{\mathbf{R}}^d}e^{-i t\cdot
x}\xi(t)dt$, where $dx$ and $dt$ are the Lebesgue
measures of $\mathbf{R}^d$ and
$\hat{\mathbf{R}}^d$, respectively.

We set $B_n:=\{x\in \mathbf{R}^d: |x|\le 1/n\},\,
n=1,2,3,...$, and let
$e_0(x):=c1_{B_4}*1_{B_4}*1_{B_4}*1_{B_4}(x) \in
C_{c, 2}(\mathbf{R}^d) $, where $1_{B_4}$ is the
characteristic function of $B_4$ and $c>0$ is
chosen so that $\|e_0\|_1=1$. We define
$e_n(x):=n^d e_0(nx)$ $(x\in \mathbf{R}^d, n\in
\mathbf{N})$. It follows that ${\rm
supp}(e_n)\subset B_n , \ \|e_n\|_1=1$ and
$\displaystyle\lim_{n\rightarrow
\infty}\widehat{e_n}(t)=\lim_{n\rightarrow\infty}
\widehat{e_0}(\frac{1}{n}t)=\widehat{e_0}(0)=1\
(t\in \hat{\mathbf{R}}^d)$.

\vspace{0.5cm}

Let $F_\phi \in\mathbb{M}({\mathcal 
L}^1(\mathbf{R}^d))$ with $\phi\in
\mathbb{M}({\mathcal A}(\hat{\mathbf{R}}^d))$. Then \[
F_\phi*e_n\in {\mathcal L}^1(\mathbf{R}^d)\ {\rm
with}\ {\mathcal F}(F_\phi*e_n)=\phi\widehat{e_n}, \
n=1,2,3,... \] Since
$\displaystyle\lim_{n\rightarrow\infty}\widehat{F_\phi*e_n}
(t)=\phi(t)\ \ (t\in \hat{\mathbf{R}}^d)$ and
$\sup_{n\in \mathbf{N}}\|F_\phi*e_n\|_{\bar{V}}\le
\|\phi\|^{\bar{V}}<\infty$, we say that
$F_\phi=\lim_{n\rightarrow \infty}F_\phi*e_n $ is
a bounded net representation of $F_\phi$ in ${\mathcal 
L}^1(\mathbf{R}^d)$ (cf. [18, Theorem 4 (i)]).

\vspace{0.5cm}

{\bf Definition 8.4.}\ Suppose $F_\phi\in
\mathbb{M}({\mathcal L}^1(\mathbf{R}^d))$ with
$\phi\in \mathbb{M}({\mathcal 
A}(\hat{\mathbf{R}}^d))$. $F_\phi$ is said to be
locally bounded if the following conditions (i)
and (ii) are satisfied:

(i) $F_\phi*e_n\in {\mathcal L}^1(\mathbf{R}^d)_{(0)}\
\ (n\in \mathbf{N})$.

(ii)\ For each $h\in C_c(\mathbf{R}^d)$,\ we have\
$\sup_{n\in
\mathbf{N}}\|h(F_\phi*e_n)\|_1<\infty.$ \ \
\vspace{0.2cm}

The symbol $\mathbb{M}({\mathcal 
L}^1(\mathbf{R}^d))_{(0)}$ is used to denote the
set of all locally bounded elements of
$\mathbb{M}({\mathcal L}^1(\mathbf{R}^d))$.

\vspace{0.5cm}

\it

{\bf Lemma 8.4.}\ Let $F_\phi\in \mathbb{M}({\mathcal 
L}^1(\mathbf{R}^d))_{(0)}$ with $\phi\in
\mathbb{M}({\mathcal A}(\hat{\mathbf{R}}^d))$. For $
k\in C_c(\mathbf{R}^d)\cap B(\mathbf{R}^d)$, we
have 
\begin{eqnarray*} &&{\rm (i)}\
\lim_{n\rightarrow\infty}
\frac{1}{(2\pi)^d}\int_{\mathbf{R}^d} e^{-ix\cdot
t} k(F_\phi*e_n)(x)dx= \hat{k}*\phi(t)\ \ \ (t\in
\hat{\mathbf{R}}^d).\ \ \ \ \\ && {\rm (ii)}\
\hat{k}*\phi\in B(\hat{\mathbf{R}}^d).
\end{eqnarray*}
 
 \rm

Proof.\ (i)\ For each $n\in \mathbf{N}$,
$F_\phi*e_n\in {\mathcal L}^1(\mathbf{R}^d)_{(0)}$,
and $F_\phi*e_n$ has a locally bounded
representation $\displaystyle
F_\phi*e_n=\lim_{j\rightarrow\infty}f_j,\, f_j\in
L^1(\mathbf{R}^d)$. So, for each $h\in
C_c(\mathbf{R}^d), \{hf_j\}_{j=1}^\infty$ is a
Cauchy sequence in $(L^1(\mathbf{R}^d), \|\cdot
\|_1)$ with a limit expressed by $h(F_\phi*e_n)$.
It follows from Lemma 3.3 that \begin{eqnarray*}
\frac{1}{(2\pi)^d}\int_{\mathbf{R}^d}e^{-ix\cdot
t}k(F_\phi*e_n)(x)dx &=&\lim_{j\rightarrow
\infty}\frac{1}{(2\pi)^d}\int_{\mathbf{R}^d}
e^{-ix\cdot t}k(x)f_j(x)dx \\ &=&
\lim_{j\rightarrow\infty}\hat{k}*\widehat{f_j}(t)
=\hat{k}*(\phi\widehat{e_n})(t)\ \ (t\in
\hat{\mathbf{R}}^d). \end{eqnarray*} Taking
$n\rightarrow\infty$ in the above, we obtain (i).

(ii)\ It follows from (i) and $\sup_{n\in
\mathbf{N}}\|k(F_\phi*e_n)\|_1 <\infty$
(Definition 8.4) that $\hat{k}*\phi\in
B(\hat{\mathbf{R}}^d)$ ([17, Theorem 1.9.2], [18,
Theorem 4\,(i)]). \ \ $\Box$ \vspace{0.5cm}

\it
 
 \vspace{0.5cm}
 
{\bf Lemma 8.5.}\ Let $F_\phi\in \mathbb{M}({\mathcal 
L}^1(\mathbf{R}^d))$ with $\phi\in
\mathbb{M}({\mathcal A}(\hat{\mathbf{R}}^d))$. Suppose
that there exists $\iota F_\phi\in
\textfrak{M}_T(\mathbf{R}^d)$ such that
$d\widehat{\iota F_\phi}(t)= \phi(t)d t$. Then,
for $n=1,2,3,...$, we have

{\rm (i)}\ $\iota F_\phi*e_n\in
\textfrak{M}_T(\mathbf{R}^d)$ with
$d\widehat{\iota F_\phi*e_n}(t)
=\widehat{e_n}(t)\phi(t)dt$.

{\rm (ii)}\ $F_\phi*e_n \in {\mathcal 
L}^1(\mathbf{R}^d)_{(0)}$ whose associated locally
integrable function is given by $\iota
F_\phi*e_n$.

{\rm (iii)}\ $h(F_\phi*e_n)=h(\iota F_\phi*e_n) \
\ (h\in C_c(\mathbf{R}^d))$. \rm

Proof.\ (i)\ Let $x\in \mathbf{R}^d$. Since $e_n(x-y)\in C_{c,
2}(\mathbf{R}^d)$ (as a function of $y$) with the
inverse Fourier transform $e^{it\cdot
x}\widehat{e_n}(t)$, we have from (23) that
\begin{equation} \iota
F_\phi*e_n(x)=\int_{\mathbf{R}^d} e_n(x-y)d\iota
F_\phi(y)=\int_{\hat{\mathbf{R}}^d} e^{i t\cdot
x}\widehat{e_n}(t)\phi(t)d t. \end{equation} Since
$\widehat{e_n}(t)\phi(t)\in
L^1(\hat{\mathbf{R}}^d)(\subset
\mathscr{I}(\hat{\mathbf{R}}^d))$\, ([1, Theorem
3.5]), it follows from (37) and [1,Theorem 3.4]
that $\iota F_\phi*e_n\in
\textfrak{M}_T(\mathbf{R}^d)$ with
$d\widehat{\iota
F_\phi*e_n}(t)=\widehat{e_n}(t)\phi(t)dt$.

(ii)\ For each $m\in \mathbf{N}$, let $u_m$ be in
$L^1(\mathbf{R}^d)$ such that $\|u_m\|_1\le 2,
\widehat{u_m}\in C_c(\hat{\mathbf{R}}^d),\,$ and\,
$\widehat{u_m}(t)=1\ \ (|t|\le m)$. We claim that
$\displaystyle F_\phi*e_n
=\lim_{m\rightarrow\infty}u_m*(F_\phi*e_n)$ is a
locally integrable representation, and hence
$F_\phi*e_n\in {\mathcal L}^1(\mathbf{R}^d)_{(0)}$.\

To show the claim, we note that, for each $m\in
\mathbf{N}$, 
\[ {\mathcal 
F}(u_m*(F_\phi*e_n))=\widehat{u_m}(\phi\widehat{e_n})\in 
A(\hat{\mathbf{R}}^d)\,
(\,\because \phi\widehat{e_n}\in
A(\hat{\mathbf{R}}^d)_{loc}), \] hence we have
$u_m*(F_\phi*e_n)\in L^1(\mathbf{R}^d)$. Moreover,
we have 
\begin{eqnarray*} 
\|\phi\widehat{e_n}-
\widehat{u_m}\phi\widehat{e_n}\|^{\bar{V}} &\le&
\sup_{t\in \hat{\mathbf{R}}^d}\|\phi\|^{\bar{V},t}
\|\widehat{e_n}-\widehat{u_m}\widehat{e_n}
\|^{\bar{V}, t} \\ &\le&
\|\phi\|^{\bar{V}}\sup_{t\in \hat{\mathbf{R}}^d:
|t|\ge
m-1}\|\widehat{e_n}-\widehat{u_m}\widehat{e_n}\|^{\bar{V},
t} \\ &\le& \|\phi\|^{\bar{V}}\sup_{t\in
\hat{\mathbf{R}}^d: |t|\ge
m-1}\biggl(\|\widehat{e_n}\|^{\bar{V},t}+\|\widehat{u_m}\|^{\bar{V},t}\|\widehat{e_n}\|^{\bar{V},
t}\biggr) \rightarrow 0\ \ (m\rightarrow \infty).
\end{eqnarray*}
 Further, since
$\phi\widehat{e_n}\in L^1(\hat{\mathbf{R}}^d)$, we
have  
\begin{eqnarray*}
 \|\widehat{u_m}(\phi\widehat{e_n})-
\widehat{u_{m'}}(\phi\widehat{e_n})\|_1&=&
\|(\widehat{u_m}\phi\widehat{e_n}-\phi\widehat{e_n}\|_1
+\|\phi\widehat{e_n}-\widehat{u_{m'}}\phi\widehat{e_n}\|_1
\\
&\le&
3\int_{|t|>m}|\phi(t)\widehat{e_n}(t)|dt 
+3\int_{|t|>m'}\phi(t)\widehat{e_n}(t)|dt
\\
&& \rightarrow 0\ \ (m, m'\rightarrow \infty), 
\end{eqnarray*}
 and
hence (with the help of the inversion theorem) we
have \[ \|u_m*(F_\phi*e_n)-
{u_{m'}}*(F_\phi*e_n)\|_\infty \rightarrow 0\, (m,
m'\rightarrow \infty). \] Therefore $\displaystyle
F_\phi*e_n= \lim_{m\rightarrow\infty}
u_m*(F_\phi*e_n)$ is a representation such that
$\{h(u_m*(F_\phi*e_n))\}_{m=1}^\infty$ is a Cauchy
sequence in $L^1(\mathbf{R}^d)$ for every $h\in
C_c(\mathbf{R}^d)$, and the claim follows.

By the claim and Theorem 7.7, there exists $\iota
(F_\phi*e_n)\in L^1_{loc, T}(\mathbf{R}^d)$ whose
Fourier transform is $\phi(t)\widehat{e_n}(t)dt.$
This result and (i) imply that
$\iota(F_\phi*e_n)=\iota F_\phi*e_n$, and this
completes the proof of (ii).

(iii)\ Let $h\in C_c(\mathbf{R}^d)$. The relation
$h(F_\phi*e_n)=h(\iota F_\phi*e_n)$ follows from
(ii) and Lemma 7.6. \ \ \ \ \ $\Box$

\vspace{0.5cm}

\it

{\bf Theorem 8.6.}\ For $F_\phi \in
\mathbb{M}({\mathcal L}^1(\mathbf{R}^d))$ with
$\phi\in \mathbb{M}({\mathcal 
A}(\hat{\mathbf{R}}^d))$, the following {\rm (a)},
{\rm (b)}, and {\rm (c)} are equivalent:

\ \ {\rm (a)}\ There exists $\iota F_\phi\in
\textfrak{M}_T(\mathbf{R}^d)$ such that
$d\widehat{\iota F_\phi}(t)=\phi(t)d t$.

\ \ {\rm (b)}\ $\hat{k}*\phi\in
B(\hat{\mathbf{R}}^d) \ (k\in
C_c(\mathbf{R}^d)\cap B(\mathbf{R}^d))$.

\ \ {\rm (c)}\ $ F_\phi \in \mathbb{M}({\mathcal 
L}^1(\mathbf{R}^d))_{(0)}$.

\rm

Proof.\ $(a)\Rightarrow (c)$: Suppose (a).  From
Lemma 8.5 (ii) and (iii), we have $F_\phi*e_n\in
{\mathcal L}^1(\mathbf{R}^d)_{(0)}$\, 
with $h(\iota F_\phi*e_n)=h(F_\phi*e_n)\ (n\in \mathbf{N}, h\in C_c(\mathbf{R}^d))$. So, for any $h\in
C_c(\mathbf{R}^d)$ with $K={\rm supp}(h)$, it
follows that \begin{eqnarray*}
 \sup_{n\in
\mathbf{N}}\|h(F_\phi*e_n)\|_1 &=&\sup_{n\in
\mathbf{N}} \|h(\iota F_\phi*e_n)\|_1 \\ &=&
\sup_{n\in \mathbf{N}} \|h((1_{K-B_1}\iota
F_\phi)*e_n)+h((1_{\mathbf{R}^d\setminus
(K-B_1)}\iota F_\phi)*e_n)\|_1 \\ &=&\sup_{n\in
\mathbf{N}} \|h((1_{K-B_1}\iota F_\phi)*e_n)\|_1\,
\\ &\le& \sup_{n\in \mathbf{N}}\|h\|_\infty
\|e_n\|_1\|1_{K-B_1}\iota F_\phi\|
=\|h\|_\infty\|1_{K-B_1}\iota F_\phi\|<\infty.
\end{eqnarray*} By Definition 8.4 we have
$F_\phi\in \mathbb{M}({\mathcal 
L}^1(\mathbf{R}^d))_{(0)}$; (c) follows. \ \ \

(c)$\Rightarrow$(b) follows from Lemma 8.4 (ii).

(b)$\Rightarrow$(a):\, Suppose (b). For $k\in
C_c(\mathbf{R}^d)\cap B(\mathbf{R}^d)$,
$\hat{k}*\phi\in B(\hat{\mathbf{R}}^d)$, and there
exists $\nu\in M(\mathbf{R}^d)$ such that
$\hat{\nu}= \hat{k}*\phi$. We denote this $\nu$ by
$kF_\phi$. It follows that $h(kF_\phi)=(hk)F_\phi\
\ (h, k\in C_c(\mathbf{R}^d)\cap B(\mathbf{R}^d))$
since $\widehat{h(kF_\phi)}=\widehat{(hk)F_\phi}
=\hat{h}*\hat{k}*\phi$.

For $H, K\in {\mathcal K}(\mathbf{R}^d)$ with $H\cap
K\neq \emptyset$, we have \begin{equation}
k(\check{u_{ H}}F_\phi)=k(\check{u_{K}}F_\phi)=k
F_\phi \ \ \ (k\in C^{H\cap K}_c(\mathbf{R}^d)
\cap B(\mathbf{R}^d)). \end{equation} Since
$C^{H\cap K}_c(\mathbf{R}^d) \cap B(\mathbf{R}^d)$
is dense in $C_c^{H\cap K}(\mathbf{R}^d)$, it
follows from (38) that \[ k(\check{u_{
H}}F_\phi)=k(\check{u_{K}}F_\phi)\ \ \ (k\in
C^{H\cap K}_c(\mathbf{R}^d)). \] Therefore the map
$h\mapsto \lim_{K\in {\mathcal 
K}(\mathbf{R}^d)}\int_{\mathbf{R}^d}
h(x)d\check{u_K}F_\phi(x)\ \ (h\in
C_c(\mathbf{R}^d))$, namely \[ h\mapsto
\int_{\mathbf{R}^d} h(x)d (\check{u_K}F_\phi)(x)\
\ \ (h\in C^K_c(\mathbf{R}^d), \ K\in {\mathcal 
K}(\mathbf{R}^d)) \] defines a linear functional
of $C_c(\mathbf{R}^d)$ which is continuous with
respect to the inductive limit topology since, for
each $K\in {\mathcal K}(\mathbf{R}^d)$, we have \[
\sup\biggl\{\biggl|\int_{\mathbf{R}^d}
h(x)d(\check{u_K}F_\phi)(x)\biggr|:\,h\in
C_c^K(\mathbf{R}^d),\ \|h\|_\infty\le 1\biggr\}\le
\|\check{u_K}F_\phi\|<\infty. \] It is a Radon
measure, which we denote $\iota F_\phi\in
\textfrak{M} (\mathbf{R}^d);\ \iota F_\phi:=
\lim_{K\in {\mathcal K}(\mathbf{R}^d)}
\check{u_K}F_\phi$ (converges in the vague
topology). Then the following formula is valid:
\begin{eqnarray} \int_{\mathbf{R}^d} k(x)d\iota
F_\phi(x) &=& \lim_{K\in {\mathcal 
K}(\mathbf{R}^d)}\int_{\mathbf{R}^d}k(x)d\check{u_K}F_\phi(x)
\nonumber \\ &=& \lim_{K\in {\mathcal 
K}(\mathbf{R}^d)}\int_{\mathbf{R}^d}
d(k\check{u_K})F_\phi(x)=\int_{\mathbf{R}^d} d
kF_\phi(x) \nonumber \\ &=&
\widehat{kF_\phi}(0)=\hat{k}*\phi(0)=
\int_{\hat{\mathbf{R}}^d}\check{k}(t)\phi(t)dt\ \
\ (k\in C_{c,2}(\mathbf{R}^d)). \end{eqnarray}
From (39), we have $\iota F_\phi\in
\textfrak{M}_T(\mathbf{R}^d)$ with
$d\widehat{\iota F_\phi}(t)=\phi(t)dt$; (a)
follows.\ \ \ $\Box$

\vspace{0.3cm}

For $F_\phi\in \mathbb{M}({\mathcal 
L}^1(\mathbf{R}^d))_{(0)}$, the measure $\iota
F_\phi \in \textfrak{M}_T(\mathbf{R}^d)$ in
Theorem 8.6 (a) (which is uniquely determined from
[1, Theorem 2.1]) will be called the associated
Radon measure for $F_\phi$. \vspace{0.5cm}

\it 

{\bf Theorem 8.7.}\ For
$\phi\in \mathbb{M}({\mathcal A} (\hat{\mathbf{R}}^d))$, the
following {\rm (a)} and {\rm (b)} are equivalent:

{\rm (a)}\ $\phi(t)dt\in \textfrak{M}_T(\hat{\mathbf{R}}^d)$.

{\rm (b)}\ $F_\phi\in \mathbf{M}({\mathcal 
L}^1(\mathbf{R}^d))_{(0)}$ whose associated Radon
measure $\iota F_\phi$ is translation bounded.

\rm 

Proof.\ Suppose (a). We have $\phi(t)dt\in
\mathscr{I}(\hat{\mathbf{R}}^d)$ from [9, Theorem 2]. So, from
[1, Theorems 2.5 and 3.4], there exists a
translation bounded Radon measure $\iota F_\phi
\in \textfrak{M}_T(\mathbf{R}^d)$ such that
$d\widehat{\iota F_\phi}(t)=\phi(t)dt$. Hence (b)
follows from Theorem 8.6.

Suppose (b). Then (a) holds from [9, Theorem 2]
and Theorem 8.6.
\ $\Box$ 
  
\vspace{0.5cm}

\rm

Finally, we will give two examples of $F_\phi\in
\mathbb{M}({\mathcal L}^1(\mathbf{R}))$; one is
$F_\phi\not\in \mathbb{M}({\mathcal 
L}^1(\mathbf{R}))_{(0)}$, and the other is
$F_\phi\in \mathbb{M}({\mathcal 
L}^1(\mathbf{R}))_{(0)}$ with $\phi\in
\mathbb{M}({\mathcal A}(\hat{\mathbf{R}}))\setminus
{\mathcal B}(\hat{\mathbf{R}})$.

\vspace{0.5cm}

{\bf Example 8.1\ ($F_\phi\in 
\mathbb{M}({\mathcal 
L}^1( \hat{\mathbf{R}}))\setminus \mathbb{M}({\mathcal L}^1(\hat{\mathbf{R}}))_{(0)}) $

\rm \ Consider $\mathbf{R}$ with the dual group
$\hat{\mathbf{R}}$, and let $V=(-1/4, 1/4)$. Let
$F_\phi\in \mathbb{M}({\mathcal L}^1(\mathbf{R}))$
with $0\le\phi\in \mathbb{M}({\mathcal 
A}(\hat{\mathbf{R}}))$, where $\phi\in
C^2(\hat{\mathbf{R}})$ such that $\phi(t)=0\ (t\le
0), $ and $=1\ (1\le t)$. That 
$\phi\in \mathbb{M}({\mathcal 
A}(\hat{\mathbf{R}}))(=C^{\bar{V}}_{BSE}(\hat{\mathbf{R}}))$
is easy to see using the relation
$C^2(\hat{\mathbf{R}})\cap
C_c(\hat{\mathbf{R}})\subset A(\hat{\mathbf{R}})$:
let $\psi\in C_c(\hat{\mathbf{R}})\cap
C^2(\hat{\mathbf{R}})$ with $\psi(t)=1\ (|t|\le
2)$, then \[ \|\phi\|^{\bar{V}}=\sup_{ t\in
\hat{\mathbf{R}}} \|\phi\|^{\bar{V},t} \le
\sup_{t\le 3/2} \|\psi\phi\|^{\bar{V,t}}+ \sup_{
3/2\le t} \|\phi\|^{\bar{V},t}\le
\|\psi\phi\|_{A(\hat{\mathbf{R}})}+1<\infty. \]

We show that $F_\phi\not\in 
\mathbb{M}({\mathcal 
L}^1(\mathbf{R}))_{(0)}$. To see this, suppose
contrary that $F_\phi\in \mathbf{M}({\mathcal 
L}^1(\mathbf{R}))_{(0)}$, and let $0\neq h\in
C_c(\mathbf{R})$ be fixed. From Theorem 8.6, we
have $|\hat{h}|^2*\phi\in B(\hat{\mathbf{R}})$. It
follows that \begin{eqnarray} \lim_{t\rightarrow
+\infty}|\hat{h}|^2*\phi(t)&=&\lim_{t\rightarrow
+\infty} \int_0^\infty |\check{h}(s-t)|^2\phi(s)ds
=\int_{-\infty}^\infty |\check{h}(s)|^2ds>0,\ \\
\lim_{t\rightarrow -\infty}|\hat{h}|^2*\phi(t)&=&
\lim_{t\rightarrow -\infty}\int_0^\infty
|\check{h}(s-t)|^2\phi(s)ds=0. \end{eqnarray} By
[17, Theorem 2.7.2], there exists $\psi\in
B(\hat{\mathbf{Z}})$ such that
$\psi(n)=|\hat{h}|^2*\phi(n)\ \ (n\in
\hat{\mathbf{Z}})$, so from (40), (41), we have 
\[
\lim_{n\rightarrow -\infty} \psi(n)=0,\ {\rm and}\
\lim_{n\rightarrow +\infty}
\psi(n)=\int_{-\infty}^\infty|\check{h}(s)|^2d
s>0, \] which contradicts a Rajchman's theorem
(cf. [2]): \[''\zeta\in B(\hat{\mathbf{Z}})\ {\rm
with}\ \lim_{n\rightarrow -\infty} \zeta(n)=0\
{\rm implies}\ \lim_{n\rightarrow +\infty}
\zeta(n)=0''. \]

{\bf Example 8.2\ ($F_\phi\in \mathbb{M}({\mathcal 
L}^1(\mathbf{R}))_{(0)}, \,\phi\in
\mathbb{M}({\mathcal A}(\hat{\mathbf{R}}))\setminus
{\mathcal B}(\hat{\mathbf{R}})))$}\

Consider the real group $\mathbf{R}$ with the dual
group $\hat{\mathbf{R}}$, and let $V=(-1/4,
1/4)\subset \hat{\mathbf{R}}$. From Lemma 5.2,
there exists a function $\psi$ on
$\hat{\mathbf{Z}}$ with $\psi=\pm 1$ and
\begin{equation} \sup_{\zeta\in
B(\hat{\mathbf{Z}})}\sup_{n\in \hat{\mathbf{Z}}}
|\psi(n)-\zeta(n) |\ge1/2. 
\end{equation}

Let $\ell\in C_{c,2}(\hat{\mathbf{R}})$ supported
on $V$ and $\ell(0)=1$. Put $\ell_0=\psi(0)\ell$.
Since $\check{\ell}\in C_0(\mathbf{R})$, we can
choose $N_1\in \mathbf{N}$ 
such that if we put
\[ \ell_1(t):=\psi(-1)e^{i2\pi N_1t}
\ell(t+1)+\psi(1)e^{-i2\pi N_1t}\ell(t-1)\ \ (t\in
\hat{\mathbf{R}}), 
\]
 the inequality $\sup_{t\in
\hat{\mathbf{R}}}\bigl (
|\check{\ell_0}(t)|+|\check{\ell}_1(t)|\bigr)\le
\|\check{\ell}\|_\infty(1+ 1/2)$ holds. 
Similarly,
we can choose $(N_1<)N_2\in \mathbf{N}$
 such that 
 if we put \[ \ell_2(t):=\psi(-2)e^{i 2\pi N_2
t}\ell(t+2)+\psi(2)e^{-i2\pi N_2 t}\ell(t-2), \]
the inequality $\sup_{x\in \hat{R}}
(|\check{\ell}_0(x)|+|\check{\ell_1}(x)|+|\check{\ell_2}(x)|)\le\|\check{\ell}\|_\infty(1+1/2+1/2^2)$
holds. Repeating this process for $n=1,2,3,...,$
we have a strictly increasing sequence of positive
integers $\{N_n\}_{n=1}^\infty$ such that if we
put $\ell_n(t):=\psi(-n)e^{i2\pi N_n
t}\ell(t+n)+\psi(n)e^{-i2\pi N_n t}\ell(t-n), n\in
\mathbf{N}, $ then $\sup_{x\in
\mathbf{R}}(\sum_{j=0}^n |\check{\ell}_j(x)|)\le
\|\check{\ell}\|_\infty\sum_{j=0}^n1/2^j$
for $n=1,2,3,....$

Then we have a convergent series
$\sum_{n=0}^\infty \ell_n(t)$ converging locally
uniformly to a function, say
$\phi(t)=\sum_{n=0}^\infty \ell_n(t)$, where
\begin{equation}
\phi(t)=
\left\{\begin{array}{@{\,}lll} \psi(n)e^{-i 2\pi
N_nt}\ell(t-n)&\ \ (t\in [n-1/2, n+1/2], n\in
\mathbf{N}), \\ \psi(0)\ell(t)&\ \ (t\in [-1/2,
1/2]), \\ \phi(n)e^{i2\pi N_nt}\ell(t-n)&\ \ (t\in
[n-1/2, n+1/2], n\in -\mathbf{N}). \end{array}
\right. 
\end{equation}
For
each $t\in \hat{\mathbf{R}}$, $ {\rm
supp}(\ell_n)\cap (t+\bar{V})\neq \emptyset$ at
most one $n\in \mathbf{N}$, so
$\|\phi\|^{\bar{V}}=
\sup_{n\in\mathbf{N}}\|\ell_n\|^{\bar{V}}=\|\ell\|^{\bar{V}}<\infty$,
and we have $\phi\in
C_{BSE}^{\bar{V}}(\hat{\mathbf{R}})=\mathbb{M}({\mathcal 
A}(\hat{\mathbf{R}}))$.

On the other hand, the series $\sum_{n=0}^\infty
\check{\ell_n}(x)$ converges absolutely and
locally uniformly on $\mathbf{R}$
from our constructions of
$\ell_n,\ n=1,2,3,...$; and hence there exists
$\iota F_\phi\in \textfrak{M}(\mathbf{R})$ such
that $\sum_{n=0}^N \frac{1}{2\pi}
\check{\ell_n}(x)dx\rightarrow \iota F_\phi\
(N\rightarrow \infty) $(converges in the vague
topology). It follows that \begin{eqnarray*}
\int_{\mathbf{R}}k(x)d\iota
F_\phi(x)&=&\lim_{N\rightarrow\infty}\int_{\mathbf{R}}k(x)\sum_{n=0}^N\frac{1}{2\pi}\check{\ell}_n(x)dx
= \lim_{N\rightarrow\infty}\int_{\hat{\mathbf{R}}}
\check{k}(t)\sum_{n=0}^N\ell_n(t)dt \\
&=&\int_{\hat{\mathbf{R}}}\check{k}(t)\phi(t)dt\ \
(k\in C_{c,2}(\mathbf{R})) . 
\end{eqnarray*}
Therefore $\iota F_\phi\in
\textfrak{M}_T(\mathbf{R})$ with $d \widehat{\iota
F_\phi}(t)=\phi(t)dt$, that is, $F_\phi\in
\mathbb{M} ({\mathcal L}^1(\mathbf{R}))_{(0)}$ by
Theorem 8.6. 

Note that $\phi(n)=\psi(n)\ \ (n\in
\hat{\mathbf{Z}})$ from (43). By [17, Theorem
2.7.2] and (42), we have 
\begin{eqnarray*}
\sup_{\xi\in
B(\hat{\mathbf{R}})}\|\xi-\phi\|^{\bar{V}}&\ge&
\sup_{\xi\in B(\hat{\mathbf{R}})}\sup_{t\in
\hat{\mathbf{R}}}|\phi(t)-\xi(t)| \ge \sup_{\xi\in
B(\hat{\mathbf{R}})}\sup_{n\in
\hat{\mathbf{Z}}}|\phi(n)-\xi(n)| \\ &=&
\sup_{\zeta\in B(\hat{\mathbf{Z}})}\sup_{n\in
\hat{\mathbf{Z}}}|\psi(n)-\zeta(n)| \ge 1/2,
\end{eqnarray*}
 which implies $\phi\not\in {\mathcal 
B}(\hat{\mathbf{R}})$.

\vspace{1cm}

\Large \begin{center} REFERENCES \end{center}
\normalsize \rm

[1]\ L. Argabright and J. Gil de Lamadrid, Fourier
Analysis of Unbounded Measures on Locally Compact
Abelian Groups, Memoirs of A.M.S. N.145, 1974.

[2]\ K. de Leeuw and Y. Katznelson, The two sides
of a Fourier-Stieltjes transform and almost
idempotent measures, Israel Journal of Math.
8(1970), 213-229.

[3]\ H. G. Feichtinger, On a new Segal algebra,
Monatsh. Math. 92(1981), 269-289.

[4]\ C.C. Graham, The Fourier transform is onto
only when the group is finite, Proc. A.M.S. 38(1973),365-366.

[5]\ C.C. Graham and O. C. McGehee, Essays in
Commutative Harmonic Analysis, Springer-Verlag New
York Inc., 1979.

[6]\ E. Hewitt and K. Ross, Abstract Harmonic
Analysis, Vol.1, Springer Verlag, Berlin 1963.

[7]\ J. Inoue and S.-E. Takahasi, On
characterizations of the image of Gelfand
transform of commutative Banach algebras, Math.
Nachr. 280, No. 1-2,(2007), 105-126.

[8]\ -----, Segal algebras in commutative Banach
algebras, Rocky Mountains J.M., 44-2 (2014),
539-589.

[9]\ -----, Constructions of Segal algebras in
$L^1(G)$ of LCA groups $G$ in which a generalized
Poissons summation formula holds, J. Korean Math.
Soc. 59(2022), No.2, 367-377.

[10]\ M. S. Jakobsen, On a (no longer) new Segal
algebra: a review of the Feichtinger algebra.
(English) J. Fourier Anal. Appl. 24,
No.6, (2018), 1579-1660.

[11]\ R. Larsen, An Introduction to the Theory of
Multipliers, Springer-Verlag Berlin Heidelberg New
York 1971.

[12]\ R. Lyons, Seventy years of Rajchman
measures, The Journal of Fourier Analysis and
Applications, Kahane Special Issue (1995), 
363-377.

[13]\ H. Reiter, Classical harmonic analysis and
locally compact groups, Oxford University Press
1968.

[14]\ -----, $L^1$-Algebras and Segal
Algebras, Lecture notes in Mathematics 231,
Springer-Verlag, Berlin, 1971.

[15]\ -----,  Metaplectic groups and Segal
algebras, Lecture notes Mathematics, Vol.1382.
Springer Verlag, Berlin, 1989.

[16]\ H. Reiter and J.D. Stegeman, Classical
Harmonic Analysis and Locally Compact Groups, 
Oxford Univ. Press Inc., New York 2000.

[17]\ W. Rudin, Fourier Analysis on Groups,
Interscience Publications, Inc., New York, 1962.

[18]\ S.-E. Takahasi and O. Hatori, Commutative
Banach algebras which satisfy a
Bochner-Shoenberg-Eberlein-type theorem
Proc.\,Amer.\,Math. Soc.110 (1990), 149-158.

\vspace{0.3cm}

Jyunji Inoue, Hokkaido University, Sapporo,
060-0808, JAPAN

e-mail address: {\sf
rqstw4a1@nifty.com} \vspace{0.3cm}

Sin-Ei Takahasi, Yamagata University, Yonezawa,
992-8510, Japan

LABORATORY of MATHEMATICS and GAMES, Chiba
273-0032, Japan

https://www.math-game-labo.com

e-mail address: {\sf sin\_ei1@yahoo.co.jp}
\end{document}